\documentclass[12pt]{amsart}
\usepackage{amsthm}
\usepackage{amsmath}
\usepackage{amssymb}
\usepackage{comment}
\usepackage{enumerate}
\usepackage{amsfonts}
\usepackage{dsfont}
\usepackage{mathtools}
\usepackage{afterpage}
\usepackage{hyperref}
\usepackage[capitalize]{cleveref}
\usepackage{todonotes}
\usepackage{algorithm,algorithmic}
\usepackage{multicol}
\usepackage{circuitikz}
\usepackage{tikz-cd}
\usetikzlibrary{math}
\usepackage{bm}
\usetikzlibrary{shapes.geometric}
\usepackage[utf8]{inputenc}
\usepackage[margin=1.1in]{geometry}
\DeclareTextFontCommand{\emph}{\color{blue}\em}

\newtheorem{theorem}{Theorem}[section]
\newtheorem{thm}[theorem]{Theorem}
\newtheorem{def-prop}[theorem]{Definition-Proposition}
\newtheorem{prop}[theorem]{Proposition}
\newtheorem{proposition}[theorem]{Proposition}

\newtheorem{lemma}[theorem]{Lemma}
\newtheorem{cor}[theorem]{Corollary}

\theoremstyle{definition}
\newtheorem{ex}[theorem]{Example}
\newtheorem{example}[theorem]{Example}

\newtheorem{question}[theorem]{Question}
\newtheorem{problem}[theorem]{Problem}
\newtheorem{defin}[theorem]{Definition}
\newtheorem{defn}[theorem]{Definition}
\newtheorem{remark}[theorem]{Remark}

\newtheorem*{claim}{Claim}

\newcommand{\centeredtab}[1]{\begin{tabular}{l} #1 \end{tabular}}

\def\F{\mathcal{F}}
\def\0{{\mathbf{0}}}
\renewcommand{\P}{\mathbf{P}}
\newcommand{\Pbb}{\mathbb{P}}
\newcommand{\Pcal}{\mathcal{P}}
\newcommand{\Ccal}{\mathcal{C}}
\newcommand{\N}{\mathcal{N}}

\newcommand{\e}{\mathbf{e}}
\newcommand{\p}{\mathbf{p}}

\newcommand{\x}{\mathbf{x}}
\newcommand{\y}{\mathbf{y}}

\newcommand{\C}{\mathbf{c}}
\newcommand{\z}{\mathbf{z}}
\newcommand{\Q}{\mathbf{Q}}
\newcommand{\R}{\mathbb{R}}
\newcommand{\PB}{\mathcal{PB}}
\newcommand{\PC}{\mathcal{PC}}
\newcommand{\VolC}{\mathrm{Vol}^{\vee}}
\newcommand{\EVolC}{\mathrm{EVol}^{\vee}}
\DeclareMathOperator{\Supp}{Supp}

\DeclareMathOperator{\Vol}{Vol}

\DeclareMathOperator{\Int}{Int}
\DeclareMathOperator{\id}{id}

\def\F{{\mathcal{F}}}
\def\M{{\mathcal{M}}}

\def\T{{\mathcal{T}}}
\def\R{{\mathbb R}}
\def\A{{\mathcal{A}}}
\def\B{{\mathcal{B}}}
\def\S{{\mathbf{S}}}

\def\Res{{\rm Res}}
\def\sp{{\rm span}}

\def\a{{\mathbf{a}}}
\def\p{{\mathbf{p}}}
\def\q{{\mathbf{q}}}
\def\z{{\mathbf{z}}}
\def\one{{\mathbf{1}}}
\def\v{{\mathbf{v}}}
\def\w{{\mathbf{w}}}

\def\tQ{{\tilde Q}}
\def\J{{\mathcal{J}}}

\def\adj{{\rm adj}}
\def\pol{\vee}
\def\dual{*}

\DeclareMathOperator{\Cone}{cone}

\DeclareMathOperator{\conv}{conv}

\newcommand{\TL}[1]{{\bf TL: #1}}

\newcommand{\Yibo}[1]{\todo[size=\small,inline,color=yellow!30]{#1 \\ \hfill --- Yibo}}

\begin{document}
\title{Dual mixed volume}
\author{Yibo Gao}
\address{\parbox{\linewidth}{Beijing International Center for Mathematical Research, Peking University, \\ \mbox{Beijing, 100871,} China}}
\email{\href{mailto:gaoyibo@bicmr.pku.edu.cn}{{\tt gaoyibo@bicmr.pku.edu.cn}}}
\author{Thomas Lam}
\address{Department of Mathematics, University of Michigan, \mbox{Ann Arbor, MI 48109,} USA}
\email{\href{mailto:tfylam@umich.edu}{{\tt tfylam@umich.edu}}}
\author{Lei Xue}
\address{Department of Mathematics, University of Michigan, \mbox{Ann Arbor, MI 48109,} USA}
\email{\href{mailto:leixue@umich.edu}{{\tt leixue@umich.edu}}}

\thanks{T.L.\ was supported by Grant No.~DMS-1953852 and DMS-2348799 from the National Science Foundation.} 
\date{\today}

\begin{abstract}
We define and study the dual mixed volume rational function of a sequence of polytopes, a dual version of the mixed volume polynomial. This concept has direct relations to the adjoint polynomials and the canonical forms of polytopes. We show that dual mixed volume is additive under mixed subdivisions, and is related by a change of variables to the dual volume of the Cayley polytope.  We study dual mixed volume of zonotopes, generalized permutohedra, and associahedra.  The latter reproduces the planar $\phi^3$-scalar amplitude at tree level. 
\end{abstract}
\maketitle
\setcounter{tocdepth}{1}
\tableofcontents

\section{Introduction}\label{sec:intro}
Let $\Vol(\cdot)$ be the normalized volume function in $\R^d$ where the coordinate simplex has volume $1$. For a sequence of convex bodies $\S =(S_1,S_2,\ldots,S_r)$ in $\R^d$ and for $x_1,\ldots,x_r >0$, we have
\begin{equation}\label{eq:MV}
\Vol(x_1 S_1 + \cdots + x_r S_r) = \sum_{i_1,i_2,\ldots,i_d=1}^r V(S_{i_1},S_{i_2},\ldots,S_{i_d}) x_{i_1} \cdots x_{i_d} =: \Vol_\S(\x),
\end{equation}
where $V(S_{i_1},S_{i_2},\ldots,S_{i_d})$ are the \emph{mixed volumes} of $S_{i_1},\ldots,S_{i_d}$, and we call $\Vol_\S(\x)$ the \emph{mixed volume polynomial}.  The mixed volumes are characterized by the following fundamental properties.
\begin{enumerate}
\item $V(S,S,\ldots,S) = \Vol(S)$,
\item $V$ is symmetric in its arguments,
\item $V$ is multilinear: $V(aS+bS',S_2,\ldots,S_d) = aV(S,S_2,\ldots,S_d)+b V(S',S_2,\ldots,S_d)$.
\end{enumerate}
Mixed volumes are nonnegative and satisfy various inequalities, among them the \emph{Alexandrov-Fenchel inequality}:
\begin{equation}\label{eq:AF}
V(S_1,S_2,S_3,\ldots,S_d)^2 \geq V(S_1,S_1,S_3,\ldots,S_d) \cdot V(S_2,S_2,S_3,\ldots,S_d),
\end{equation}
a deep result in the geometry of convex bodies, with exciting applications in combinatorics \cite{stanley-combinatorial-Aleksandrov-Fenchel}. We refer the reader to \cite{Schneider} for more on mixed volumes.

In this work, we study the \emph{dual mixed volume function}, defined by
$$
m_\S (\x) = m_\S (x_1,\ldots,x_r): = \Vol((x_1S_1 + \cdots + x_r S_r)^\vee) 
$$
where $S^\vee$ denotes the polar of a convex set $S$.  For general convex bodies, $m_\S(\x)$ is a complicated analytic object, but for a sequence of polytopes $\P = (P_1,\ldots,P_r)$, the dual volume function $m_\P(\x)$ is a rational function.  We focus on polytopes in this work, leaving the general dual mixed volume function for future investigation.

Besides the natural parallel with the mixed volume polynomial, the dual mixed volume function is motivated by recent developments in \emph{positive geometry} \cite{ABL,lam2022invitation}, occurring at the interface of combinatorial algebraic geometry and the physics of scattering amplitudes.  The dual mixed volume function specializes to the \emph{canonical form} of a polytope, and produces in a special case the scalar $\phi^3$-amplitude (see \cref{sec:associahedra}), the field theory limit of the open string amplitude.

\medskip 

We now discuss the results of this work in more detail.
\subsection{Dual volumes}
In \cref{sec:dual-volume}, we define the \emph{dual volume} $\VolC(P)$ of a polyhedron $P \subseteq \R^d$.  When $P$ contains the origin $\0 \in \R^d$ in its interior, the dual volume $\VolC(P)$ is simply the volume of the polar of $P$.  However, this function can naturally be extended to all polyhedra $P$ satisfying the non-codegeneracy condition that $\0$ does not belong to the affine span of any of the facets of $P$.  
We also study a rational function $\VolC_\z(P):= \VolC(P - \z)$, called the \emph{dual volume function}.  In \cref{thm: dual vol is integral}, we give an integral formula for $\VolC_{\z}(P)$.

In \cref{section: valuation}, we show that $\VolC(P)$ and $\VolC_\z(P)$ are valuative.  This result generalizes duality results of Filliman~\cite{filliman} and Kuperberg~\cite{Kuperburg}; see \cref{rmk: compare with filliman and kuperburg}.  We deduce the valuative property from the classical result of Lawrence stating that the algebra of indicator functions of cones has a polarity involution.  While our results on dual volumes would not be surprising to experts, we hope that our presentation clarifies the understanding of the dual volume function.

\subsection{Adjoint and canonical form}

In \cref{sec:adjoint}, we show that the numerator of the dual volume function $\Vol_\z(P)$, suitably normalized, coincides with the \emph{adjoint polynomial} of (the cone of) the dual polytope $C(P)^*$, originally defined by \cite{Warren}.

In \cref{sec:canform}, we show that the canonical form $\Omega(P)$ of a polytope $P$, in the sense of positive geometry, is given by $\Omega(P) = \Vol_\z(P) dz_1 \cdots dz_d$. As a consequence, we recover the fact that the adjoint hypersurface is the zero set of the canonical form (\cref{cor:adjointzero}).  We show that the dual volume function behaves well under projective transformations (\cref{cor:proj}).

\subsection{Dual mixed volumes}
We move on to consider a sequence $\P = (P_1,\ldots,P_r)$ of polytopes in $\R^d$.  In \cref{sec:Minkowski} and \cref{sec:DMV}, we define the dual mixed volume function $m_\P(\x)$ of $\P$, show that it is a rational function of degree $-d$, and establish an integral formula.  We study the behavior of the ``polytope" $\x \P := x_1 P_1 + \cdots + x_r P_r$ as a function of $x_1,\ldots,x_r$, even allowing the $x_i$ to be negative.

In Sections~\ref{sec:fine}--\ref{sec:Cayley}, we study the behavior of the dual mixed volume function $m_\P(\x)$ under \emph{mixed subdivisions}: those subdivision of $P_1+\cdots+P_r$ that are compatible with the Minkowski sum.  Our main result (\cref{thm: mixed subdivision formula}) states that dual mixed volume is additive under mixed subdivisions.  Applying this to a fine mixed subdivision, we obtain a decomposition of $m_\P(\x)$ as a sum of monomials (\cref{prop:fine-mixed-dual-volume}). This generalizes the analogous property that canonical forms are additive under subdivisions (\cref{thm:canformsub}).  Mixed subdivisions of $\P$ are in bijection with subdivisions of the corresponding Cayley polytope $\Ccal(\P)$.  We show in \cref{thm:Cayleysub} that the dual volume function of $\Ccal(\P)$ and the mixed volume funciton of $\P$ are related by a change of variables.

\subsection{Zonotopes, generalized permutohedra, and associahedra}
We consider dual mixed volume functions for some important examples.  In \cref{sec:affine}, we develop a formalism to discuss dual volume functions when the affine span of the polytope $P = P_1+ P_2+ \dots +P_r$ belongs to a hyperplane, which is often the case in examples.

In \cref{sec:zonotope}, we consider \emph{zonotopes}, polytopes obtained as Minkowski sums of a number of intervals.  We develop a deletion-contraction formalism for the dual mixed volumes of zonotopes $\P = (P_1, \dots, P_r)$, showing that the rational function $m_\P(\x)$ is uniquely determined dual mixed volumes of the deletion $(P_1,\ldots,P_{r-1})$, and contraction $\P/\P_r$.  By the Bohne-Dress Theorem (see \cref{thm:BD}), mixed subdivisions of zonotopes, or zonotopal tilings, are in bijection with one-element liftings of the corresponding oriented matroid.  Thus we obtain a formula for the dual mixed volume $m_\P(\x)$ for each such lifting.

In \cref{sec:genperm}, we consider \emph{generalized permutohedra}, polytopes obtained by parallel translating the facet hyperplanes of a permutohedron.  We give an explicit formula (\cref{prop:genperm}) for the dual mixed volume of a generalized permutohedron.  The combinatorics of fine mixed subdivisions of generalized permutohedra has been studied in \cite{postnikov-permutohedra}.  Our results lead to some curious identities (\cref{cor:Jm1} and \cref{cor:Jm2}) that appear to be new.

In \cref{sec:associahedra}, we give an explicit formula (\cref{prop:assoc}) for the dual mixed volume of an \emph{associahedron} and relate it to the planar $\phi^3$-amplitude.

\subsection{Further directions}
In this paper, we begun exploring some analogies between mixed volumes and dual mixed volumes.  There are many other natural directions to pursue, and we indicate some in \cref{sec:further}.

\subsection*{Acknowledgements}
We thank Alexander Barvinok for a number of enlightening conversations, especially in relation to \cref{thm:valuation} and \cref{prop:dualBM}.  T.L. was supported by the Simons Foundation and by the National Science Foundation under grants DMS-1953852 and DMS-2348799.

\section{The dual volume function of polyhedra}\label{sec:dual-volume}
\subsection{Cones, fans, and polyhedra}
A \emph{polyhedral cone} $C \subseteq \R^d$ is a non-empty intersection of finitely many closed half spaces, each passing through the origin $\mathbf{0}$.  A polyhedral cone $C$ is called \emph{pointed} if it does not contain any line.  A \emph{face} of $C$ is a subset $F\subseteq C$ such that there exists $\y\in \R^d$ and $b\in \R$ that satisfy 
\begin{eqnarray*}
    \langle \x,\y\rangle &=&b \quad \text{for all } \x\in F, \\
\langle \x,\y\rangle &<&b \quad \text{for all } \x\in C\setminus F. 
\end{eqnarray*}
A \emph{(polyhedral) fan} $\mathcal{F}$ is a finite set of polyhedral cones such that if $C\in\mathcal{F}$ and $F\subseteq C$ is a face, then $F \in\mathcal{F}$, and that if $C_1,C_2\in\mathcal{F}$, then $C_1\cap C_2$ is a face of both $C_1$ and $C_2$.  The support of a fan $\F$ is the union of its cones.  A fan $\F$ in $\R^d$ is called \emph{complete} if the support of $\F$ is the whole space $\R^d$.  We say that a polyhedral fan $\F$ is \emph{pure} of dimension $r$ if all maximal (under inclusion) cones of $\F$ are of dimension $r$.

A \emph{polyhedron} $P \subseteq \R^d$ is a non-empty intersection of finitely many closed half spaces (not necessarily passing through the origin).  A polyhedron $P$ is a \emph{polytope} if it is a bounded subset of $\R^d$.  Faces are defined for polyhedra and polytopes as they are for cones.

For any non-empty closed convex set $S \subseteq\R^d$, the \emph{support function} of $S$ describes the (signed) distances from its supporting hyperplanes to the origin. It is given by
\begin{eqnarray}\label{eq:support}
h_S:\R^d &\rightarrow &\R \cup \{\infty\} \nonumber \\ 
\v &\mapsto & -\min_{\y\in S}\langle\v,\y\rangle.
\end{eqnarray}
The support of this function is given by $\Supp(h_S):= \{\v \in \R^d \mid h_S(\v) \in \R\}$.  The function $h_S$ is a piecewise-linear function on its support.  Since $h_S(\0) = 0$, the support of $h_S$ always includes $\0$ and so is non-empty.


For a polyhedron $P$, the \emph{normal fan} $\N(P)$ consists of the cones 
$$
C_F:=\{\v\in\R^d\:|\: h_P(\v)=-\langle\v,\y\rangle \neq \infty \text{ for every }\y\in F\}
$$ for each face $F$ of $P$. In particular, $h_P$ is a linear function on each $C_F$.  The support of $\N(P)$ is equal to the support of $h_P$.  When $P$ is a polytope, its normal fan $\N(P)$ is a complete fan. 

The support function $h_S$ of a closed convex set $S$ uniquely determines $S$ since we can recover $S$ as
\[ S = \{\y \; |\;\langle \v, \y\rangle \leq -h_P(\v) \, \text{for all } \v\in \R^d \}.  \]

In this section, we develop the notion of the \emph{dual volume} $\VolC(P)$ for polyhedra.  

\subsection{Triangulation of the dual}
\begin{defin}\label{def:universal-dual-volume}
Let $\mathcal{F}$ be a polyhedral fan in $\R^d$ pure of dimension $d$. Each $1$-dimensional face is called a \emph{extreme ray} in $\mathcal{F}$. For each extreme ray, we pick a vector $\mathbf{v}_i$ that spans it ($\mathbf{v}_i$ is called a \emph{generator} of the ray), and collect these vectors as $\mathbf{v}=(\mathbf{v}_1,\ldots,\mathbf{v}_g)$ and associate variables $u_1,\ldots,u_g$ correspondingly. Let $\mathcal{T}=\{C_1,\ldots,C_N\}$ be any triangulation of $\mathcal{F}$ into full-dimensional \emph{simplicial cones}, i.e., cones whose generators are linearly independent. Define the following rational function \[f_{\mathcal{F},\mathbf{v},\mathcal{T}}(u_1,\ldots,u_g):=\sum_{C=\mathrm{span}_{\R_{\geq0}}(\mathbf{v}_{j_1},\ldots,\mathbf{v}_{j_d})\in \mathcal{T}}\frac{|\det(\mathbf{v}_{j_1},\ldots,\mathbf{v}_{j_d})|}{u_{j_1}\cdots u_{j_d}}.\]

If $\F$ is a polyhedral fan in $\R^d$ pure of dimension $r < d$, then we define $$f_{\F,\v,\T}(u_1,\ldots,u_g) := 0.$$
\end{defin}

When $u_1,u_2,\ldots,u_g$ are all positive, $|\det(\mathbf{v}_{j_1},\ldots,\mathbf{v}_{j_d})|/(u_{j_1}\cdots u_{j_d})$ equals the normalized volume of the simplex formed by $\mathbf{0}$, $\mathbf{v}_{j_1}/u_{j_1}$, $\ldots$, $\mathbf{v}_{j_d}/u_{j_d}$.  While the formula for $f_{\mathcal{F},\mathbf{v},\mathcal{T}}$ may initially look daunting, it is simply the volume of some regions, when the $u_i$-s satisfy the positivity condition.
\begin{lemma}\label{lem:universal-dual-volume-independent-of-triangulation}
The rational function $f_{\mathcal{F},\mathbf{v},\mathcal{T}}$ does not depend on the triangulation $\mathcal{T}$. 
\end{lemma}
\begin{proof}
By additivity, it suffices to prove the equality when $\mathcal{F}$ contains a single maximal dimensional cone. Now consider the polytope \[Q=Q(u_1,\ldots,u_g)=\conv\left(\mathbf{0},\frac{\mathbf{v}_1}{u_1},\ldots,\frac{\mathbf{v}_g}{u_g}\right).\]
If all of $\mathbf{v}_j/u_j$'s are vertices of $Q$, then $\Vol(Q)=f_{\mathcal{F},\mathbf{v},\mathcal{T}}$ by construction. When $u_j=|\mathbf{v}_j|$, then all of $\mathbf{v}_j/u_j$'s are vertices of $Q$ since they lie on the unit sphere. There exists $\epsilon>0$ such that whenever $u_j\in(|\mathbf{v}_j|-\epsilon,|\mathbf{v}_j|+\epsilon)$ for all $j$, $Q$ has $\mathbf{v}_j/u_j$'s as its vertices. This means that for two distinct triangulations $\mathcal{T}$ and $\mathcal{T}'$, the two rational functions $f_{\mathcal{F},\mathbf{v},\mathcal{T}}$ and $f_{\mathcal{F},\mathbf{v},\mathcal{T}'}$ agree on a set with positive measure, and therefore, they must be equal. 
\end{proof}

From now on, we write $f_{\mathcal{F},\mathbf{v}}(u_1,\ldots,u_g):=f_{\mathcal{F},\mathbf{v},\mathcal{T}}(u_1,\ldots,u_g)$ for any triangulation $\mathcal{T}$ of $\mathcal{F}$. We also make the following definition for simplicity of notations.
\begin{defin}\label{def:universal-fan-dual-volume}
Let $h:\R^d\rightarrow\R$ be a function such that $h(c\x)=c\cdot h(\x)$ for all $\x\in\R^d$ and $c\in\R$. Let $\F$ be a polyhedral fan. Write \[f_{\F}(h):=f_{\F,\v}(u_1,\ldots,u_g)\] where $\v=(\v_1,\ldots,\v_g)$ are rays in $\F$ and $u_i=h(\v_i)$.
\end{defin}
Since $h(c\x)=c\cdot h(\x)$, scalings of $\v_i$'s do not matter and thus $f_{\F}(h)$ is well-defined. 

We say that $P$ is \emph{non-codegenerate with $\0$}, or simply \emph{non-codegenerate}, if the origin $\0$ is not contained in the affine span of any of the facets of $P$.  

\begin{lemma}
A polyhedron $P$ is non-codegenerate if and only if for any ray $\R_{\geq 0} \cdot \v$ of $\N(P)$, we have $h_P(\v) \neq 0$.
\end{lemma}
\begin{proof}
The origin $\0$ is not contained in the affine span of any facet of $P$. Equivalently, any ray of $\N(P)$ cannot be orthogonal to any vertex of $P$. This is the same as $h_P(\v) \neq 0$ for any ray $\R_{\geq 0}\cdot \v$. 
\end{proof}

The following is the key definition of this section.
\begin{defin}\label{def:dual-volume}
Let $P\subseteq\R^d$ be a full-dimensional non-codegenerate polyhedron.
For each ray in $\N(P)$, pick a vector that spans it, and collect them as $\mathbf{v}=(\mathbf{v}_1,\ldots,\mathbf{v}_g)$. Define the \emph{dual volume} \[\VolC(P):=f_{\N(P)}(h_{P})=f_{\N(P),\mathbf{v}}(h_P(\mathbf{v_1}),\ldots,h_P(\mathbf{v_g})).\] When $P$ is not full dimensional, set $\VolC(P)=0$. Also define the \emph{dual volume function}
 \[\VolC_{\z}(P):=\VolC(P-\z)\]
viewed as a rational function in the coordinates $z_1,\ldots,z_d$ of $\mathbf{z}$.
\end{defin}

The following result is clear from the definitions.
\begin{lemma}
The dual volume $\VolC(P)$ and dual volume function $\VolC_{\z}(P)$ do not depend on the choice of $\mathbf{v}=(\mathbf{v}_1,\ldots,\mathbf{v}_g)$.
\end{lemma}

The name ``dual volume" comes from its connection to the polar dual. Recall 
\begin{defin}\label{def:dual}
For a polyhedron $P\subseteq\R^d$, its \emph{polar dual} is \[P^\pol:=\{\v\in\R^d\:|\: h_P(\v)\leq 1\} = \{\v\in\R^d\:|\: \langle \v,\y \rangle \geq -1 \text{ for all } \y \in P\}. \]
\end{defin}
\begin{prop}\cite[Chapter 4 (1.2)]{Barvinokconvexity}\label{prop:polar-dual-basics} 
If $P$ is a polytope with $\mathbf{0}$ in its interior, then $P^{\pol}$ is also a polytope with $\mathbf{0}$ in its interior. In this case $(P^{\pol})^{\pol}=P$. 
\end{prop}

%
%

\begin{lemma}\label{lem:same-function-interior}
If $P\subseteq\mathbb{R}^d$ is a full dimensional polytope and $\mathbf{0}$ is in its interior, then we have $\VolC(P)=\Vol(P^\pol).$
\end{lemma}
\begin{proof} 
By \Cref{def:dual} and \Cref{prop:polar-dual-basics}, when $\mathbf{0}$ is in the interior of $P$, the vertices of $P^{\pol}$ are precisely $\mathbf{v_j}/h_P(\mathbf{v_j})$ for $j=1,\ldots,g$. Now, $\VolC(P)=$ $f_{\N(P),\mathbf{v}}(h_P(\mathbf{v_1}),\ldots,h_P(\mathbf{v_g}))$ computes $\Vol(P^{\pol})$ by decomposing $P^\pol$ as the union of the cones over its facets, with the cone point at the origin. 
\end{proof}
\Cref{def:dual-volume} is motivated by \Cref{lem:same-function-interior}, and intuitively, one views $\VolC$ as a volume function. To be precise, the notion of the dual volume $\VolC$ is much more powerful, in the sense that $\VolC(P)$ is always defined whenever $\mathbf{0}$ is not contained in the affine span of any facets of $P$, and the rational function $\VolC_{\z}(P)$ is always well-defined. However, the polar dual $P^\pol$ is very often unbounded, including the situation when $P$ is a polytope that does not contain $\mathbf{0}$ in its interior. 
\begin{ex}\label{ex:one-dimensional-dual-volume}
The simplest example is that of a line segment $[a,b]\subset\R^1$. Pick $\v_1=1$ and $\v_2=-1$ so that $h_{P}(\v_1)=-a$ and $h_{P}(\v_2)=b$. Also, $h_{P-z}(\v_1)=-a+z$ and $h_{P-z}(\v_2)=b-z$. This gives \[\VolC_{z}(P)=\frac{1}{z-a}-\frac{1}{z-b}, \qquad \VolC(P) = \VolC_{z}(P)|_{z=0} = \frac{1}{b} - \frac{1}{a} \text{ when $a,b \neq 0$}.\]
\end{ex}

\begin{ex}\label{ex: quadrilateral-dual-volume}
Consider the following polytope $P\subset\R^2$ as the convex hull of $(1,1)$, $(2,1)$, $(3,-1)$, $(1,-1)$, with its normal fan $\N(P)$ shown in Figure~\ref{fig:dual-volume-polytope-example-1}. We can pick $\mathbf{v_1}=(0,1)$, $\mathbf{v_2}=(1,0)$, $\mathbf{v_3}=(0,-1)$ and $\mathbf{v_4}=(-2,-1)$. Summing cyclically with $\mathbf{v_5}=\mathbf{v_1}$, \[f_{\N(P),\mathbf{v}}(u_1,u_2,u_3,u_4)=\sum_{i=1}^4\frac{|\det(\mathbf{v_i},\mathbf{v_{i+1}})|}{u_iu_{i+1}}=\frac{1}{u_1u_2}+\frac{1}{u_2u_3}+\frac{2}{u_3u_4}+\frac{2}{u_4u_1}.\]
Now $\VolC(P)$ is obtained from $f_{\N(P),\mathbf{v}}$ by assigning $u_i=h_P(\mathbf{v_i})$ for all $i$. We have $h_P(\mathbf{v_1})=1$, $h_P(\mathbf{v_2})=-1$, $h_P(\mathbf{v_3})=1$ and $h_P(\mathbf{v_4})=5$. Therefore, \[\VolC(P)=(-1)+(-1)+\frac{2}{5}+\frac{2}{5}=-\frac{6}{5}.\]
We can also calculate $\VolC_{\z}(P)=\VolC(P-\z)$ which is obtained from $f_{\N(P),\mathbf{v}}$ by assigning $u_i=h_{P-\z}(\mathbf{v_i})=h_P(\mathbf{v_i})+\langle \z,\mathbf{v_i}\rangle$. We then have \[\VolC_{\z}(P)=\frac{1}{(1{+}z_2)({-}1{+}z_1)}+\frac{1}{({-}1{+}z_1)(1{-}z_2)}+\frac{2}{(1{-}z_2)(5{-}2z_1{-}z_2)}+\frac{2}{(5{-}2z_1{-}z_2)(1{+}z_2)}.\]
\begin{figure}[h!]
\centering
\begin{tikzpicture}[scale=1.0]
\draw[fill=black,opacity=0.2](1,1)--(2,1)--(3,-1)--(1,-1)--(1,1);
\draw[thick](1,1)--(2,1)--(3,-1)--(1,-1)--(1,1);
\node at (1,1) {$\bullet$};
\node at (2,1) {$\bullet$};
\node at (3,-1) {$\bullet$};
\node at (1,-1) {$\bullet$};
\node[above] at (1,1) {$(1,1)$};
\node[above] at (2.5,1) {$(2,1)$};
\node[below] at (3,-1) {$(3,-1)$};
\node[below] at (1,-1) {$(1,-1)$};
\draw[thin](0,-2)--(0,2);
\draw[thin](-1,0)--(4,0);
\end{tikzpicture}
\qquad\qquad
\begin{tikzpicture}[scale=1.0]
\draw[fill=black,opacity=0.2](0,1.6)--(1.6,0)--(0,-1.6)--(-1.6,-0.8);
\node at (0,0) {$\bullet$};
\draw[->,thin](0,0)--(0,2);
\draw[->,thin](0,0)--(0,-2);
\draw[->,thin](0,0)--(2,0);
\draw[->,thin](0,0)--(-2,-1);
\node[right] at (0,1.7) {$\mathbf{v_1}$};
\node[below] at (1.7,0) {$\mathbf{v_2}$};
\node[left] at (0,-1.7) {$\mathbf{v_3}$};
\node[above] at (-2,-1) {$\mathbf{v_4}$};
\end{tikzpicture}
\caption{A polytope $P$ and its normal fan $\N(P)$}
\label{fig:dual-volume-polytope-example-1}
\end{figure}
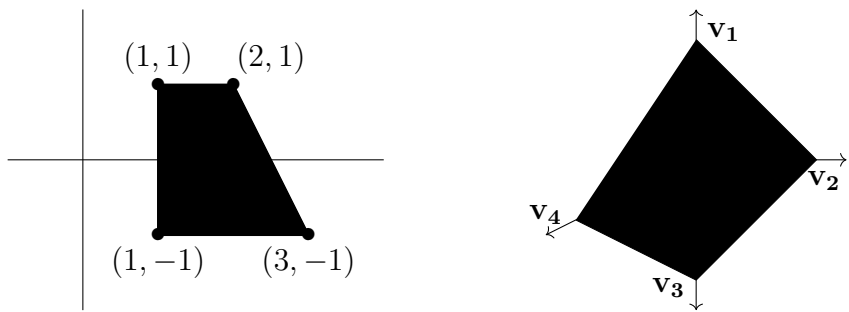
\end{ex}
\begin{ex}
Consider the following unbounded polyhedron $P\subset\R^2$ defined by the inequalities $3x_1+x_2+3\geq0$, $x_1+x_2+1\geq0$, $-2x_1+x_2+4\geq0$ with its normal fan $\N(P)$, which is not a complete fan, shown in Figure~\ref{fig:dual-volume-polyhedron-example-1}. We can pick $\mathbf{v_1}=(-2,1)$, $\mathbf{v_2}=(1,1)$ and $\mathbf{v_3}=(3,1)$. Compute that $h_P(\mathbf{v_1})=4$, $h_P(\mathbf{v_2})=1$, $h_P(\mathbf{v_3})=3$ and \[f_{\N(P),\mathbf{v}}(u_1,u_2,u_3)=\frac{3}{u_1u_2}+\frac{2}{u_2u_3}.\]
Substituting $u_i=h_P(\mathbf{v_i})$ for all $i$ gives $\VolC(P)$ and substituting $u_i=h_P(\mathbf{v_i})+\langle\z,\mathbf{v_i}\rangle$ gives $\VolC_\z(P)$. We therefore obtain $\VolC(P)=17/12$ and \[\VolC_{\z}(P)=\frac{3}{(4-2z_1+z_2)(1+z_1+z_2)}+\frac{2}{(1+z_1+z_2)(3+3z_1+z_2)}.\]
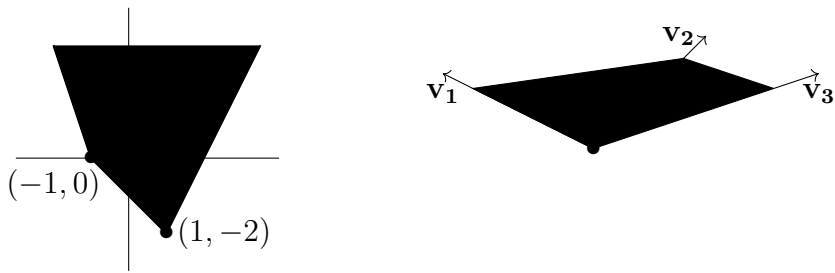
\begin{figure}[h!]
\centering
\begin{tikzpicture}[scale=0.5]
\draw[fill=black,opacity=0.2](-2,3)--(-1,0)--(1,-2)--(3.5,3);
\draw[thick](-2,3)--(-1,0)--(1,-2)--(3.5,3);
\draw[thin](-3,0)--(4,0);
\draw[thin](0,-3)--(0,4);
\node[below] at (-2,0) {$(-1,0)$};
\node[right] at (1,-2) {$(1,-2)$};
\node at (-1,0) {$\bullet$};
\node at (1,-2) {$\bullet$};
\end{tikzpicture}
\qquad\qquad
\begin{tikzpicture}[scale=0.5]
\draw[fill=black,opacity=0.2](0,0)--(-3.2,1.6)--(2.4,2.4)--(4.8,1.6);
\node at (0,0) {$\bullet$};
\draw[->,thin](0,0)--(-4,2);
\draw[->,thin](0,0)--(3,3);
\draw[->,thin](0,0)--(6,2);
\node[below] at (-4,2) {$\mathbf{v_1}$};
\node[left] at (3,3) {$\mathbf{v_2}$};
\node[below] at (6,2) {$\mathbf{v_3}$};
\node at (0,-3) {};
\node at (0,4) {};
\end{tikzpicture}
\caption{An unbounded polyhedron $P$ and its normal fan $\N(P)$}
\label{fig:dual-volume-polyhedron-example-1}
\end{figure}
\end{ex}

\subsection{An integral formula}
\begin{theorem}\label{thm: dual vol is integral}
Let $P\subset\R^d$ be a polyhedron and $\z\in \R^d$ be any point such that $P-\z$ is non-codegenerate. Then \begin{align}\label{eqn: dual volume = exp integral}
\VolC_{\z}(P)=\int_{\R^d}\exp(-h_P(\v) - \langle\v,\z\rangle)d\v.
\end{align}
\end{theorem}

We begin with the following basic formula of exponential integrals over a simplicial cone.  For a cone $C \subseteq\R^d$, define the \emph{dual cone} $C^{\dual}$ by
 \begin{equation}
 \label{eq:dualcone}
 C^\dual := \{\v\in \R^d\; | \; \langle \v,\y \rangle \geq 0\; \text{ for all } \y\in C \}. 
 \end{equation}
Note that we use ${}^*$ for polarity/duality on cones and ${}^\vee$ for polarity on polytopes.

\begin{prop}\label{prop: integral formula for simplex}
 Let $C$ be a full dimensional simplicial cone in $\R^d$ spanned by the linearly independent vectors ${\v}_1, {\v}_2, \dots, {\v}_d$. 
 Then for all $\C$ in the the interior of $C^*$, we have
 \[\int_C e^{-\langle \v, \C \rangle} d\v = \frac{|\det (\v_1, \v_2, \dots, \v_d)|}{\prod\limits_{i=1}^d \langle \v_i,-\C \rangle}. \]
\end{prop}
\begin{proof}
Changing integration variables to $\alpha_1,\ldots,\alpha_d$ given by $\v = \alpha_1 \v_1 + \cdots + \alpha_d \v_d$, we have
 \begin{align*}\int_C e^{-\langle \v, \C \rangle} d\v &= \int_{\R_{\geq 0}^d} e^{-\sum_i \alpha_i \langle v_i, \C \rangle} |\det (\v_1, \v_2, \dots, \v_d)| d\alpha_1 \cdots d\alpha_d \\
 &= |\det (\v_1, \v_2, \dots, \v_d)| \prod_{i=1}^d \int_0^\infty e^{-\alpha_i \langle v_i, \C \rangle}d\alpha_i.
 \end{align*}
 The result then follows from the equality $\int_0^\infty e^{-\alpha z }d\alpha = 1/z$ for any $z > 0$.
\end{proof}
\begin{proof}[Proof of Theorem \ref{thm: dual vol is integral}]

 Let $\mathcal{N}(P)$ be the normal fan of the polyhedron $P$, then $\mathcal{N}(P)$ is the collection of normal cones associated with the vertices of $P$, i.e.,
 \[\mathcal{N}(P) = \bigcup_{\p \in V(P)} C_\p, \]
 where $C_\p := \{\v\in \R^d \; |\; h_P(\v) = -\langle \v, \p\rangle \}$ and $V(P)$ is the vertex set of P.
 Now consider $P-\z$. We have $\mathcal{N}(P-\z)=\mathcal{N}(P)$ since translation of a polyhedron preserves the normal fan. Its support function is
 \[h_{P-\z}(\v) = h_P(\v) + \langle \v, \z\rangle.\]
For and arbitrary point $\x\in \R^d$, 
\begin{align*}
h_{P-\z}(\v) = \begin{cases}
-\langle \v, \p\rangle + \langle \v, \z\rangle = \langle \v, \z -\p \rangle & \text{if }\v\in C_\p \text{ for some } \p\in V(P); \\
\infty & \text{ otherwise.} \\
\end{cases}
\end{align*}
Therefore 
\begin{eqnarray*}
\int_{\R^d}e^{-h_P(\v)-\langle\v,\z\rangle}d\v 
&=& \int_{\R^d}e^{-h_{P-\z}(\v)}d\v \\
&=& \sum_{\p \in V(P)}\int_{C_\p}e^{-h_{P-\z}(\v)}d\v \\
&=& \sum_{\p \in V(P)}\int_{C_\p}e^{-\langle \v, \z -\p \rangle}d\v.
\end{eqnarray*}

Next, we will subdivide $\mathcal{N}(P)$ into simplicial cones to apply Proposition \ref{prop: integral formula for simplex}. Let $\mathcal{T}$ be any such subdivision, by Proposition \ref{prop: integral formula for simplex} we have

\begin{eqnarray*} 
\sum_{\p \in V(P)}\int_{C_\p}e^{-\langle \v, \z -\p \rangle}d\v 
&=& \sum_{\p \in V(P)} \sum_{\substack{K\in \mathcal{T}, C\subseteq C_\p \\ C =  \mathrm{span}_{\R_{\geq0}}\{ \v_1, \dots, \v_d\}} } \int_{C}e^{-\langle \v, \z -\p \rangle}d\v \\
&=& \sum_{\p \in V(P)} \sum_{\substack{C\in \mathcal{T}, C\subseteq C_\p \\ C =  \mathrm{span}_{\R_{\geq0}} \{ \v_1, \dots, \v_d\}} } 
\frac{|\det (\v_1, \v_2, \dots, \v_d)|}{\prod\limits_{i=1}^d \langle \v_i, \z-\p\rangle} \\
&=& \sum_{\substack{C\in \mathcal{T} \\ C =  \mathrm{span}_{\R_{\geq 0}} \{ \v_1, \dots, \v_d\}} } 
\frac{|\det (\v_1, \v_2, \dots, \v_d)|}{\prod\limits_{i=1}^d h_{P-\z}(\v_i)}. 
\end{eqnarray*}
Comparing the right-hand side with the formula of $\VolC (P-\z)$ in Definition \ref{def:dual-volume} yields that $\VolC_{\z} (P)= \VolC (P-\z) =  \int_{\R^d}e^{-h_{P-\z}(\v)}d\v$.
\end{proof}

If $P$ is a non-codegenerate polyhedron, the following corollary holds as we take the limit on both sides of (\ref{eqn: dual volume = exp integral}) by letting $\z$ go to the origin.

\begin{cor}\label{cor:intform}
Let $P\subset\R^d$ be a non-codegenerate polyhedron. Then \begin{align}
\VolC(P)=\int_{\R^d}e^{-h_P(\v)}d\v.
\end{align}
\end{cor}

\begin{ex}
Consider the polytope $P = [1,3]$ in $\R$. Its normal fan $\mathcal{N}(P)$ consists of two $1$-dimensional cones (rays) $C_1$ spanned by $v_1 =1$ and $C_2$ spanned by $v_2=-1$. For any $z\in \R^d$, the support function $h_{P-z}$ is
\[h_{P-z} (v) = \begin{cases}
-(1-z)v & v\in C_1,\\
-(3-z)v & v\in C_2.
\end{cases}
\]
We can evaluate the integral as
\begin{eqnarray*}  \int_{\R}e^{-h_{P-\z}(v)}dv &=& \int_0^\infty e^{(1-z)v }dv + \int_{-\infty}^0 e^{(3-z)v}dv\\
&=& \frac{1}{h_{P-z} (v_1)} + \frac{1}{h_{P-z}(v_2)}\\
&=& \frac{1}{z-1} + \frac{1}{3-z}.
\end{eqnarray*}
The dual volume of $P$ is therefore
\[ \VolC(P) = \lim_{z\to 0} \big(\frac{1}{z-1} + \frac{1}{3-z} \big) = -1 +\frac{1}{3} = - \frac{2}{3}.\]
\end{ex}

\section{Valuative property of dual volume}\label{section: valuation}

We have defined $\VolC(P)$ for non-codegenerate polyhedra and $\VolC_\z(P)$ for polyhedra $P \subseteq\R^d$.  In this section, we show that $\VolC$ and $\VolC_\z$ are valuations.  

For any subset $S\subseteq \R^d$, the \emph{indicator function} from $\R^d$ to $\{0,1\}$ is defined as follows.
\begin{eqnarray*}
    [S](\x) = 
    \begin{cases}
    1, & \x\in S, \\
    0, & \x\notin S.
    \end{cases}
\end{eqnarray*}
Let $[P]$ be the indicator function of a polyhedron $P$.  Let $\Pcal$ denote the space of functions spanned by indicator functions of polyhedra.


The following result is also known to Alexander Barvinok \cite{Barvinokprivatecommunication}.

\begin{theorem}\label{thm:valuation}
Suppose that we have
\begin{equation}\label{eq:sumP}
\sum_{i=1}^r \alpha_i [P_i] = 0
\end{equation}
in $\Pcal$, for polyhedra $P_1,\ldots,P_r\in\R^d$ and $\alpha_1,\ldots,\alpha_r \in \R$. Then
\begin{equation}\label{eq:valz}
\sum_{i=1}^r \alpha_i \Vol_\z^\vee(P_i) = 0.
\end{equation}
If each $P_i$ is non-codegenerate, then
\begin{equation}\label{eq:noncodeg}
\sum_{i=1}^r \alpha_i \Vol^\vee(P_i) = 0.
\end{equation}
\end{theorem}

\begin{remark}\label{rmk: compare with filliman and kuperburg}
The case of Theorem \ref{thm:valuation} where the $P_i$'s are all simplices containing the origin was established by Filliman~\cite{filliman} in 1990, and is known as Filliman's duality. In 2003, Kuperburg \cite{Kuperburg} extends Filliman's duality to an involution on the space of polytopal measures of non-codegenerate polytopes, and we have borrowed the terminology ``non-codegenerate" from Kuperberg. Our results extends both results to include all (possibly unbounded) polyhedra; see also discussion at the end of \cite{Kuperburg}.  \cref{thm:valuation} is also known to Alexander Barvinok~\cite{Barvinokprivatecommunication}. 
\end{remark}

We prove \cref{thm:valuation} using the polytope algebra.
\begin{lemma}\label{lem:Vol0}
If $P$ is not full-dimensional or $P$ contains a line, then we have $\Vol^\vee_\z(P) = 0$ and $\Vol^\vee(P) = 0$. 
\end{lemma}
\begin{proof}
When $P$ is not full-dimensional this is part of the definition.  When $P$ contains a line, then $\N(P)$ only consists of lower-dimensional cones, so it follows from \cref{def:universal-dual-volume} that $\VolC(P) = \VolC_\z(P) = 0$.
\end{proof}

Lemma \ref{lem:Vol0} implies that in \eqref{eq:sumP}, lower-dimensional polytopes can be ignored and the conclusion of \cref{thm:valuation} is still valid.

We deduce \cref{thm:valuation} from the polarity map of the cone algebra.  

Let $\Ccal$ denote the space of functions spanned by indicator functions $[C]$ of cones in $\R^d$.  For a cone $C \subseteq\R^{d+1}$, the dual cone $C^*$ is defined in \eqref{eq:dualcone}.
If $C$ is a pointed polyhedral cone of full-dimension then so is $C^*$.  If $C$ is not full-dimensional (resp. $C$ contains a line) then $C^*$ contains a line (resp. $C^*$ is not full-dimensional).

The following result is a consequence of Lawrence's duality theorem.
\begin{theorem}[{\cite[Theorem IV.1.5 and Problem 3 on p.149]{Barvinokconvexity}}] \label{thm:Bar}
There is a polarity involution $\chi \mapsto \chi^*$ on $\Ccal$ that sends $[C]$ to $[C^*]$ for every cone $C \subseteq\R^{d+1}$.
\end{theorem}

Identify $\R^d$ with the affine hyperplane $x_0 = 1$ in $\R^{d+1} = \{(y_0,y_1,\ldots,y_d)\}$.  
For a polyhedron $P \subseteq\R^d$, define the cone over $P$ to be
$$
C(P) := \overline{\{(t,t\y) \mid t \in \R_{\geq 0} \text{ and } \y \in P\}}.
$$
If $P \subseteq\R^d$ is a polytope, then the closure operation is not necessary.  If $P$ is unbounded and contains a ray $\R_{\geq 0} \cdot \y$, then $C(P)$ contains the ray $\R_{\geq 0 } \cdot (0,\y)$.  If $P$ contains a line $\R \cdot \y$, then $C(P)$ contains the line $\R \cdot (0, \y)$.  

\begin{lemma}\label{lem:Cval}
There is a linear map $C(-): \Pcal \to \Ccal$ sending $[P]$ to $[C(P)]$ for every polyhedron $P \subseteq\R^d$.
\end{lemma}
\begin{proof}
Recall that the Euler characteristic $\mu: \Pcal \to \R$ (\cite{mcmullen-polytope-algebra}) is a linear map satisfying $\mu([P]) = 1$ for any non-empty closed polyhedron $P$.
For $\chi \in \Pcal$, define $C(\chi) \in \Ccal$ by $C(\chi)(t,t\x) = \chi(\x)$ for $t \neq 0$, and $C(\chi)({\bf 0}) = \mu(\chi)$.  It is clear that this has the required properties.
\end{proof}

\begin{lemma}
There is a linear map $P(-): \Ccal \to \Pcal$ sending $[C]$ to $C \cap \R^d$ for every cone $C$.
\end{lemma}
\begin{proof}
    For $\chi \in \Ccal$, define $P(\chi) \in \Pcal$ by $P(\chi)(\x) = \chi(1,\x)$.
\end{proof}

\begin{lemma}
We have $P \circ C = \id$. 
\end{lemma}
\begin{proof}
    Follows immediately from the definitions of $P$ and $C$.
\end{proof}

Define a polarity map $*: \Pcal\to \Pcal$ by 
$$
\chi \mapsto C(\chi) \mapsto C(\chi)^* \mapsto P(C(\chi)^*).
$$

\begin{lemma}\label{lem:Pdual}
 If $P \subseteq\R^d$ is a polyhedron containing $\0$ in its interior, then $P^\pol= C(P)^* \cap \R^d$.  Thus $[P]^* = [P^\pol]$.
\end{lemma}
\begin{proof}
Suppose that $\x \in P^\pol$.  Then $\langle \x, \y \rangle \geq -1$ for $\y \in P$.  For $t,s >0$, we have
$$
\langle (t,t\x), (s,s\y) \rangle = ts + ts \langle \x, \y \rangle \geq ts - ts = 0 \qquad \text{ for all } \y \in P.
$$
The inequality still holds after taking closure, so $(t,t\x) \in C(P)^\dual$, giving $C(P^\pol) \subseteq C(P)^\dual$, and similarly we deduce that $C(P)^\dual \subseteq C(P^\pol)$.  It follows that $C(P)^\dual = C(P^\pol)$ and thus $P^\pol= C(P)^\dual \cap \R^d$, and the result follows.
\end{proof}

Let $\Ccal^\ell \subseteq \Ccal$ be the subspace spanned by indicator functions of cones that contain lines.  Let $\Ccal^s \subseteq \Ccal$ be the subspace spanned by indicator functions of cones that are not full-dimensional.  Let $\Ccal^u\subseteq \Ccal$ be the subspace spanned by indicator functions of cones $C(P)$ where $P \subseteq \R^d$ is a polytope.

\begin{lemma} \label{lem:lsu}
Any function in $(\Ccal^\ell + \Ccal^s) \cap \Ccal^u$ belongs to $\Ccal^s$.  Equivalently, $(\Ccal^\ell + \Ccal^s)/\Ccal^s \cap (\Ccal^u + \Ccal^s)/\Ccal^s = 0$.
\end{lemma}
\begin{proof}
We will work modulo $\Ccal^s$ throughout, or equivalently, we can work with measures on $\R^{d+1}$ instead of functions.

Let $\chi \in \Ccal/\Ccal^s$ and $H \subseteq \R^{d+1}$ be a hyperplane through the origin with open subspaces $H_{>0}$ and $H_{<0}$ on either side of $H$.  We define the \emph{gradient} $\partial_H \chi$ of $\chi$ along $H$ to be the function on $H$ given by
$$
(\partial_H \chi)(\z) = \chi(U_\z \cap H_{>0}) - \chi(U_\z \cap H_{<0}) 
$$
where $\z \in H$ and $U_\z$ is a sufficiently small neighborhood of $\z$ in $\R^{d+1}$.  For a dense subset of $H$, the function $\chi$ is constant on $U_\z \cap H_{>0}$ and $U_\z \cap H_{<0}$.  We obtain a function $\partial_H \chi \in \Ccal(H)$, where $\Ccal(H)$ denotes the cone algebra in $H$.  For any expression $\chi = \sum_i \alpha_i [C_i]$ for $\chi$ in terms of indicator functions of full-dimensional cones, we have
$$
\partial_H \chi (\z) = \sum_{H \text{ is a facet of } C_i}  \pm \alpha_i [C_i \cap H]
$$
where the sign is $+$ (resp. $-$) if $C_i$ belongs to $H_{>0}$ (resp. $H_{<0}$).  Note that if $H$ is a facet hyperplane of a cone $C$, then any line contained in $C$ is contained in $H$.  In particular, if $\chi \in \Ccal^\ell/\Ccal^s$ then we must have $\partial_H \chi \in \Ccal^\ell/\Ccal^s$.

Now let $\chi \in (\Ccal^\ell + \Ccal^s)/\Ccal^s \cap (\Ccal^u + \Ccal^s)/\Ccal^s$.  We suppose that we have $\chi = \sum_i \alpha_i [C(P_i)] \neq 0$ where $P_i \subseteq \R^d$ are polytopes of full-dimension in $\R^d$.  Since the summation is finite, there is a hyperplane $H$ such that $\partial_H \chi \neq 0$, and since we must have that $H$ intersects $\R^d$, we have $\partial_H \chi \in (\Ccal^\ell(H) + \Ccal^s(H))/\Ccal^s(H) \cap (\Ccal^u(H) + \Ccal^s(H))/\Ccal^s(H)$.  By repeating this, we reduce to the case that $d = 0$.  In this case, the only possible full-dimensional (that is, one-dimensional cones in $\R^1$) cones are the positive ray $\R_{\geq 0}$, the negative ray $\R_{\leq 0}$, and the whole line $\R$.  We have that $(\Ccal^\ell(\R) + \Ccal^s(\R))/\Ccal^s(\R)$ is spanned by $[\R]$ and $(\Ccal^u(\R) + \Ccal^s(\R))/\Ccal^s(\R)$ is spanned by $[\R_{\geq 0}]$, and we arrive at a contradiction.
%
%
%
\end{proof}

Let $\Ccal^{w}$ be the subspace of functions supported on the (weak) upper halfspace; that is, those functions $\chi$ satisfying $\chi((t,\y)) \neq 0$ only if $t \geq 0$.  It is clear that $\Ccal^u \subseteq \Ccal^w$.

\begin{lemma}\label{lem:wl}
We have $\Ccal^w + \Ccal^\ell = \Ccal$.
\end{lemma}
\begin{proof}
By triangulation, it suffices to show that $[C] \in \Ccal^w + \Ccal^\ell$ for any simplicial cone $C$.  Let $C$ have generators $r_1,\ldots,r_a$.  Suppose that $r_i = (t_i,r'_i)$.  If $t_i \geq 0$ for all $i$, then $[C] \in \Ccal^w$.  Otherwise, suppose that $t_i < 0$ and consider the simplicial cone $C'$ with generators $r_1,\ldots, -r_i, \ldots, r_a$.  Then we have $[C] + [C'] -[D] \in \Ccal^\ell$ where $D$ is a lower-dimensional cone.  By induction on dimension, we may suppose that $[D] \in \Ccal^w + \Ccal^\ell$.  The simplicial cone $C'$ has more generators in the upper halfspace than $C$.  By induction on the number of such generators, we conclude that $[C] \in \Ccal^w + \Ccal^\ell$.
\end{proof}

\begin{lemma}\label{lem:codeg} If $P \subseteq \R^d$ is non-codegenerate then $[C(P)^*] \in \Ccal^u + \Ccal^\ell + \Ccal^s$.
\end{lemma}
\begin{proof}
$C(P)$ is codegenerate if and only if $C(P)^*$ has a generating ray that lies on the hyperplane $H_0 = \{(0,\y)\} \subseteq \R^{d+1}$.  If $C(P)^*$ has no such rays, then by the argument in the proof of \cref{lem:wl}, modulo $\Ccal^\ell$ we may write $[C(P)^*]$ as an element of $\Ccal^u + \Ccal^s$.
\end{proof}
\begin{proof}[Proof of \cref{thm:valuation}]
Define a linear function $\Vol: (\Ccal^u + \Ccal^\ell + \Ccal^s)/(\Ccal^\ell + \Ccal^s) \to \R,$ by
$$
\Vol([C(P)]) = \Vol(P) \mbox{ for any polytope $P \subseteq \R^d$}.
$$
By \cref{lem:lsu}, and the fact that $\Vol$ is a valuation for polytopes, this map is well-defined.  Note that if $\chi$ has lower-dimensional support than $\Vol(\chi) = 0$.

We claim that $\Vol^\vee(P) = \Vol([C(P)]^*)$ when $P$ is non-codegenerate.  To see this, consider the rational functions $\VolC_\z(P)$ and $\Vol([C(P-\z)]^*)$.  It follows from \cref{lem:Pdual} and \cref{lem:codeg} that these rational functions agree when $\z$ is in the interior of $P$, and thus $\VolC_\z(P) = \Vol([C(P-\z)]^*)$ as rational functions.  Substituting $\z = 0$, we obtain the desired identity.

The valuative property of $\Vol^\vee$ then follows from the valuative property of taking cones (\cref{lem:Cval}), the valuative property for the duality map of $\Ccal$ (\cref{thm:Bar}), and the valuative property of $\Vol$ on $(\Ccal^u + \Ccal^\ell + \Ccal^s)/(\Ccal^\ell + \Ccal^s)$.

Note that if $P_i$ is not full-dimensional, then $C(P_i)$ is not full-dimensional, and $C(P_i)^*$ contains a line.  If $P_i$ contains a line, then so does $C(P_i)$, and thus $C(P_i)^*$ is lower-dimensional.  In either case we have $[C(P_i)^*] \in \Ccal^\ell + \Ccal^s$, and thus
$$
\Vol^\vee(P_i) = 0, \qquad \text{agreeing with} \qquad \Vol([C(P_i)^*] )=0.
$$
The equality \eqref{eq:noncodeg} follows immediately, while the equality \eqref{eq:valz} follows since the left hand side is a rational function in $\z$, and for generic $\z$, all the polyhedra $P_i - \z$ are non-codegenerate.
%
%
%
%
%
\end{proof}

\section{Adjoint polynomial}\label{sec:adjoint}

In this section, we discuss the relation of $\Vol^\vee_\z(P)$ to the \emph{adjoint} of a polytope.

Let $C$ be a polyhedral cone in $\R^{d+1}$ generated by a set of extreme rays $V(C)$, i.e., each ray $\v\in V(C)$ cannot be written as a convex hull of the rest. A \emph{triangulation} of $C$ is a collection $T$ of simplicial cones of the same dimension satisfying the following conditions:
\begin{itemize}
\item The union of all simplicial cones in $T$ is $C$.
\item The intersection of two cones is a face in both.
\item The extreme rays of each cone are within $V(C)$.
\end{itemize}

The adjoint polynomial of a polytope was first introduced by Warren in \cite{Warren}.  We use the version given in Aluffi \cite{Aluffi}; see also \cite{Kohn-Ranestad}.

\begin{defin}\cite[Definition/Theorem 4.1]{Aluffi}\label{def-adjoint}
Let $C$ be a polyhedral cone in $\R^{d+1}$ generated by the extreme rays $V(C)$, and let $T$ be a triangulation of $C$. The \emph{adjoint polynomial} of $C$ is given by
\begin{equation}\label{eqn: def-adjoint}
\adj_C(\z) = \sum_{F \in T} |\det(F)| \prod_{\v \in V(C) \setminus V(F)} \langle \v, \z\rangle
\end{equation}
where $F$ varies over simplicial cones $F = \sp_{\R_{\geq 0}}(\v_1,\ldots,\v_{d+1})$ in $T$, and $\det(F) = \det(\v_1,\ldots,\v_{d+1})$.
\end{defin}

This definition is independent of the choice of the triangulation $T$.

In \cite{Warren}, Warren gave essentially the same definition of the adjoint but required all the rays in $V(C)$ to be normalized to the length $1$. Since in each summand of (\ref{eqn: def-adjoint}), each extreme ray in $V(C)$ shows up exactly once, we have that Warren's adjoint polynomial differs from our version by the normalizing factor $\frac{1}{\mu}$ with $\mu = \prod_{\v\in V(C)} |\v|$.



By \cref{lem:Pdual}, if $P$ contains the origin in its interior, the intersection of $C(P)^*$ and the hyperplane $x_1=1$ is $P^\vee$.

We will compare the adjoint polynomial $\adj_{C(P)^*}$ with the dual volume function $\VolC_{\z}(P)$. We start with the following useful property of the adjoint function.

\begin{proposition}\label{prop: adj shifting}
    For any vector $\z \in \R^d$, let $\bar{\z} := (1, \z) \in \R^{d+1}$. Then for all $\z, \z' \in \R^d$, 
    \begin{equation}\label{eqn: adj shifting}
        \adj_{C(P)^*}(\overline{\z+\z'}) = \adj_{C(P-\z')^*}(\bar{\z}).
    \end{equation}
\end{proposition}

\begin{proof} 
    We let $C:= C(P)^*$. Each extreme ray of $C$ is of the form $\w_i =(h_P(\v_i),\v_i)$ for some $\v_i\in N(P)$. Let $C':= C(P-\z')^*$.  Since translating a polytope does not change its normal fan, the extreme rays of $C'$ are given by
    \[ \w_i' = (h_{P-\z'}(\v_i), \v_i) = (h_P(\v_i) +\langle \v_i, \z'\rangle ,\v_i) = \w_i + \langle \v_i, \z'\rangle \e_1. \]
   Since every ray $\w' \in C'$ is a ray $\w$ from $C$ shifted only in the first coordinate by a linear term, every triangulation $T$ of $C$ induces a triangulation $T'$ of $C'$. For each simplicial cone $F = \Cone \{ \w_1, \w_2, \dots, \w_{d+1} \} \in T$, we have
   \begin{eqnarray*}
    \det (F') 
    &=& \det(\w_1', \w_2', \dots, \w_{d+1}')\\
    &=& \det(\w_1, \w_2, \dots, \w_{d+1})
        + \underbrace{\begin{vmatrix}
            \langle \v_1, \z'\rangle  & \langle \v_2, \z'\rangle & \dots & \langle \v_{d+1}, \z'\rangle\\
            \vert & \vert & \dots & \vert\\
            \v_1   & \v_2 & \dots & \v_{d+1}  \\
            \vert & \vert & \dots & \vert
        \end{vmatrix}}_{= 0}\\
    &=& \det(F)
   \end{eqnarray*}
   Also notice that 
    \begin{equation*}
        \langle \w_i', \bar{\z} \rangle = \langle \w_i, \bar{\z} \rangle + \langle \v_i, \bar{\z'} \rangle = \langle \w_i, \overline{\z+\z'} \rangle.
    \end{equation*}
Together, we obtain
    \begin{eqnarray*}
    \adj_{C'}(\bar{\z}) 
    &=& 
    \sum_{F'\in T'} |\det (F')| \prod_{\w' \in V(C')\setminus V(F')} \langle \w', \bar{\z} \rangle\\
    &=&
    \sum_{F\in T} |\det (F)| \prod_{\w \in V(C)\setminus V(F)} \langle \w_i, \overline{\z+\z'} \rangle
    = 
    \adj_{C}(\overline{\z+\z'}) 
    \end{eqnarray*}
as desired.
\end{proof}

To state the relation between $\adj_{C(P)^*}$ and the dual volume function $\VolC_\z(P)$, we first combine the summands in the definition of $\VolC_{\z}(P)$ into a single fraction.

\begin{align}\label{eqn: VolC one fraction}
\begin{split}
 \VolC_\z(P) &=f_{\N(P),\mathbf{v}}(h_{P-\z}(\mathbf{v_1}),\ldots,h_{P-\z} (\mathbf{v_g})) 
 \\
 &= \frac{\displaystyle\sum\limits_{C\in \mathcal{T}}|\det(C)| \displaystyle \prod\limits_{\v\in V(\N(P))\setminus V(C)} h_{P-\z}(\v)}{\displaystyle\prod_{\v\in V(\N(P))} h_{P-\z}(\v)}.
 \end{split}
 \end{align}
where $\mathcal{T}$ is any triangulation of $\N(P)$, $C$ is a simplicial cone in $\mathcal{T}$ generated by its extreme rays $V(C) = \{\mathbf{v}_{1},\ldots,\mathbf{v}_{d}\}$, and $\det(C) = \det(\v_1 , \dots, \v_d)$.

Define 
\begin{equation} \label{eq:denominator of VolC}
\begin{split}
A_{\z}(P) &:= \displaystyle\sum\limits_{C\in \mathcal{T}}|\det(C)| \displaystyle \prod\limits_{\v\in V(\N(P))\setminus V(C)} h_{P-\z}(\v),   \\
B_{\z}(P) &:=  \displaystyle\prod_{\v\in V(\N(P))} h_{P-\z}(\v).
\end{split}
 \end{equation}
For the rest of the paper we will call $A_{\z}(P)$ (resp. $B_{\z}(P)$) the \emph{numerator} (resp. \emph{denominator}) of $\VolC_{\z}(P)$.

 
\begin{theorem}\label{thm:adjoint} For any vector $\z \in \R^d$ we let $\bar{\z} = (1, \z) \in \R^{d+1}$. Then we have
\begin{equation*} \VolC_{\z} (P) =\frac{\adj_{C(P)^*}(\bar{\z} ) }{B_{\z}(P)}, \qquad \text{or equivalently,} \qquad  \adj_{C(P)^*}(\bar{\z}) = A_{\z}(P).
\end{equation*}
\end{theorem}

\begin{proof}
Shifting by an arbitrary vector $\z'\in \R^d$, the statement is equivalent to
\[ \adj_{C(P)^*}(\overline{\z+\z'}) = B_{\z+\z'}(P) \cdot \VolC_{\z+\z'}(P) \quad \text{ for all } \z\in \R^d.\]

By definition, $\VolC_{\z+\z'}(P) = \VolC_{\z}(P-\z')$ and $B_{\z+\z'}(P) = B_{\z}(P-\z')$. On the other hand, in \cref{prop: adj shifting} we get $\adj_{C(P)^*}(\overline{\z+\z'}) = \adj_{C(P-\z')^*}(\bar{\z})$. Therefore the statement above is equivalent to 
\[ \adj_{C(P-\z')^*}(\bar{\z}) = B_{\z}(P-\z') \cdot \VolC_{\z}(P-\z') \quad \text{ for all } \z\in \R^d.\]
Hence it suffices to assume that $P$ (in $\R^d$) contains the origin in its interior. By \cref{lem:Pdual}, in this case $C:= C(P)^* = C(P^\vee)$. 

As before, each extreme ray of $C$ is of the form $\w = (h_P(\v), \v)$, where $\v \in \mathcal{N}(P)$, and so 
\[\langle \w, \bar{\z} \rangle = h_P(\v) + \langle \v, \z\rangle = h_{P-\z}(\v).\]
Therefore
\begin{equation}\label{eqn: adj/B_z}
\frac{\adj_C(\bar{\z})}{B_\z(P)}
=\frac{\sum_{F\in T} |\det(F)| \displaystyle\prod_{\w \in V(C) \setminus V(F)} h_{P-\z}(\v)}{\displaystyle\prod_{\v\in \mathcal{N}(P)} h_{P-\z}(\v)} =
\displaystyle\sum_{F\in T} \frac{|\det(F)| }{\displaystyle\prod_{\w \in V(F)} h_{P-\z}(\v)}.
\end{equation}

Intersecting every simplicial cone in $T$ with the $\{x_1=1\}$ hyperplane gives us a triangulation $T'$ of $P^\vee$. Since $P$ and $P-\z$ have the same normal fan, the vertices of $(P-\z)^\vee$ are just scaling of the vertices of $P^\vee$, and so $T$ also gives rise to a triangulation $T''$ of $(P-\z)^\vee$. In particular, each extreme ray $\w = (h_P(\v), \v)$ of $F$ corresponds to a vertex $\v'' = \frac{1}{h_{P-\z}(\v)} \v$ of $F''\in T''$. Since $\Vol((P-z)^\vee) = \sum_{F''\in T''} \Vol(F'')$, comparing with (\ref{eqn: adj/B_z}), it suffices to show that the following identity holds for each $F \in T$:
\begin{equation}\label{eq:volsimplex}  \frac{|\det(F)| }{\displaystyle\prod_{\w \in V(F)} h_{P-\z}(\v)} = \Vol(F''). 
\end{equation}

Let $F = \Cone \{\w_1, \dots, \w_{d+1} \}$, then $F'' = \Cone \{\bar{\v}''_1, \dots, \bar{\v}''_{d+1} \}$. Using the determinant formula to represent the volume of $F''$ (up to a sign), we obtain
\begin{eqnarray*}
\begin{vmatrix}
                1 & 1 & \dots & 1\\
            \vert & \vert & \dots & \vert\\
            \v''_1   & \v''_2 & \dots & \v''_{d+1}  \\
            \vert & \vert & \dots & \vert
        \end{vmatrix}
&=& \prod h_{P-\z}(\v_i) \begin{vmatrix}
                h_{P}(\v_1) &  h_{P}(\v_2) & \dots &  h_{P}(\v_{d+1})\\           
                \vert & \vert & \dots & \vert\\
            \v_1   & \v_2 & \dots & \v_{d+1}  \\
            \vert & \vert & \dots & \vert
        \end{vmatrix}\\
&=& \big(\prod h_{P-\z}(\v_i) \big) \det(F).
\end{eqnarray*}
Assuming $P-\z$ contains the origin in its interior, we have $\prod h_{P-\z}(\v_i) >0$, and so $\Vol(F'') = \big(\prod h_{P-\z}(\v_i) \big) |\det(F)|$, as desired.
\end{proof}

\begin{ex}
Consider the polytope $P$ from \cref{ex: quadrilateral-dual-volume}. The polyhedral cone $C(P)$ is generated by $V(C(P)) = \{ (1,-1,1), (1,0,1), (1,1,-1), (1,-1,-1)\}$. Its dual cone $C := C(P)^*$ is generated by 
\[V(C) = \{\w_1 := (1,0,1), \w_2 := (-1,1,0), \w_3 := (1.0,-1), w_4 := (5,-2,-1) \}. \]
There is a triangulation of $C$ in to $F_1 = \Cone(\w_1, \w_2, \w_4)$ and $F_2 = \Cone(\w_2, \w_3, \w_4)$. Let's compute the adjoint of $C$ using this triangulation:
\begin{eqnarray*}
\adj_C(\bar{\z}) &=& |\det(F_1)|\langle \w_3, \bar{\z}\rangle + |\det(F_2)|\langle \w_1, \bar{\z}\rangle \\
&=& 4(1-z_2) + 2(1-z_2)\\
&=& 6-2z_2.
\end{eqnarray*}
On the other hand, the dual volume function we computed in \cref{ex: quadrilateral-dual-volume} is
\begin{eqnarray*}
\VolC_{\z}(P) &=&\frac{1}{(1{+}z_2)({-}1{+}z_1)}+\frac{1}{({-}1{+}z_1)(1{-}z_2)}+\frac{2}{(1{-}z_2)(5{-}2z_1{-}z_2)}+\frac{2}{(5{-}2z_1{-}z_2)(1{+}z_2)}\\
&=& \frac{6-2z_2}{(1{+}z_2)({-}1{+}z_1)(1{-}z_2)(5{-}2z_1{-}z_2)} =: \frac{A_\z(P)}{ B_\z(P)}.
\end{eqnarray*}
Comparing the two formulae, we observe that $\adj_C(\bar{\z}) = A_{\z}(P)$, consistent with \cref{thm:adjoint}.
\end{ex}

\section{Canonical forms}\label{sec:canform}
\def\PP{{\mathbb{P}}}
In this section, we connect the dual volume of a polyhedron to positive geometry.  Positive geometries are semialgebraic sets endowed with a distinguished meromorphic form called the \emph{canonical form}.

Full-dimensional projective polytopes $\bar P \subseteq \PP^d(\R)$ are examples of positive geometries \cite{ABL,lam2022invitation}.  Let $P \subseteq \R^d$ be a full-dimensional polyhedron not containing any lines.  View the projective space $\PP^d = \R^d \sqcup H_\infty$ as the union of the affine space and a plane at infinity.  Then the closure $\bar P \subseteq \PP^d(\R)$ is a \emph{projective polytope}.  If $P$ is bounded, that is, it is a polytope, then $\bar P = P$.  If $P$ is an unbounded polyhedron, then $\bar P$ is the union of $P$ together with the face $\bar P \cap H_\infty$ at infinity, which could be a vertex, a facet, or a face of intermediate dimension.

The following result states that (oriented) projective polytopes are positive geometries.
\begin{theorem}[{\cite[Section 6.1]{ABL}}]\label{thm:canformexist}
For each full-dimensional oriented projective polytope $\bar P \subseteq \PP^d(\R)$, there exists a unique nonzero meromorphic $d$-form $\Omega(\bar P)$ on $\PP^d$, called the \emph{canonical form}, satisfying the following recursive property.  If $d > 0$, the canonical form $\Omega(\bar P)$ has poles exactly along the facets of $\bar P$, all these poles are simple, and for each facet $F$, the residue of $\Omega$ along $F$ is given by
$$
\Res_F(\Omega(\bar P)) = \Omega(F).
$$
If $d = 0$, then $\bar P$ is a point and we have $\Omega = \pm 1$ depending on the orientation. 
\end{theorem}

\begin{remark}
    Since all the projective polytopes we consider in this paper arise as closures of polyhedra $P \subseteq \R^d$, they inherit a natural orientation from any chosen orientation of $\R^d$.  Thus in the following we will ignore the orientation.
\end{remark}

Our main result in this section connects dual volumes and canonical forms of polyhedra.
\begin{theorem}\label{thm:canformdualvol}
Let $P \subseteq \R^d$ be a full-dimensional polyhedron that does not contain lines.  Then the canonical form $\Omega(\bar P)$ of the projective polytope $\bar P$ is given by 
\[\Omega(\bar P) = \Vol^\vee_\z(P) dz_1 dz_2 \cdots dz_d.\]
\end{theorem}

We shall use the following basic result.
\begin{thm}[\cite{ABL}]\label{thm:canformsub}
Suppose that $\{Q_1,\ldots,Q_n\}$ is a subdivision of a projective polytope $\bar P$.  Then \[\sum_{i=1}^n \Omega(Q_i) = \Omega(\bar P).\]
\end{thm}
\begin{proof}[Sketch proof.]
Let $\Omega = \sum_{i=1}^n \Omega(Q_i)$.  To check that $\Omega$ is the canonical form of $\bar P$, we need to check the recursive properties in \cref{thm:canformexist}. The possible poles of $\Omega$ are along the facets $F$ of the $Q_i$'s.  If $F$ is not a facet of $\bar P$, then one checks that the residues of $\Omega(Q_i)$ along $F$ cancel out.  If $F$ is a facet of $\bar P$, then we calculate
$$
\Res_F \Omega = \sum_{i=1}^n \Omega(F \cap Q_i),
$$
and we note that (keeping the correct dimensional pieces of) $\{F \cap Q_1,\ldots,F \cap Q_n\}$ gives a subdivision of $F \cap \bar P$.  By induction on the dimension, we thus have $\Res_F \Omega = \Omega(F \cap \bar P)$, confirming that $\Omega$ has the correct residues.
\end{proof}

\begin{proof}[Proof of \cref{thm:canformdualvol}]
Taking a (projective) triangulation of $\bar P$, and using \cref{thm:valuation} and \cref{thm:canformsub}, we may reduce to the case that $\bar P$ is a projective simplex.  First suppose that $P \subseteq \R^d$ is a (bounded) simplex.  Then $(P-\z)^\vee$ is also a simplex whenever $\z$ is in the interior of $P$, and its volume is given in \eqref{eq:volsimplex} by the ratio $\frac{|\det(\w_1,\ldots,\w_{d+1})| }{\prod_{i=1}^{d+1} h_{P-\z}(\v_i)}$, where the $\v_i$'s are the generators of the rays of the normal fan of $P$, and $\w_i = (h_P(\v_i),\v_i) \in \R^{d+1}$ for each $i$.  Thus $\VolC_\z(P)$ has exactly $d+1$ simple poles in $\z$, each of the form
$$
h_{P-\z}(\v_i) = h_P(\v_i) + \langle \v_i, \z\rangle.
$$
The vanishing of $h_{P-\z}(\v_i)$ is exactly the condition that $\z$ lies in the supporting hyperplane of the facet $F_i$ of $P-\z$ that is normal to $\v_i$.  We have
\begin{equation}\label{eq:resF}
\Res_{F_i}\frac{|\det(\w_1,\ldots,\w_{d+1})| }{\prod_{i=1}^{d+1} h_{P-\z}(\v_i)} dz_1 dz_2 \cdots dz_d = \left.\frac{|\det(\w_1,\ldots,\w_{d+1})| }{\prod_{j \neq i}^{d+1} h_{P-\z}(\v_j)}\right|_{h_{P-\z}(\v_i) =0} \frac{dz_1 dz_2 \cdots dz_d}{d(h_{P-\z}(\v_i))}.
\end{equation}
The meaning of $\frac{dz_1 dz_2 \cdots dz_d}{d(h_{P-\z}(\v_i))}$ is that it is the $(d-1)$-form $\eta$ on $F_i$ such that $d(h_{P-\z}(\v_i)) \wedge \eta = dz_1 dz_2 \cdots dz_d$.  Now, for $j \neq i$ the vectors $\v_j$'s, viewed as vectors in $\R^d/\v_i$, are exactly the normal vectors to the simplex $F_i$ in its supporting hyperplane.  It follows that the right hand side of \eqref{eq:resF} is the dual volume function of $F_i$ in its supporting hyperplane, where the hyperplane is equipped with a volume measure scaled according to 
$$
|\det(\w_1,\ldots,\w_{d+1})|\frac{dz_1 dz_2 \cdots dz_d}{d(h_{P-\z}(\v_i))}.
$$
By induction, it follows that $\VolC_\z(P) dz_1 dz_2 \cdots dz_d$ has the same poles and residues as $\Omega(\bar P)$ and by the uniqueness of the canonical form (\cref{thm:canformexist}) we conclude that $\Omega(\bar P) = \VolC_\z(P) dz_1 dz_2 \cdots dz_d$.  Essentially the same argument holds if $\bar P$ is a projective simplex such that the face at infinity is not a facet.  

When $\bar P$ has a facet at infinity, then the normal fan of $P$ only has $d$ rays, so $\VolC_\z(P)$ has only $d$ linear factors in the denominator.  However, in this case the form $\VolC_\z(P) dz_1 dz_2 \cdots dz_d$ has an additional simple pole at infinity, and a similar residue calculation holds.
\end{proof}

\begin{remark}
    An alternative proof is to directly compare \eqref{eq:volsimplex} with the canonical form of a projective simplex, given in \cite[(5.7)]{ABL}.
\end{remark}

\begin{remark}
Note that if $P \subseteq \R^d$ contains a line then the closure $\bar P \subseteq \Pbb^d(\R)$ is not a projective polytope.  Nevertheless, it can be considered a pseudo-positive geometry \cite{ABL}, and has canonical form $0$, agreeing with \cref{lem:Vol0}.
\end{remark}

Combining \cref{thm:canformdualvol} with \cref{thm:adjoint}, we obtain the following; see \cite{Gaetz} for an alternative approach.
\begin{cor}\label{cor:adjointzero}
Let $P \subsetneq \R^d$ be a full-dimensional polyhedron not containing lines.  Then $\Omega(\bar P) = \frac{\adj_{C(P)^*}(\bar \z)}{B_\z(P)} dz_1 dz_2 \cdots dz_d$.  Thus the zero set of $\Omega(\bar P)$ is the vanishing set of $\adj_{C(P)^*}$, known as the \emph{adjoint hypersurface} of $P^\vee$.
\end{cor}

\def\PGL{{\rm PGL}}
The definition of the canonical form $\Omega(\bar P)$ for a projective polytope is invariant under automorphisms of $\PP^d$.  In particular, if $\bar P$ and $\bar P'$ are related by a projective transformation $g \in \PGL(d+1)$, we have $\Omega(\bar P) = g^* \Omega(\bar P')$.  This leads to the following (perhaps unexpected) property of the dual volume function.

\begin{cor}\label{cor:proj}
    Let $P \subseteq \R^d$ and let $P' = g \cdot P$ where $g \in \PGL(d+1)$ is a projective transformation.  Let $\z' = (z'_1,\ldots,z'_d)$ be the image of the coordinates $z_1,\ldots,z_d$ under $g$.  Then 
    $$
    \Vol_\z(P) =  \det\left(\frac{\partial z'_j}{\partial z_i}\right) \Vol_{\z'}(P').
    $$
\end{cor}

\begin{example}
    Suppose that $P = [-1,2]$ and $P' = [-2,4]$.  Then we have $z' = 2z$.  So
    $$
    \Vol_z([-1,2]) = 2 \Vol_{2z}([-2,4])  = \frac{1}{4-(2z)} - \frac{1}{-2-(2z)} = \frac{1}{2-z} - \frac{1}{-1-z}.
    $$
\end{example}

\begin{example}
Let $P \subseteq \R^2$ be the triangle with vertices $(0,3),(0,1),(1,1)$ and $P'$ be the triangle with vertices $(0,3/2),(0,1/2),(1,1/2)$.  The projective transformation $(z_1,z_2)\mapsto (z'_1,z'_2) = (z_1/(z_2+1),z_2/(z_2+1))$ takes $P$ to $P'$.  We have
$$
\VolC_\z(P)  =\frac{1}{(1+z_2)^3} \frac{1}{z'_1(z'_2 - \frac{1}{2})(\frac{3}{2} - z'_1 - z'_2)}= \frac{4}{z_1 (3 - 2 z_1 + z_2) ( z_2-1)}.
$$
\end{example}

\section{Minkowski sum}\label{sec:Minkowski}
For convex bodies $P,Q \subseteq \R^d$, the \emph{Minkowski sum} is defined by
$$
P+Q := \{\p+\q \;|\; \p \in P, \q\in Q\}.
$$
More generally, for closed convex bodies $S_1,\ldots,S_r \subseteq \R^d$ and positive real numbers $x_1,\ldots,x_r$, we define
$$
\x \S = x_1 S_1 + \cdots + x_r S_r := \left\{\p_1 + \p_2 + \cdots + \p_r \;|\; \frac{1}{x_i}\p_i \in S_i \text{ for all } i\right\}.
$$

Recall in \eqref{eq:support} the support function $h_S$ of $S$ is defined for, and uniquely determines, any closed convex body $S$.  It follows directly from the definition that, the support function of the Minkowski sum $\x\S$ is given by 
\begin{equation}\label{eq:hxS}
h_{\x\S} = x_1 h_{S_1} + x_2 h_{S_2} + \cdots + x_r h_{S_r}.
\end{equation}
We use this formula to extend the definition of $\x \S$ to arbitrary real parameters $x_1,\ldots,x_r \in \R$.  Namely, $\x\S$ is the ``object" with support function $h_{\x\S}$ given by \eqref{eq:hxS}.  It is an interesting question to characterize the parameters $\x$ such that $h_{\x\S}$ is the support function of a closed, convex set.  Note that the support functions of closed, bounded, convex sets are exactly the positive homogeneous, convex, real-valued functions on $\R^d$.

Suppose that $h_{\x\S} = h_{T}$ for some closed convex set $T$. Furthermore, suppose that $x_1,\ldots,x_s \geq 0$ and $x_{s+1},\ldots,x_r \leq 0$.  It follows then that we have 
$$
x_1 h_{S_1} + \cdots + x_s h_{S_s} = h_T + (-x_{s+1}) h_{S_{s+1}} + \cdots + (-x_r) h_{S_r}.
$$
Equivalently, $T + \sum_{i=s+1}^r (-x_i) S_i = \sum_{j=1}^s x_j S_j$, where we are taking usual Minkowski sums of convex bodies.  Thus our definition of $\x\S$ is consistent with Minkowski sums when $\x\S = T$ is a closed convex set.

For the remainder of this section, we focus our attention on the case of a sequence $\P = (P_1,\ldots,P_r)$ of polyhedra in $\R^d$.  In this case, each $h_{P_i}$ is a piecewise-linear function, and the common domains of linearity of $h_{P_1},\ldots,h_{P_r}$ give a fan $\F$ in $\R^d$, and this fan is equal to the normal fan $\N(P)$ of the Minkowski sum $P = P_1+ \cdots + P_r$.  The fan $\F$ is complete exactly when $P$ (and thus $P_1,\ldots,P_r$) is a polytope.    In the following, we make the assumption that 
\begin{equation}\label{eq:full}
\mbox{the Minkowski sum $P = P_1 + \cdots + P_r$ is full-dimensional in $\R^d$.}
\end{equation}
Thus the maximal cones of $\F$ are pointed.

The following lemma is straightforward.
\begin{lemma}\label{lem:hxP}
For any $x_i \in \R$, the piecewise-linear function $h_{\x\P}$ restricts to a linear function on any of the maximal cones of $\F$.  In particular, the function $h_{\x\P}$ is uniquely determined by its values on the rays of $\F$.
\end{lemma}

Call $\P$ \emph{regular} if $\P$ satisfies \eqref{eq:full} and there exists a point $\x \in \R^r$ so that $h_{\x\P}(\v) \neq 0$ for all rays $\R_{\geq 0} \cdot \v$ of $\F$.  Not all sequences $\P$ of polytopes are regular, even when \eqref{eq:full} is satisfied.  For example if $\P = (P)$ is a single full-dimensional polytope, then $\P$ is regular if and only if $P$ is non-codegenerate. 

\begin{proposition}\label{prop:reg}
The sequence $\P = (P_1,\ldots,P_r)$ is regular if and only if for each ray $\R_{\geq 0} \cdot \v$ of $\F$, there exists a $j$ such that $h_{P_j}(\v) \neq 0$.
\end{proposition}
\begin{proof}
For any ray $\v$ of $\mathcal{F}$, if $h_{\x\P}(\v) = x_1 h_{P_1}(\v) + \dots + x_r h_{P_r}(\v) \neq 0$, then one of the summands has to be nonzero, therefore $h_{P_j}(\v) \neq 0$ for some $j$. One the other hand, if $h_{P_j}(\v) \neq 0$, then letting $\x = \e_j$ yields that $h_{\x\P} (\v) = h_{P_j}(\v)\neq 0$.
\end{proof}

Since for any ray $\R_{\geq 0}\cdot \v$ there exists some standard unit vector $P= \e_i$ such that $h_P(\v) = -\langle \v , \e_i\rangle \neq 0$. The following proposition follows from Proposition \ref{prop:reg}.

\begin{proposition}\label{prop:mpxz-regular}
Suppose that \eqref{eq:full} is satisfied and $\P = (P_1,P_2,\ldots,P_r, -\e_1,-\e_2,\ldots,-\e_d)$.  Then $\P$ is regular.
\end{proposition}


\section{Dual mixed volume}\label{sec:DMV}

\subsection{Definition of dual mixed volume}
\begin{defin}\label{def:dual-mixed-volume}
Let $\P = (P_1,\ldots,P_r)$ be a regular sequence of polyhedra with normal fan $\N(P)$ and let $\x = (x_1, \dots, x_r)$. The \emph{dual mixed volume rational function} $m_\P(\x)$ is
$$
m_\P(\x) := \VolC(h_{\x\P}) = f_{\N(P),\mathbf{v}}(h_{\x\P}(\mathbf{v_1}),\ldots,h_{\x\P}(\mathbf{v_g}))
$$
where $\v_1, \dots, \v_g$ are the generating rays of $\mathcal{N}(P)$, with notation as in \cref{def:dual-volume}.  For any sequence $\P$ that is full-dimensional, also define
$$
m_\P(\x,\z) := m_{(P_1,P_2,\ldots,P_r,-\e_1,\ldots,-\e_r)}(x_1,\ldots,x_r,z_1,\ldots,z_d) = \VolC(h_{\x\P-\z}).
$$
\end{defin}
By \cref{prop:mpxz-regular}, $m_{\P}(\x,\z)$ is always well-defined. If $\P$ is not full-dimensional,  we set $m_\P(\x) := 0$ and $m_{\P}(\x,\z):=0$.

By \cref{def:universal-dual-volume} we can write $m_{\P}$ as a rational function in $x_1, x_2, \dots, x_r$, with degree $-d$. The denominator is a product of linear factors, each corresponding to a ray $\v_i$ of $\mathcal{N}(P)$.  As in \eqref{eqn: VolC one fraction}, we have
\begin{align} \label{eq:DMVratio}
 m_{\P}(\x)
 &= \frac{\displaystyle\sum\limits_{C\in \mathcal{T}}|\det(C)| \displaystyle \prod\limits_{\v\in V(\N(P))\setminus V(C)} h_{\x\P}(\v)}{\displaystyle\prod_{\v\in V(\N(P))} h_{\x\P}(\v)}
 \end{align}
where $T$ is any triangulation of $\N(P)$, $C$ is a simplicial cone in $T$ generated by the rays $V(C) = \{\mathbf{v}_{i_1},\ldots,\mathbf{v}_{i_d}\}$, and $\det(C) = \det(\v_{i_1} , \dots, \v_{i_d})$.

Note that $m_\P(\x)$ generalizes both $\VolC$ and $\VolC_\z$.  For a full-dimensional polytope $P$, we have the specializations
$$
m_P(1) = \VolC(P), \qquad m_{(P,-\e_1,-\e_2,\ldots,-\e_d)}(1,z_1,z_2,\ldots,z_d) = \VolC_\z(P).
$$

\subsection{Mixed polar sets}

We next consider the mixed analogue of \cref{lem:same-function-interior}. For the polar dual polyhedron $(\x\P)^\pol$ to exist, we need to have the origin inside the interior of $\x\P$. The support function $h_{\x\P}$ can help us characterize such a condition.

\begin{proposition}\label{prop: support function and containing origin}
    Let $P \subseteq \R^d$ be a polyhedron.  Then $h_P \geq 0$ if and only if $\0 \in P$, and $h_P > 0$ on $\R^d \setminus \0$ if and only if $\0$ is  in the interior of $P$.
\end{proposition}

From now on, we assume that $\P$ is a regular sequence of polyhedra. We define the following polyhedral cones in $\R^r$:
\begin{align*}
A_\P &:= \{ \x \mid h_{\x\P} \geq 0 \text{ on } \R^d\},\text{ and} \\
A_\P^\circ &= \{\x \mid h_{\x\P} > 0 \text{ on } \R^d\setminus \0\}.
\end{align*} 
By \cref{lem:hxP}, $\x \in \P$ if and only if $h_{\x\P}(\v_i) \geq 0$ for each ray $\R_{\geq 0} \cdot \v_i$ of the fan $\F = \mathcal{N}(P)$.  Thus the cone $A_\P$ is a nonempty closed polyhedral cone, and $A_\P^\circ$ is an open (possibly empty) polyhedral cone, the interior of $A_\P$.  By Proposition \ref{prop: support function and containing origin} the cones $A_\P$ and $A_\P^\circ$ correspond to polyhedra $\x\P$ that contain the origin and those that contain the origin in the interior, respectively.

When $\x\in A_\P^\circ$, we define the polar dual of $\x\P$ to be the set 
$$
(\x\P)^\pol:= \{\w \in \R^d \mid h_{\x\P}(\w) \leq 1\} \subseteq \R^d,
$$ which is clearly closed and bounded, and contains the origin.  When $\x\P$ is a polytope, then this definition agrees with the usual polar of $\x\P$.
\begin{proposition}
Suppose that $\x \in A_\P^\circ$.  Then 
$$
m_\P(\x) = \Vol((\x\P)^\pol).
$$
\end{proposition}
\begin{proof}
This result follows in the same way as \cref{lem:same-function-interior}.
\end{proof}

\subsection{Integral formulae}

By the definition of $m_\P(\x, \z)$, the following is a straightforward corollary of Theorem \ref{thm: dual vol is integral}.

\begin{theorem}
Let $\P = (P_1, \dots, P_r ) $ be a sequence of polyhedra in $\R^d$, $\x = (x_1, \dots, x_r)$ in $\R^r$, and $\z \in \R^d$ be any point such that $\x\P - \z$ is non-codegenerate. Then
\begin{align*}
 m_\P(\x,\z) = \int_{\R^d} \exp(-h_{\x\P-\z}(\v)) d\v.
\end{align*}
\end{theorem}

In the case that $\x\P$ is non-codegenerate, then we can let $\z$ go to the origin and take limits on both sides to get:

\begin{cor}
    Let $\P = (P_1, \dots, P_r ) $ be a sequence of polyhedra in $\R^d$ and $\x = (x_1, \dots, x_r) \in \R^r$. If $\x\P$ is non-codegenerate then
\begin{align*}
 m_\P(\x) = \int_{\R^d} \exp(-h_{\x\P}(\v)) d\v.
\end{align*}
\end{cor}

\section{Fine mixed cells}\label{sec:fine}
We consider the dual mixed volumes of fine mixed cells.  These are the ``minimal" dual mixed volume functions that serve as our building blocks. 
\begin{defin}\label{def:fine-mixed-cell}
A \emph{fine mixed cell} is a sequence $\P = (P_1,P_2,\ldots,P_r)$ of simplices in $\R^d$ satisfying both of the conditions:
\begin{itemize}
\item $\dim(P_1) + \dim(P_2) + \cdots + \dim(P_r) = d$, and
\item $\dim(P_1+P_2+ \cdots + P_r) = d$.
\end{itemize}
\end{defin}
Because of the dimension condition, it is clear that $P_1+\cdots+P_r$ is non-codegenerate if and only if all of the $P_1,\ldots,P_r$ are non-codegenerate. In this case, we say that the fine mixed cell $\P$ is non-codegenerate.

Let $\P$ be a fine mixed cell. Write $d_i=\dim(P_i)$ for $i=1,\ldots,r$ and write $P=P_1+\cdots+P_r$. We first describe the normal fan $\mathcal{N}(P)$. For each $P_i$ with $d_i>0$, let $\v_{i,1},\v_{i,2},\ldots \v_{i,d_i+1}$ be rays in the normal fan of $P_i$ chosen to be normal to $\mathrm{span}(P_k)$ for each $k\neq i$, where $\mathrm{span}(P_k)$ is the affine span of the simplex $P_i$. We scale each $\v_{i,a}$ so that \[\max_{\y\in P_i}\langle\y,\v_{i,a}\rangle-\min_{\y\in P_i}\langle\y,\v_{i,a}\rangle=1.\]
If $d_i=\dim(P_i)=0$, we just assume $\v_{i,1}=\mathbf{0}$.
\begin{lemma}\label{lem:fine-mixed-cell-normal-fan}
With notations as above, $\sum_{a=1}^{d_i+1}\v_{i,a}=0$ for each $i=1,\ldots,r$. Moreover, $\mathcal{N}(P)$ has maximal cones $\{C_{\mathbf{a}}\:|\:\mathbf{a}=(a_1,\ldots,a_r)\in[d_1{+}1]\times\cdots\times[d_r{+}1]\}$ where \[C_{\mathbf{a}}=C_{a_1,\ldots,a_r}=\mathrm{span}_{\geq0}\bigcup_{i=1}^r\{\v_{i,a}\:|\:a\neq a_i\}.\]
\end{lemma}
\begin{proof}
Let the vertices of the simplex $P_i$ be $\p_{i,1},\ldots,\p_{i,d_i+1}$ labeled so that the facet of $P_i$ with respect to $\v_{i,a}$ is the convex hull of $\{\p_{i,j}\:|\: j\neq a\}$. In other words, $\min_{\y\in P_i}\langle\y,\v_{i,a}\rangle=\langle\p_{i,j},\v_{i,a}\rangle$ for all $j\neq a$, and also $\max_{\y\in P_i}\langle\y,\v_{i,a}\rangle=\langle\p_{i,a},\v_{i,a}\rangle$. By construction of $\v_{i,a}$'s, we see that $\langle\p_{i,a},\v_{i,a}\rangle-\langle\p_{i,j},\v_{i,a}\rangle=1$ for all $j\neq a$. Let $c=\sum_{j=1}^{d_i+1}\langle\p_{i,j},\v_{i,j}\rangle$, and let $\v=\sum_{j=1}^{d_i+1}\v_{i,j}$. The above equality gives \[c-\langle\p_{i,a},\v\rangle=\sum_{j=1}^{d_i+1}\langle\p_{i,j}-\p_{i,a},\v_{i,j}\rangle=d_i\] for all $a\in[d_i{+}1]$. In particular, this says that the vector $\v$ has the same inner product with all vertices in $P_i$, meaning that $\v$ is normal to the affine span of $P_i$. But each $\v_{i,a}$ is also normal to the affine span of $P_k$ for each $k\neq i$. Thus, $\v$ is normal to the affine span of $P_1,\ldots,P_r$, and as $\P$ is full-dimensional, $\v=0$. This concludes the first part.

Next, each maximal cone of $\N(P)$ corresponds to a vertex of $P$, which can be uniquely written as $\p_{1,a_1}+\cdots+\p_{r,a_r}$ for some $\a=(a_1,\ldots,a_r)\in[d_1{+}1]\times\cdots\times[d_r{+}1]$. For every $\v\in\R^d$, in order for $\min_{\y\in P}\langle\y,\v\rangle=\langle\p_{1,a_1}+\cdots+\p_{r,a_r},\v\rangle$, it is equivalent to split $\v=\v_1+\cdots+\v_{r}$ where $v_k$ is normal to $\mathrm{span}(\cup_{j\neq k}P_j)$ and require that $\min_{\y\in P_i}\langle\y,\v_i\rangle=\langle\p_{i,a_i},\v_i\rangle$. This precisely means $\v_i\in\mathrm{span}_{\geq0}\{\v_{i,a}\:|\:a\neq a_i\}$ by construction of $\v_{i,a}$'s and $\p_{i,a}$'s. Combining these statements, we obtain the description of the maximal cones in $\N(P)$.
\end{proof}
\begin{defin}\label{def:kappa-determinant}
Let $\P=(P_1,\ldots,P_r)$ be a fine mixed cell in $\R^d$ and let $\v_{i,a}$ be as above, for $i\in[r]$ and $a\in[d_i+1]$. Define $\kappa(\P):=|\det(\v_{i,j})_{i\in[r],j\in[d_i]}|$.
\end{defin}
Note that since $\v_{i,1}+\v_{i,2}+\cdots+\v_{i,d_i+1}=0$ (Lemma~\ref{lem:fine-mixed-cell-normal-fan}), we can choose any of the $d_i$ vectors among $\{\v_{i,1},\ldots,\v_{i,d_i+1}\}$ for each $i$, for a total of $d_1+\cdots+d_r=d$ vectors, to form a determinant and take the absolute value to obtain $\kappa(\P)$. In other words, $\kappa(\P)$ is well-defined.
\begin{prop}\label{prop:fine-mixed-dual-volume}
Let $\P=(P_1,\ldots,P_r)$ be a fine mixed cell in $\R^d$ and let $\v_{i,a}$ be as above. Then its dual mixed volume function is \[m_{\P}(\x,\z)=\kappa(\P)\frac{\prod_{i:d_i>0}x_i}{\prod_{i:d_i>0}\prod_{a=1}^{d_i+1}(h_{\x\P}(\v_{i,a})+\langle\v_{i,a},\z\rangle)}\]
where $h_{\x\P}(\v)=\sum_{i=1}^r x_ih_{P_i}(\v)$. In particular, when $\P$ is non-codegenerate, \[m_{\P}(\x)=\kappa(\P)\frac{\prod_{i:d_i>0}x_i}{\prod_{i:d_i>0}\prod_{a=1}^{d_i+1}h_{\x\P}(\v_{i,a})}\]
\end{prop}
\begin{proof}
Let the vertices of $P_i$ be $\{\p_{i,1},\ldots,\p_{i,d_i+1}\}$ as in the proof of Lemma~\ref{lem:fine-mixed-cell-normal-fan}. Let's analyze $h_{\x\P}(\v_{i,a})$, which is a linear function in $x_1,\ldots,x_r$. Note that for $j\neq i$, $\langle\y,\v_{i,a}\rangle$ is constant for $\y\in P_j$ by construction of $\v_{i,a}$'s. In other words, $h_{P_j}(\v_{i,a})=-\langle\y_j,\v_{i,a}\rangle$ for any fixed $\y_j\in P_j$. Also recall that $h_{P_i}(\v_{i,a})=-\langle \p_{i,j},\v_{i,a}\rangle$ for any $j\neq a$, and that $\langle\p_{i,a},\v_{i,a}\rangle-\langle\p_{i,j},\v_{i,a}\rangle=1$ for all $j\neq a$. We then have, for $d_i>0$,
\begin{align*}
\sum_{a=1}^{d_i+1}h_{\x\P}(\v_{i,a})=&\sum_{j\neq i}\left(\sum_{a=1}^{d_i+1}\langle \y_j,\v_{i,a}\rangle\right)\cdot x_j+\sum_{a=1}^{d_i+1}\big(-\min_{\y\in P_i}\langle \y,\v_{i,a}\rangle\big)x_i\\
=&0+\big(-\langle \p_{i,2},\v_{i,1}\rangle-\langle\p_{i,1},\v_{i,2}\rangle-\langle\p_{i,1},\v_{i,3}\rangle-\cdots-\langle\p_{i,1},\v_{i,d_i+1}\rangle\big)x_i\\
=&\big(1-\langle \p_{i,1},\v_{i,1}\rangle-\langle\p_{i,1},\v_{i,2}\rangle-\langle\p_{i,1},\v_{i,3}\rangle-\cdots-\langle\p_{i,1},\v_{i,d_i+1}\rangle\big)x_i\\
=&x_i,
\end{align*}
where we have repeatedly used the fact that $\sum_{a=1}^{d_i+1}\v_{i,a}=0$. Also,
\begin{equation}\label{eqn:sum-of-hs-in-cell}
\sum_{a=1}^{d_i+1}h_{\x\P-\z}(\v_{i,a})=\sum_{a=1}^{d_i+1}h_{\x\P}(\v_{i,a})+\left\langle \sum_{a=1}^{d_i+1}\v_{i,a},\z \right\rangle=x_i.
\end{equation}
Using Definition~\ref{def:dual-volume} and the description of the normal fan $\N(P)$ in Lemma~\ref{lem:fine-mixed-cell-normal-fan},
\begin{align*}
m_{\P}(\x,\z)=&\VolC(\x\P)=\sum_{(a_1,\ldots,a_r)\in[d_1{+}1]\times\cdots\times[d_r{+}1]}\kappa(\P)\prod_{i:d_i>0}\prod_{a\neq a_i}\frac{1}{h_{\x\P-\z}(\v_{i,a})}\\
=&\kappa(\P)\prod_{i:d_i>0}\sum_{a_i=1}^{d_i+1}\prod_{a\neq a_i}\frac{1}{h_{\x\P-\z}(\v_{i,a})}=\kappa(\P)\prod_{i:d_i>0}\frac{\sum_{a=1}^{d_i+1}h_{\x\P-\z}(\v_{i,a})}{\prod_{a=1}^{d_i+1}h_{\x\P-\z}(\v_{i,a})}\\
=&\kappa(\P)\prod_{i:d_i>0}\frac{x_i}{\prod_{a=1}^{d_i+1}h_{\x\P-\z}(\v_{i,a})}=\kappa(\P)\frac{\prod_{i:d_i>0}x_i}{\prod_{i:d_i>0}\prod_{a=1}^{d_i+1}h_{\x\P-\z}(\v_{i,a})}\\
=&\kappa(\P)\frac{\prod_{i:d_i>0}x_i}{\prod_{i:d_i>0}\prod_{a=1}^{d_i+1}(h_{\x\P}(\v_{i,a})+\langle\v_{i,a},\z\rangle)}
\end{align*}
as desired.
\end{proof}
Proposition~\ref{prop:fine-mixed-dual-volume} says that our building blocks, fine mixed cells, have dual mixed volume being essentially a monomial. In the remainder of the text, we compute many dual mixed volumes using \cref{prop:fine-mixed-dual-volume} and \cref{thm: mixed subdivision formula}. 
\begin{example}
Consider the following fine mixed cell $\P=(P_1,P_2)$ that is non-codegenerate with parameters $d_1=2$, $d_2=1$, $d=3$. That is, $P_1$ is a $2$-simplex and $P_2$ is a line segment:
\begin{align*}
P_1=&\text{conv}\big(\p_{1,1}=(0,2,1),\p_{1,2}=(1,-1,1),\p_{1,3}=(-1,0,1)\big),\\
P_2=&\text{conv}\big(\p_{2,1}=(0,0,-2),\p_{2,2}=(-1,1,-1)\big).
\end{align*}
We will not plot the polytope $\x\P$ here since it will not be helpful towards the computation, but the readers can view $\x\P$ as a slanted triangular prism. 

To compute the ray $\v_{1,1}$, we require that \[\langle\p_{1,1},\v_{1,1}\rangle-1=\langle\p_{1,2},\v_{1,1}\rangle=\langle\p_{1,3},\v_{1,1}\rangle,\ \text{and }\langle\p_{2,1},\v_{1,1}\rangle=\langle\p_{2,2},\v_{1,1}\rangle.\]
Thus, $\v_{1,1}=\frac{1}{5}(1,2,-1)$. Similarly, $\v_{1,2}=\frac{1}{5}(2,-1,3)$, $\v_{1,3}=\frac{1}{5}(-3,-1,-2)$, and $\v_{2,1}=(0,0,-1)$, $\v_{2,2}=(0,0,1)$. The determinant $\det(\v_{1,1},\v_{1,2},\v_{2,1})=\frac{1}{5}$ so $\kappa(\P)=\frac{1}{5}$. Next,
\begin{align*}
h_{\x\P}(\v_{1,1})=&(-\min_{\y\in P_1}\langle\y,\v_{1,1}\rangle)x_1+(-\min_{\y\in P_2}\langle\y,\v_{1,1}\rangle)x_2=\frac{1}{5}(2x_1-2x_2),\\
h_{\x\P}(\v_{1,2})=&(-\min_{\y\in P_1}\langle\y,\v_{1,2}\rangle)x_1+(-\min_{\y\in P_2}\langle\y,\v_{1,2}\rangle)x_2=\frac{1}{5}(-x_1+6x_2),\\
h_{\x\P}(\v_{1,3})=&(-\min_{\y\in P_1}\langle\y,\v_{1,3}\rangle)x_1+(-\min_{\y\in P_2}\langle\y,\v_{1,3}\rangle)x_2=\frac{1}{5}(4x_1-4x_2),\\
h_{\x\P}(\v_{2,1})=&(-\min_{\y\in P_1}\langle\y,\v_{2,1}\rangle)x_1+(-\min_{\y\in P_2}\langle\y,\v_{2,1}\rangle)x_2=(x_1-x_2),\\
h_{\x\P}(\v_{2,2})=&(-\min_{\y\in P_1}\langle\y,\v_{2,2}\rangle)x_1+(-\min_{\y\in P_2}\langle\y,\v_{2,2}\rangle)x_2=(-x_1+2x_2).\\
\end{align*}
Proposition~\ref{prop:fine-mixed-dual-volume} gives us that \[m_{\P}(x_1,x_2)=\frac{1}{5}\cdot\frac{x_1x_2}{\frac{1}{5}(2x_1-2x_2)\frac{1}{5}(-x_1+6x_2)\frac{1}{5}(4x_1-4x_2)(x_1-x_2)(-x_1+2x_2)}.\]
\end{example}

\begin{example}\label{ex:par}
Suppose that $d_i=\dim(P_i) = 1$ for all $i$.  This is the case of a \emph{parallelotope}, a special case of a zonotope. Up to a translation, we may choose $P_i=[-\p_i/2,\p_i/2]$. Then $\v_i=-\v_{i,1}=\v_{i,2}$ so that $\langle \p_i,\v_j\rangle=\delta_{i,j}$. In other words, $\{\p_i\}_{i=1}^d$ and $\{\v_i\}_{i=1}^d$ are dual basis of $\R^d$. Here, $\kappa(\P)=\det(\v_1,\ldots,\v_d)$. We see that $h_{\x\P}(-\v_i)=\frac{x_i}{2}$ and $h_{\x\P}(\v_i)=\frac{x_i}{2}$. Therefore, by Proposition~\ref{prop:fine-mixed-dual-volume}, \[m_{\P}(\x)=\kappa(\P)\frac{4^d}{x_1x_2\cdots x_d}.\]
If we include translation and consider $\P'=(P_1,\ldots,P_d,-e_1,\ldots,-e_d)$ with \[\x\P'=x_1P_1+\cdots+x_dP_d-y_1e_1-\cdots-y_de_d.\]
Our choices of $\v_i$'s do not change. We now have $h_{\x\P'}(-\v_i)=\frac{x_i}{2}-\langle\y,\v_i\rangle$ and $h_{\x\P'}(\v_i)=\frac{x_i}{2}+\langle\y,\v_i\rangle$. Therefore, we calculate that \[m_{\P'}(\x,\y)=\kappa(\P)\frac{x_1\cdots x_d}{\prod_{i=1}^d(x_i/2-\langle\y,\v_i\rangle)(x_i/2+\langle\y,\v_i\rangle)}=\kappa(\P)\frac{x_1\cdots x_d}{\prod_{i=1}^d(x_i^2/4-\langle\y,\v_i\rangle^2)}.\]
\end{example}

%
%

\section{Mixed subdvisions}\label{sec:mixedsub}


To simplify the presentation, we henceforth focus on the case of polytopes.  The case of unbounded polyhedra does not present an essential difficulty.

\begin{defin}
Let $P$ be a full-dimensional polytope in $\R^d$. A \emph{cell} of $P$ is any full-dimensional convex hull of a subset of the vertices of $P$. Two cells $Q$ and $Q'$ intersect properly if their intersection $Q\cap Q'$ is a proper face of both. A \emph{subdivision} of $P$ is a collection of cells that cover $P$ and intersect properly. 

Now consider a sequence of polytopes  $\P = (P_1, P_2, \dots, P_r)$ in $\R^d$ such that $P := P_1+ P_2 + \dots + P_r$ is full dimensional.  A \emph{Minkowski cell} $\Q$ of $\P$ is a sequence $\Q = (Q_1, Q_2, \dots, Q_r)$ with $Q = Q_1+ Q_2+\dots + Q_r$ such that $\dim(Q) = d$ and each $Q_i$ is a (possibly not-full dimensional) convex hull of some vertices of $P_i$ for each $i$. Let $\Q = (Q_1,\dots, Q_r)$ and $\Q' = (Q'_1,\dots, Q'_r)$ be two Minkowski cells of $\P$. Then they \emph{intersect properly as Minkowski sums} if $Q_i$ intersects properly with $Q'_i$ for each $i = 1,2,\dots, r$.
A Minkowski cell $\Q$ is called \emph{fine} if $\sum_{i=1}^r \dim(Q_i) = d$ and each $Q_i$ is a simplex.  

A collection of Minkowski cells $\mathcal{S} = \{\Q^{(1)} , \Q^{(2)}, \dots, \Q^{(n)}\}$ is a \emph{mixed subdivision} of $P$ if we have $\bigcup_i Q^{(i)} = P$ and for $i \neq j$ we have that $\Q^{(i)}$ and $\Q^{(j)}$ intersect properly as Minkowski sums.  If, in addition, every Minkowski cell in it is fine, then $\mathcal{S}$ is called a \emph{fine mixed subdivision}.
\end{defin}


Notice that we use the $r$-tuple $\Q = (Q_1, \dots, Q_r)$ rather than just the polytope $Q$ which is the Minkowski sum $\sum Q_i$ in the definitions above. We do so to emphasize that mixed subdivisions not only depends on the polytopes $Q$, but also the specific Minkowski summands of $Q$ in the given order. For more discussion on this, see \cite[Remark 1.2]{Santos}.

Mixed subdivisions of $\P$ form a poset under refinement. 

\begin{defin}
Let $\Q = (Q_1, Q_2, \dots, Q_r)$ 
 and $\Q' = (Q'_1, Q'_2, \dots, Q'_r)$ be two Minkowski cells of $\P$. If $Q_i \subseteq Q'_i$ for every $1\leq i\leq r$, then we say that $\Q$ is contained in $\Q'$ and denote it by $\Q \leq \Q'$. Let $\mathcal{S}$ and $\mathcal{S}'$ be two mixed subdivisions of $\P$. We call $\mathcal{S}$ a \emph{refinement} of $\mathcal{S}'$ and denote it by $\mathcal{S} \leq \mathcal{S}'$ if for each $\Q\in \mathcal{S}$, there is a cell $\Q'\in \mathcal{S}'$ such that $\Q\leq \Q'$. 
 The poset of all mixed subdivisions under refinement has a unique maximal element which is the trivial subdivision $\mathcal{S} = \{\P\}$. The minimal elements are the fine mixed subdivisions. 
\end{defin}



\begin{ex}
The following \Cref{fig:fine-mixed-subdivision-two-triangles} shows a fine mixed subdivision $\mathcal{S}=\{\Q^{(1)},\ldots,\Q^{(5)}\}$ of $\P=(P_1, P_2)$, where the position of the origin is omitted. 
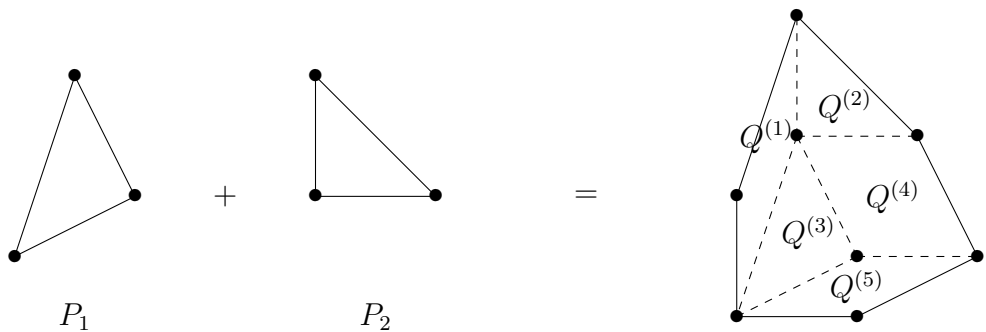
\begin{figure}[h!]
\centering
\begin{tikzpicture}[scale=0.8]
\def\a{-8};
\def\b{-12};
\def\h{-1};
\node at (1+\b,0) {$\bullet$};
\node at (\b,2) {$\bullet$};
\node at (-1+\b,-1) {$\bullet$};
\draw(1+\b,0)--(\b,2)--(-1+\b,-1)--(1+\b,0);
\node at (\b,-2) {$P_1$};
\node at (\a,0) {$\bullet$};
\node at (\a+2,0) {$\bullet$};
\node at (\a,2) {$\bullet$};
\draw(\a,0)--(\a+2,0)--(\a,2)--(\a,0);
\node at (\a+1,-2) {$P_2$};
\node at (0.5+0.5*\a+0.5*\b,0) {$+$};
\node at (0,4+\h) {$\bullet$};
\node at (2,2+\h) {$\bullet$};
\node at (3,\h) {$\bullet$};
\node at (1,-1+\h) {$\bullet$};
\node at (-1,-1+\h) {$\bullet$};
\node at (-1,1+\h) {$\bullet$};
\draw(0,4+\h)--(2,2+\h)--(3,\h)--(1,-1+\h)--(-1,-1+\h)--(-1,1+\h)--(0,4+\h);
\node at (0,2+\h) {$\bullet$};
\node at (1,0+\h) {$\bullet$};
\draw[dashed] (0,2+\h)--(0,4+\h);
\draw[dashed] (0,2+\h)--(2,2+\h);
\draw[dashed] (0,2+\h)--(1,0+\h);
\draw[dashed] (0,2+\h)--(-1,-1+\h);
\draw[dashed] (1,0+\h)--(3,0+\h);
\draw[dashed] (1,0+\h)--(-1,-1+\h);
\node at (0.5*\a+0.5,0) {$=$};
\node at (-0.5,2+\h) {$Q^{(1)}$};
\node at (0.8,2.5+\h) {$Q^{(2)}$};
\node at (0.2,0.4+\h) {$Q^{(3)}$};
\node at (1.6,1+\h) {$Q^{(4)}$};
\node at (1,-0.5+\h) {$Q^{(5)}$};
\end{tikzpicture}
\caption{A fine mixed subdivision of the Minkowski sum of two triangles}
\label{fig:fine-mixed-subdivision-two-triangles}
\end{figure}
\end{ex}

The following is the main theorem in this section.  It generalizes \cref{thm:valuation} and \cref{thm:canformsub} to the mixed setting.

\begin{theorem}\label{thm: mixed subdivision formula}
Let $\mathcal{S} = \{\Q^{(1)} , \Q^{(2)}, \dots, \Q^{(n)}\}$ be a mixed subdivision of $\P = (P_1, \dots, P_r)$, and let $\x = (x_1, x_2, \dots, x_r)$. Then the dual mixed volume of $\P$ can be written as a sums of the dual mixed volumes of the cells in  $\mathcal{S}$:
\begin{equation}\label{eqn: mixed subdivision formula}
m_{\P} (\mathbf{x},\z) = \sum_{ \Q \in \mathcal{S} } m_{\Q} (\mathbf{x},\z).
\end{equation}
\end{theorem}

To prove Theorem \ref{thm: mixed subdivision formula} we use the following known result; see for example \cite{HRS}. 

\begin{lemma}\label{lem:xS}
Suppose that $\x = (x_1, x_2, \dots\, x_r)$ satisfies $x_i >0$ for all $i\in [r]$.  Then $\x\mathcal{S} = \{\x \Q^{(1)} , \x \Q^{(2)}, \dots , \x \Q^{(n)}  \}$ is a subdivision of $\x\P$. 
\end{lemma}

We will also provide proof of \cref{lem:xS} in \cref{sec:Cayley} below.


\begin{proof}[Proof of Theorem \ref{thm: mixed subdivision formula}]
Let $\P = (P_1 \dots, P_r)$ and suppose that $\x = (x_1, x_2, \dots\, x_r)$ satisfies $x_i >0$ for all $i\in [r]$. 
By \cref{lem:xS}, we have an equality of indicator functions
    \[ 
    [\x\P] = \sum_{ \Q \in \mathcal{S}} [\x\Q] + (\text{an alternating sum of indicators of lower dimensional polytopes})
    \]
    By \cref{thm:valuation}, we can take dual volumes of both sides of this equation to obtain \eqref{eqn: mixed subdivision formula}, using the fact that $\VolC(Q) = 0$ for a polytope $Q$ whose dimension is lower than $d$.
\end{proof}
%


\section{Formulae from the Cayley polytope}\label{sec:Cayley}

 \begin{defin}
 Let $\P = (P_1, P_2, \dots, P_r)$ be a sequence of polytopes in $\R^d$. Let $\mu_i: \R^d \to \R^{d+r}$ be the inclusion function such that $\mu_i(\z) = (\z, \mathbf{e}_i)$ for $i \in [r]$, where $\mathbf{e}_1, \mathbf{e}_2, \dots, \mathbf{e}_r$ form the standard basis of $\R^r$. The \emph{Cayley embedding} of $\P$ is the map $\mathcal{C}$ sending the sequence of polytopes $\P=(P_1, P_2, \dots, P_r)$ to following polytope in $\R^{d+r}$.
\[ 
\mathcal{C}(\P):=\mathcal{C}(P_1, P_2, \dots, P_r) := \conv \{ \mu_1(P_1) \cup \mu_2(P_2) \cup \dots \cup \mu_r(P_r) \}.
\]

We call $\mathcal{C}(P_1, P_2, \dots, P_r)$ the \emph{Cayley polytope} of $(P_1, P_2, \dots, P_r)$.

\end{defin}

For full dimensional $\P=(P_1,\ldots,P_r)$, the Cayley polytope $ \mathcal{C}(P_1, P_2, \dots, P_r)$ has codimension $1$ in $\R^{d+r}$, living inside the affine hyperplane $\langle\x,\one\rangle=1$, where we use $\z=(z_1,\ldots,z_d)$ to label the first $d$ coordinates of $\R^{d+r}$ and $\x=(x_1,\ldots, x_r)$ to label the last $r$ coordinates of $\R^{d+r}$. Throughout this section, write $\mathbf{u}=(0,\ldots,0,1,\ldots,1)\in\R^{d+r}$ whose first $d$ coordinates are $0$ and last $r$ coordinates are $1$ to be the normal of the hyperplane where $\mathcal{C}(\P)$ resides. One primary reason to consider the Cayley polytope in our context is that its intersection with $\R^r \times (x_1, x_2, \dots, x_r)$ is exactly $\x\P = x_1P_1 + x_2P_2 + \dots + x_r P_r$ when each $x_i>0$ and $\sum_{i=1}^r x_i = 1$. 

The following proposition is known as the Cayley trick.
\begin{prop}{\cite[Theorem 3.1]{HRS}}\label{prop: Cayley trick}
There is a bijection between the mixed subdivisions of $P_1 + P_2 + \dots + P_r$ and the subdivisions of the Cayley embedding $\mathcal{C}(P_1, P_2, \dots, P_r)$. Furthermore, it restricts to a bijection between the fine mixed subdivisions of $P_1 + P_2 + \dots + P_r$ and the triangulations of $\mathcal{C}(P_1, P_2, \dots, P_r)$. The bijection is given by applying the Cayley embedding on each Minkowski cell in a mixed subdivision $\mathcal{S}$ of $P= P_1 + P_2 + \dots + P_r$. The Cayley polytopes of these cells form a subdivision of $\mathcal{C}(P_1, P_2, \dots, P_r)$.
\end{prop}

\begin{proof}[Proof of  \cref{lem:xS}]
%
By \cref{prop: Cayley trick}, there is a subdivision of the Cayley polytope $\mathcal{C} (P_1, P_2, \dots, P_r)$ that corresponds to $\mathcal{S}$. We call this subdivision $\mathcal{C}(\mathcal{S})$ and the cells in it $\mathcal{C}(\Q^{(1)}) , \mathcal{C}(\Q^{(2)}), \dots , \mathcal{C}(\Q^{(r)})$. Consider the $d$-dimensional affine subspace of $\R^{d+r}$:
\[ W_{\x} := \R^d \times (x_1, x_2, \dots, x_r).\]

The intersection between the Cayley polytope $\mathcal{C}(\P)$ and $W_{\x}$ is $\x\P = x_1 P_1+ x_2 P_2+\dots + x_r P_r$. By the Cayley trick (\Cref{prop: Cayley trick}), every cell of a subdivision of a Cayley polytope is again a Cayley polytope. Therefore for a cell $\mathcal{C}(\Q) \in \mathcal{C}(\mathcal{S})$ where $\Q = (Q_1, Q_2, \dots, Q_r)$, the intersection between $\mathcal{C}(\Q)$ and $W_{\x}$ is $\x\Q$, which is a full dimensional polytope in $W_{\x}$. Since intersecting with an affine subspace does not affect the proper intersection condition of the cells in $\mathcal{C}(\mathcal{S})$, any two cells in $\{\x \Q^{(1)}, \x \Q^{(2)}, \dots, \x \Q^{(n)} \}$ intersect properly as a Minkowski sum. Therefore the subdivision $\mathcal{C}(\mathcal{S})$ of the Cayley polytope $\mathcal{C}(\P)$ produces a subdivision of $\x\P$, which is exactly $\x\mathcal{S}$. \end{proof}

Let $\P= (P_1, P_2, \dots, P_r)$ be a sequence of polytopes in $\R^d$ that is full-dimensional, $\mathcal{C}(\P)$ be the Cayley polytope of $\P$, and $t \mathcal{C}(\P)$ be its dilation by $t\in \R$. The Cayley trick allows us to write the dual mixed volume function as a single dual volume function calculated in an affine hyperplane, where volumes are normalized by the normal vector $\mathbf{u}=(0,\ldots,0,1,\ldots,1)$ (see Section~\ref{sec:affine} for details).
\begin{theorem}\label{thm:Cayleysub}
Let $\P=(P_1,\ldots,P_r)$ be full-dimensional. For $t=\langle\one,\x\rangle$, \[m_{\P}(\x,\z)=\frac{\prod_{i=1}^n x_i}{\sum_{i=1}^n x_i}\cdot\VolC_{(\z,\x)}(t\mathcal{C}(\P))\]
where the right hand side is calculated in the affine hyperplane $\{(\z,\x)\:|\: \langle\one,\x\rangle=t\}\subset\R^{d+r}.$
\end{theorem}
\begin{proof}
By Theorem~\ref{thm: mixed subdivision formula}, we can take a fine mixed subdivision $\mathcal{S}$ of $\P$ so that $m_{\P}(\x,\z)=\sum_{\Q\in\mathcal{S}}m_{\Q}(\x,\z)$. By Proposition~\ref{prop: Cayley trick}, we correspondingly have a triangulation of $\mathcal{C}(\P)$ into $\mathcal{C}(\Q)$'s, and since the dual volume function is a valuation (Theorem~\ref{thm:valuation}), it suffices to prove the current theorem in the case when $\P$ is a fine mixed cell. Correspondingly, $\mathcal{C}(\P)$ is a simplex.

Recall notations from Section~\ref{sec:fine}. Let $d_i=\dim(P_i)$ such that $\sum_{i=1}^r d_i=d$. For $i\in[r]$ and $a\in[d_i+1]$, we define $\v_{i,a}\in\R^d$ so that $\langle\p_{i,a},\v_{i,a}\rangle-\langle\p_{i,a},\v_{i,b}\rangle=1$ for $b\neq a$, and that $\v_{i,a}$ is normal to the span of $P_j$ for all $j\neq i$, where $\p_{i,1},\ldots,\p_{i,d_i+1}$ are the vertices of $P_i$. When $d_i=0$, we simply write $\v_{i,1}=\mathbf{0}$. Now we extend $\v_{i,a}$'s to rays in $\N(\mathcal{C}(\P))$.

For $i\in[r]$ and $a\in[d_i+1]$, we define $\mathbf{u}_{i,a}\in\R^r$ by specifying its coordinates as follows: 
\begin{align*}
\langle\mathbf{u}_{i,a},\e_i\rangle=&-\langle\v_{i,a},\p_{i,b}\rangle-c_{i,a},\text{ for }b\neq a,\\
\langle\mathbf{u}_{i,a},\e_j\rangle=&-\langle\v_{i,a},\p_{j,b}\rangle-c_{i,a},\text{ for }j\neq i,\\
\langle\mathbf{u}_{i,a},\one\rangle=&0,
\end{align*}
for some $c_{i,a}\in\R$, chosen just to make sure the last condition $\langle\mathbf{u}_{i,a},\one\rangle=0$ holds. Note that construction works when $d_i=0$ by the convention that $\v_{i,a}=\mathbf{0}$. Write $\tilde\v_{i,a}=(\v_{i,a},\mathbf{u}_{i,a})\in\R^{d+r}$. It is then clear by construction that $\langle\tilde\v_{i,a},(\p_{i,a},\e_i)\rangle=1-c_{i,a}$ and that $\langle\tilde\v_{i,a},\p\rangle=-c_{i,a}$ for all other vertices $\p$ of the Cayley polytope $\mathcal{C}(\P)$. Therefore, $\tilde\v_{i,a}$'s are all the extremal rays in the normal fan $\N(\mathcal{C}(\P))$. Note that since $\mathcal{C}(\P)$ is a simplex, its normal fan has $d+r$ maximal cones, generated by all $(d+r-1)$-subsets of $\tilde\v_{i,a}$'s.

Fix $i_0\in[r]$ and $a_0\in[d_i+1]$. Let's calculate the determinant formed by all $\tilde\v_{i,a}$'s except $\tilde\v_{i_0,a_0}$, together with $\mathbf{u}=(0,\ldots,0,1,\ldots,1)$. Recall from Lemma~\ref{lem:fine-mixed-cell-normal-fan} that $\sum_{a=1}^{d_i+1}\v_{i,a}=\textbf{0}$. For $i\neq i_0$, let $\mathbf{u}_i=\sum_{a=1}^{d_i+1}\mathbf{u}_{i,a}$. We calculate that \[\langle\mathbf{u}_i,\e_j\rangle=-\left\langle\sum_{a=1}^{d_i+1}\v_{i,a},\p_{i,b}\right\rangle-\sum_{a=1}^{d_i+1}c_{i,a}=-\sum_{a=1}^{d_i+1}c_{i,a}\text{ for }j\neq i\]
and that \[\langle\mathbf{u}_i,\e_i\rangle=1-\langle\v_{i,1},\p_{i,1}\rangle-\left\langle\sum_{a=2}^{d_i+1}\v_{i,a},\p_{i,1}\right\rangle-\sum_{a=1}^{d_i+1}c_{i,a}=1-\sum_{a=1}^{d_i+1}c_{i,a}\]
as $\langle\v_{i,a},\p_{i,a}\rangle=1+\langle\v_{i,a},\p_{i,b}\rangle$ for $b\neq a$. In the determinant calculation, for each $i\neq i_0$, we can thus replace one of $\tilde\v_{i,1},\ldots,\tilde\v_{i,d_i+1}$ by $\sum_{a=1}^{d_i+1}\tilde\v_{i,a}+\sum_{a=1}^{d_i+1}c_{i,a}\mathbf{u}=\e_i$, without changing its absolute value. After that, we can take away the row and column corresponding to these $\e_i$'s. Next, $\mathbf{u}$ becomes $\e_{i_0}$, which can be taken away as well. We are then left with a square matrix of size $d$, whose determinant after taking the absolute value is precisely $\kappa(\P)$ as in Definition~\ref{def:kappa-determinant}.

For $d_i>0$, we examine $\min_{\y\in\mathcal{C}(\P)}\langle\y,(\v_{i,a},\mathbf{u}_{i,a})\rangle$. It suffices to consider vertices $\y=(\p_{j,b},\e_j)$ of $\mathcal{C}(\P)$ for some $j\in[r]$ and $b\in[d_j+1]$. By construction, if $j\neq i$, $\langle\y,\tilde\v_{i,a}\rangle=\langle\v_{i,a},\p_{j,b}\rangle+\langle\mathbf{u}_{i,a},\e_j\rangle=-c_{i,a}$. Similarly, if $j=i$ and $b\neq a$, we also have $\langle\y,\tilde\v_{i,a}\rangle=-c_{i,a}$ and if $j=i$ and $b=a$, we have $\langle\y,\tilde\v_{i,a}\rangle=1-c_{i,a}$. The minimum is $-c_{i,a}$. Now,
\begin{align*}
h_{t\mathcal{C}(\P)-(\z,\x)}(\tilde\v_{i,a})=&-t\min_{\y\in\mathcal{C}(\P)}\langle\y,\tilde\v_{i,a}\rangle+\langle\v_{i,a},\z\rangle+\langle\mathbf{u}_{i,a},\x\rangle\\
=&tc_{i,a}+\langle\v_{i,a},\z\rangle+\big(-c_{i,a}-\langle\v_{i,a},\p_{i,b}\rangle\big)x_i+\sum_{j\neq i}\big(-c_{i,a}-\langle\v_{i,a},\p_{j,1}\rangle\big)x_j\\
=&-\langle\v_{i,a},\p_{i,b}\rangle x_i+\sum_{j\neq i}-\langle\v_{i,a},\p_{j,1}\rangle x_j+\langle\v_{i,a},\z\rangle\\
=&\sum_{j=1}^r x_j h_{P_j}(\v_{i,a})+\langle\v_{i,a},\z\rangle=h_{\x\P-\z}(\v_{i,a})
\end{align*}
where $b\neq a$ and $t$ is substituted to be $x_1+\cdots+x_r$. For $d_i=0$, $\v_{i,1}=\mathbf{0}$ and $\tilde\v_{i,1}=\e_i$. Then $h_{t\mathcal{C}(\P)-(\z,\x)}(\tilde\v_{i,1})=x_i$. Also recall from Equation~\eqref{eqn:sum-of-hs-in-cell} that $\sum_{a=1}^{d_i+1}h_{\x\P-\z}(\v_{i,a})=x_i$.

Finally, putting all the pieces together, we obtain that 
\begin{align*}
\VolC_{(\z,\x)}(t\mathcal{C}(\P))=&\sum_{i=1}^r\sum_{a=1}^{d_i+1}\frac{\kappa(\P)}{\prod_{(j,b)\neq(i,a)}h_{t\mathcal{C}(\P)-(\z,\x)}(\tilde\v_{j,b})}=\kappa(\P)\frac{\sum_{i=1}^r\sum_{a=1}^{d_i+1}h_{t\mathcal{C}(\P)-(\z,\x)}(\tilde\v_{i,a})}{\prod_{i=1}^r\prod_{a=1}^{d_i+1}h_{t\mathcal{C}(\P)-(\z,\x)}(\tilde\v_{i,a})}\\
=&\kappa(\P)\frac{\sum_{i=1}^r x_i}{\prod_{d_i=0}x_i\prod_{d_i>0}\prod_{a=1}^{d_i+1}h_{\x\P-\z}(\v_{i,a})}=\frac{\sum_{i=1}^r x_i}{\prod_{i=1}^r x_i}\frac{\prod_{d_i>0}x_i}{\prod_{d_i>0}\prod_{a=1}^{d_i+1}h_{\x\P-\z}(\v_{i,a})}\\
=&\frac{\sum_{i=1}^r x_i}{\prod_{i=1}^r x_i} m_{\P}(\x,\z)
\end{align*}
by Proposition~\ref{prop:fine-mixed-dual-volume}, as desired.
\end{proof}

\begin{ex}
Let $P_1=\conv((0,0),(1,0))$ and $P_2=\conv((0,0),(1,2))$ in $\R^2$. We can take the extremal rays in $\N(P)$, where $P=P_1+P_2$, to be \[\v_1=\left(-1,\frac{1}{2}\right),\ \v_2=\left(1,-\frac{1}{2}\right),\ \v_3=\left(0,-\frac{1}{2}\right),\ \v_4=\left(0,\frac{1}{2}\right),\]
with $\kappa(\P)=|\det(\v_1,\v_3)|=\frac{1}{2}$. We then see that 
\begin{align*}
h_{\x\P-\z}(\v_1)=x_1-z_1+\frac{z_2}{2},\quad& h_{\x\P-\z}(\v_2)=z_1-\frac{z_2}{2},\\
h_{\x\P-\z}(\v_3)=x_2-\frac{z_2}{2},\quad& h_{\x\P-\z}(\v_4)=\frac{z_2}{2}.
\end{align*}
Proposition~\ref{prop:fine-mixed-dual-volume} yields \[m_{\P}(\x,\z)=\frac{1}{2}\frac{x_1x_2}{(x_1-z_1+\frac{z_2}{2})(z_1-\frac{z_2}{2})(x_2-\frac{z_2}{2})(\frac{z_2}{2})}.\]

On the other side, $\mathcal{C}(\P)=\conv\big((0,0,1,0),(1,0,1,0),(0,0,0,1),(1,2,0,1)\big).$ We can extend $\v_1,\ldots,\v_4$ to be extremal rays of the normal fan $\N(\mathcal{C}(\P))$ as follows:
\[\tilde\v_1=\left({-}1,\frac{1}{2},\frac{1}{2},{-}\frac{1}{2}\right),\ \tilde\v_2=\left(1,{-}\frac{1}{2},0,0\right),\ \tilde\v_3=\left(0,{-}\frac{1}{2},{-}\frac{1}{2},\frac{1}{2}\right),\ \tilde\v_4=\left(0,\frac{1}{2},0,0\right).\]
We calculate that $|\det(\v_i,\v_j,\v_k,(0,0,1,1))|=\frac{1}{2}=\kappa(\P)$ for $1\leq i<j<k\leq 4$. Also,
\begin{align*}
h_{t\mathcal{C}(\P)-(\z,\x)}(\v_1)=\frac{t}{2}-z_1+\frac{z_2}{2}+\frac{x_1}{2}-\frac{x_2}{2},\quad& h_{t\mathcal{C}(\P)-(\z,\x)}(\v_2)=z_1-\frac{z_2}{2},\\
h_{t\mathcal{C}(\P)-(\z,\x)}(\v_3)=\frac{t}{2}-\frac{z_2}{2}-\frac{x_1}{2}+\frac{x_2}{2},\quad& h_{t\mathcal{C}(\P)-(\z,\x)}(\v_4)=\frac{z_2}{2}.
\end{align*}
These four sum to $t$. The dual mixed volume function of $t\mathcal{C}(\P)$ is then \[\VolC_{(\z,\x)}(t\mathcal{C}(\P))=\frac{1}{2}\frac{t}{(\frac{t}{2}-z_1+\frac{z_2}{2}+\frac{x_1}{2}-\frac{x_2}{2})(z_1-\frac{z_2}{2})(\frac{t}{2}-\frac{z_2}{2}-\frac{x_1}{2}+\frac{x_2}{2})(\frac{z_2}{2})}\]
which differs from $m_{\P}(\x,\z)$ by a factor of $x_1x_2/(x_1+x_2)$ when $t=x_1+x_2$.
\end{ex}

\section{Polytopes in an affine hyperplane}\label{sec:affine}
Many polytopes of interest in combinatorics, such as matroid polytopes, alcoved polytopes, and generalized permutohedra, are naturally defined to live in an affine hyperplane.
\subsection{Normalization in an affine hyperplane}\label{sub:normalization-affine}
For an affine hyperplane $H\subset\R^d$ normal to $\v$ and a point $\z_0\in H$, we identify $H$ with $\R^{d-1}$, with $\z_0$ as the origin as follows. Its vector space structure is denoted as $(+_H,\cdot_H)$, which is given by
$$
\alpha \cdot_H \y := \alpha(\y - \z_0) + \z_0 = \alpha \y + (1-\alpha) \z_0, \quad \y +_H \y' := (\y - \z_0) + (\y' - \z_0) + \z_0 = \y + \y' - \z_0,
$$
for $\alpha \in \R$, and $\y,\y' \in H$. Its dual space $H^{*}$ can be identified with the space of linear functions that vanish on $\z_0$:
\[H^* = \{f: \R^d \to \R \mid f(\z_0) = 0\}.\]
Indeed, when $f(\z_0) =0$, we have
\[f(\alpha \cdot_H y) = f( \alpha \y + (1-\alpha) \z_0) = \alpha f(\y), \quad f(\y +_H \y') =  f(\y + \y' - \z_0) = f(\y) + f(\y').\]

We normalize the volume in the dual space $\R^d/\v$ so that the simplex spanned by vectors $\v_1,\ldots,\v_{d-1}$ has volume $|\det(\v_1,\ldots,\v_{d-1},\v)|$, which is a quantity that does not depend on representatives in $\R^d/\v$. This definition depends on scaling of $\v$ but for our purposes, choices of $\v$ will be fairly natural. For the rest of this section, this normal vector is chosen to be $\v = \one$ and we consider affine hyperplanes of the following form:
$$
H_k:= \{\y = (y_1,\ldots,y_d) \mid y_1 + y_2 + \cdots + y_d = k\}.
$$

Let $P \subseteq H_k$ be full-dimensional in $H_k$.  Let $\Cone(P)$ be the \emph{cone over $P$}, 
$$
\Cone(P):= \{t \cdot \y \mid \y \in P, \;\; t \in \R_{\geq 0}\}.
$$

The normal fan of $\N(P)$ is a complete fan in $\R^d$, where every cone $C \in \N(P)$ has lineality space $\R \cdot \one$, where $\one = (1,1,\ldots,1)$.  The reduced normal fan $\N'(P)$ is the image of $\N(P)$ in $\R^d/\one$, and every cone of $\N'(P)$ is a pointed, convex, polyhedral cone.

Let $\Cone(P)^* \subseteq \R^d$ be the dual cone of $\Cone(P)$.  Then each cone $C \in \N(P)$ is the direct product of $\R \cdot \one$ with a face of $\Cone(P)^*$.

\subsection{Dual volume in an affine hyperplane}
The following \emph{hyperplane dual volume function} will play the role of the \emph{dual volume function} for full dimensional polytopes that live inside the affine hyperplane $H_k$, including all important polytopes that appear in later sections.
\begin{defin}\label{def:hyperplane-dual-volume}
Fix a triangulation $\T$ of the boundary of $\Cone(P)^*$, and let $\{\v_i\}_{i=1}^g$ be generators for the one-dimensional cones that appear in $\T$. We write $F \in \T$ for a top-dimensional ($(d-1)$-dimensional) simplicial cone $F$ in this triangulation. The \emph{hyperplane dual volume function} of the ``affine polytope" $P$ is 
\begin{equation}\label{eq:affinedual}
\EVolC_\z(P) := \sum_{F=\mathrm{span}(\v_{i_1},\ldots,\v_{i_{d-1}})\in\T}\frac{|\det(\v_{i_1},
\v_{i_2},\ldots,\v_{i_{d-1}},\one)|}{\langle \v_{i_1},\z \rangle \langle \v_{i_2},\z \rangle \cdots \langle \v_{i_{d-1}},\z \rangle}.
\end{equation}
If $P$ does not have full dimension in $H_k$, i.e. $\dim(P)<d-1$, we define $\EVolC_\z(P) :=  0$.
\end{defin}
We reiterate that $\T$ is a triangulation of the boundary of $\Cone(P)^*$, not of $\Cone(P)^*$ itself. Also note that $\EVolC_\z(P)$ depends only on $\Cone(P)$; it is unchanged under replacing $P$ by a scalar multiple $tP$. Similar to the case of the dual volume (Definition~\ref{def:dual-volume}), we will see that $\EVolC_{\z}(P)$ does not depend on the triangulation $\T$. 

\begin{proposition}\label{prop:enrich}
The hyperplane dual volume function $\EVolC_\z(P)$ is a rational function of degree $1-d$ that does not depend on the choice of $\T$. Moreover, we can realize the hyperplane dual volume function as a dual volume function in $H_t$ with $\z_0\in\R^d$ as the origin (see Section~\ref{sub:normalization-affine} for the identification) where $t:=\langle\one,\z_0\rangle$. Then
\begin{equation}\label{eq:hyperplane-function-is-actual-dual-function}
\EVolC_{\z}(P)=\Vol^\vee_\z((t/k)P) \mbox{ for $\z \in H_t$},
\end{equation}
and
if the dual volume $\Vol^\vee((t/k)P)$ is well-defined, that is, $\z_0$ is not in the affine span of the facets of $(t/k)P$, we have
\begin{equation}\label{eq:hyperplane-volume-is-actual-dual-volume}
\EVolC_{\z}(P)\mid_{\z=\z_0}=\VolC((t/k)P).
\end{equation}
\end{proposition}
\begin{proof}
Given $\v\in\R^d$ viewed as a linear function on $\R^d$, define \[\tilde\v:=\v-\frac{\langle\v,\z_0\rangle}{t}\one.\]
Then $\langle\tilde\v,\y\rangle=\langle\v,\y\rangle-\langle\v,\z_0\rangle$ for $\y\in H_t$, and in particular, $\langle\tilde\v,\z_0\rangle=0$. With this notation, $\tilde\v$ is the unique linear function in $H_t^*$ with the same image as $\v$ in the quotient $\R^d/\one$. Therefore, for a fixed a triangulation $\T$ of the boundary of $\Cone(P)^{*}$, we naturally obtain a triangulation $\tilde T$ of $\N((t/k)P)$ in $H_t^*$ via the map $\v\mapsto\tilde\v$ on rays. We calculate the support function. By definition of $\Cone(P)^{*}$, $\langle\v,\y\rangle\geq0$ for all $\y\in P$ and there exists a point $\p\in P$ such that $\langle\v,\p\rangle=0$. Thus,
\[h_{(t/k)P}(\tilde\v)=-\min_{\y\in\P}\langle\tilde\v,\y\rangle=-\min_{\y\in\P}\langle\v,\y\rangle+\langle\v,\z_0\rangle=\langle\v,\z_0\rangle\]
and for $\z\in H_t$, we then have that \[h_{(t/k)P-\z}(\tilde\v)=h_{(t/k)P}(\tilde\v)+\langle\tilde\v,\z\rangle=\langle\v,\z_0\rangle+\langle\v,\z\rangle-\langle\v,\z_0\rangle=\langle\v,\z\rangle.\]

With Definition~\ref{def:dual-volume}, for $\z\in H_t$, we can write
\begin{align*}
\VolC_{\z}((t/k)P)=&\sum_{\tilde F=\mathrm{span}(\tilde\v_{i_1},\ldots,\tilde\v_{i_{d-1}})\in\tilde T}\frac{|\det(\tilde\v_{i_1},\tilde\v_{i_2},\ldots,\tilde\v_{i_{d-1}},\one)|}{h_{(t/k)(P)-\z}(\tilde\v_{i_1})h_{(t/k)(P)-\z}(\tilde\v_{i_2})\cdots h_{(t/k)(P)-\z}(\tilde\v_{i_{d-1}})}\\
=&\sum_{F=\mathrm{span}(\v_{i_1},\ldots,\v_{i_{d-1}})\in\T}\frac{|\det(\v_{i_1},
\v_{i_2},\ldots,\v_{i_{d-1}},\one)|}{\langle \v_{i_1},\z \rangle \langle \v_{i_2},\z \rangle \cdots \langle \v_{i_{d-1}},\z \rangle}=\EVolC(P).
\end{align*}
where the equality on the numerator is due to the fact that $\tilde\v$ and $\v$ differ by a multiple of $\one$ so the determinant remains unchanged. This shows Equation~\eqref{eq:hyperplane-function-is-actual-dual-function}, which implies Equation~\eqref{eq:hyperplane-volume-is-actual-dual-volume}, and also implies that $\EVolC$ does not depend on the triangulation $\T$.
\end{proof}

\begin{example}
Let $P$ be the simplex that is the convex hull of $(1,0,0),(0,1,0),(0,0,1) \in \R^3$, which lies in $H_1$. So the cone $\Cone(P)$ is the nonnegative octant of $\R^3$ and we can identify the dual cone $\Cone(P)^*$ as the nonnegative octant of $\R^3$ as well. Take rays $\v_1=(1,0,0)$, $\v_2=(0,1,0)$, $\v_3=(0,0,1)$.  The facets of $\Cone(P)^*$ are \[\sp_{\geq 0}(r_1,r_2), \qquad \sp_{\geq 0}(r_1,r_3), \qquad \sp_{\geq 0}(r_2,r_3),\]
each of which is a simplicial cone of dimension $2$ already. Then \[\EVolC_{\z}(P)=\frac{|\det(\v_1,\v_2,\one)|}{\langle\v_1,\z\rangle\langle\v_2,\z\rangle}+\frac{|\det(\v_1,\v_3,\one)|}{\langle\v_1,\z\rangle\langle\v_3,\z\rangle}+\frac{|\det(\v_2,\v_3,\one)|}{\langle\v_2,\z\rangle\langle\v_3,\z\rangle}=\frac{1}{z_1z_2}+\frac{1}{z_1z_3}+\frac{1}{z_2z_3}.\]
Now let us suppose that we take $\z_0 = (0,0,1) = \v_3$.  Then $H_1$ can be identified with $\R^2$ by the projection $(y_1,y_2,y_3) \mapsto (y_1,y_2)$, under which $P$ is identified with the convex hull of $(0,0), (1,0),(0,1)$ in $\R^2$.  Then one can calculate
\[\Vol^\vee_\y(P) = \frac{1}{y_1y_2} + \frac{1}{y_1(1-y_1-y_2)} + \frac{1}{y_2(1-y_1-y_2)},\]
which is equal to $\EVolC_\z(P)$ under the substitution $(z_1,z_2,z_3) = (y_1,y_2,1-y_1-y_2)$ for a typical point in $H_1$.
\end{example}

\begin{thm}\label{thm:affinevaluation}
Suppose that we have the identity
$$
\sum_{i=1}^r \alpha_i [P_i] = 0
$$
for polytopes $P_1,\ldots,P_r\in H_k$ and $\alpha_1,\ldots,\alpha_r\in\R$, then
$$
\sum_{i=1}^r \alpha_i \EVolC_\z(P_i) = 0.
$$
\end{thm}
\begin{proof}
Let $t=\langle\z,\one\rangle$ and by Proposition~\ref{prop:enrich}, it suffices to show that $\sum_i \alpha_i\VolC_{\z}((t/k)P_i)=0$. This follows from Theorem~\ref{thm:valuation} applied to $\sum_{i=1}^r\alpha_i[(t/k)P_i]=0$ in $H_t$ because $\sum_{i=1}^r\alpha_i[P_i]=0$ in $H_k$.
\end{proof}

\subsection{Dual mixed volume in an affine hyperplane}\label{sec:mixedaffine}
Let $P_1,\ldots,P_r \subseteq H_1$ and $P:=P_1+\cdots+P_r$ be full-dimensional in $H_r$. As before, we fix a triangulation $\T$ of the boundary of $\Cone(P)^*$, and let $\{\v_i\}_{i=1}^g$ be generators for the one-dimensional rays that appear in $\T$. As before, write $\x\P$ for $x_1P_1+\cdots+x_rP_r$, and $h_{\x\P}(\v)=\sum_{i=1}^r x_ih_{P_i}(\v)$. 

\begin{defn}\label{def:hyperplane-dual-mixed}
The \emph{hyperplane dual mixed volume function} of $\P=(P_1,\ldots,P_r)$ is
\begin{equation}\label{eq:affinedualmixed}
\tilde m_{\P}(\x,\z) := \sum_{F=\mathrm{span}(\v_{i_1},\ldots,\v_{i_{d-1}})\in\T}\frac{|\det(\v_{i_1},\v_{i_2},\ldots,\v_{i_{d-1}},\one)|}{h_{\x\P-\z}(\v_{i_1})h_{\x\P-\z}(\v_{i_2}) \cdots h_{\x\P-\z}(\v_{i_{d-1}})},
\end{equation}
a rational function in $\x=(x_1,\ldots,x_r)$ and $\z=(z_1,\ldots,z_d)$. As before, \[h_{\x\P}(\v)=\sum_{i=1}^r x_ih_{P_i}(\v)\text{ and }h_{\x\P-\z}(\v)=h_{\x\P}(\v)+\langle\v,\z\rangle.\]
\end{defn}


The following fact follows directly from Proposition~\ref{prop:mpxz-regular}, where $\x\P$ is considered as a polytope in $H_t$ with $t=\langle\z_0,\one\rangle=\langle\x,\one\rangle$.
\begin{proposition}\label{prop:hyperplane-mixed}
The hyperplane dual mixed volume function $\tilde m_{\P}(\x,\z)$ is a rational function of degree $1-d$ that does not depend on the choice of $\T$.  For $\z_0 \in \R^d$ satisfying $\langle\z_0,\one\rangle=\langle\x,\one\rangle=t$, we have $\tilde m_{\P}(\x,\z)|_{\z=\z_0}=\VolC_{\z_0}(\x\P)$, where $\x\P$ is viewed as living in $H_t$, identified as a vector space with origin $\z_0$.
\end{proposition}

\subsection{Fine mixed cells in an affine hyperplane}
Our setup in this section resembles that of Section~\ref{sec:fine}, but there are some differences. We will highlight the distinction.
\begin{defin}
A \emph{fine mixed cell in the affine hyperplane $H_1$} is a sequence $\P=(P_1,\ldots,P_r)$ of simplices in $H_1$ satisfying both of the conditions
\begin{itemize}
\item $\dim(P_1)+\cdots+\dim(P_r)=d-1$, and
\item $\dim(P_1+P_2+\cdots+P_r)=d-1$.
\end{itemize}
\end{defin}
Let $\P$ be a fine mixed cell and write $d_i=\dim(P_i)$ for $i=1,\ldots,r$ as before, and write $P=P_1+\cdots+P_r$. We describe $\Cone(P)^*$ as follows. Label vertices of $P_i$ as $\p_{i,1},\ldots,\p_{i,d_i+1}$ and for each $P_i$ with $d_i>0$, let $\v_{i,1},\ldots,\v_{i,d_{i}+1}$ be rays in $\R^d$ to be normal to $\mathrm{span}(P_j)$ for each $j\neq i$, and for each $a\in[d_i+1]$, 
\begin{align*}\label{eq:affine-fine-mixed-cell-rays}
&\max_{\y\in P_i}\langle\y,\v_{i,a}\rangle=\langle\p_{i,a},\v_{i,a}\rangle=1-\sum_{j\neq i}\langle \z_j,\v_{i,a}\rangle,\\
&\min_{\y\in P_i}\langle\y,\v_{i,a}\rangle=\langle\p_{i,b},\v_{i,a}\rangle=-\sum_{j\neq i}\langle\z_j,\v_{i,a}\rangle,
\end{align*}
for any $\z_j\in P_j$ and any $b\neq a$. The first equality signs determine the direction of $\v_{i,a}$ and the second equality sign tells us how to scale it. By counting dimensions, we see that $\v_{i,a}$'s are uniquely defined. As before, if $d_i=\dim(P_i)=0$, we just assume $\v_{i,1}=\frac{1}{r}\one$. This vector will not appear in any meaningful calculation. We do this scaling to ensure that $\langle\v_{i,1},\y\rangle=1$ for all $\y\in P=P_1+\cdots+P_r$. The idea of how $\v_{i,a}$'s are scaled is the same as in Section~\ref{sec:fine}, but the technical details are different.
\begin{lemma}\label{lem:affine-fine-mixed-cell-normal-fan}
With notations as above, $\sum_{a=1}^{d_i+1}\v_{i,a}=\frac{1}{r}\one$. Moreover, the boundary of $\Cone(P)^*$ has maximal cones $\{C_{\mathbf{a}}\:|\:\mathbf{a}=(a_1,\ldots,a_r)\in[d_1{+}1]\times\cdots\times[d_r{+}1]\}$ where \[C_{\mathbf{a}}=C_{a_1,\ldots,a_r}=\mathrm{span}_{\geq0}\bigcup_{i=1}^r\{\v_{i,a}\:|\:a\neq a_i\}.\]
\end{lemma}
\begin{proof}
By construction of $\v_{i,a}$'s, $\min_{\y\in P}\langle\y,\v_{i,a}\rangle=0$ and all the $\y$'s that satisfy $\langle\y,\v_{i,a}\rangle=0$ can be written as $\z_1+\cdots+\z_r$ where $\z_j\in P_j$ is arbitrary and $\z_i$ lies in the facet of $P_i$ which does not contain $\p_{i,a}$. This means that $\v_{i,a}$'s are precisely the rays in $\Cone(P)^*$. As $P$ is simple, maximal cones in the boundary of $\Cone(P)^*$ must be simplicial and are in bijection with vertices of $P$. All vertices of $P$ are indexed by $\a=(a_1,\ldots,a_r)\in[d_1+1]\times\cdots\times[d_r+1]$ where $\p_{\a}=\p_{1,a_1}+\cdots+\p_{r,a_r}$. As each $\v_{i,a}$ with $a\neq a_i$ satisfies that $\v_{i,a}\in\Cone(P)^*$ and $\langle\v_{i,a},\p_{\a}\rangle=0$, they are generators of the maximal cone in the boundary of $\Cone(P)^*$ corresponding to $\p_{\a}$, and by counting the total number $d_1+\cdots+d_r=d-1$, this cone is precisely $C_{\a}$ as defined above. This completes the description of $\Cone(P)^*$.

We have seen that $\langle\v_{i,a},\p_{1,a_1}+\cdots+\p_{r,a_r}\rangle$ equals $0$ if $a\neq a_i$ and equals $1$ if $a=a_i$. For $\v_i:=\v_{i,1}+\cdots+\v_{i,d_i+1}$, we then have that $\langle\v_i,\p_{\a}\rangle=1$ for all vertices $\p_{\a}$ of $P$, and thus $\langle\v_i,\y\rangle=1$ for all $\y\in P$. This means that $\v_i$ is normal to $H_1$ and can be written as $\gamma\one$ for some $\gamma$. As $P=P_1+\cdots+P_r\in H_r$ and $\langle\gamma\one,\y\rangle=1$ for all $\y\in P$, we see $\gamma=\frac{1}{r}$.
\end{proof}
As before, we define $\kappa(\P)$ to be the absolute value of the determinant with $d$ vectors $\v_{i,a}$'s together with $\one$, where $i\in[r]$ and $a\in A_i$ for any subset $A_i\subsetneq [d_i+1]$ with cardinality $d_i$. By Lemma~\ref{lem:affine-fine-mixed-cell-normal-fan}, as $\v_{i,1}+\cdots+\v_{i,d_i+1}$ is a multiple of $\one$, $\kappa(\P)$ is well-defined.
\begin{proposition}\label{prop:{lemma}affine-fine-mixed-cell-dual-volume}
Let $\P=(P_1,\ldots,P_r)$ be a fine mixed cell in $H_1$ and let $\v_{i,a}$ be as above. Then its hyperplane dual volume function is \[\tilde m_{\P}(\x,\z)=\kappa(\P)\prod_{i:d_i>0}\frac{x_i-\langle\x,\one\rangle/r+\langle\z,\one\rangle/r}{\prod_{a=1}^{d_i+1}(h_{\x\P}(\v_{i,a})+\langle\v_{i,a},\z\rangle)}\]
In particular, when $\langle\x,\one\rangle=\langle\z,\one\rangle$, this formula simplifies to \[\tilde m_{\P}(\x,\z)|_{\langle\x,\one\rangle=\langle\z,\one\rangle}=\kappa(\P)\frac{\prod_{i:d_i>0} x_i}{\prod_{i:d_i>0}\prod_{a=1}^{d_i+1}(h_{\x\P}(\v_{i,a})+\langle\v_{i,a},\z\rangle)}\]
\end{proposition}
\begin{proof}
Let $\p_{i,a}$ be vertices of $P_i\subsetneq H_1$ for $i\in[r]$ and $a\in[d_i+1]$. Using Lemma~\ref{lem:affine-fine-mixed-cell-normal-fan}, we calculate that for each $d_i>0$,
\begin{align*}
\sum_{a=1}^{d_i+1}h_{\x\P}(\v_{i,a})=&\sum_{j\neq i} x_j\left(\sum_{a=1}^{d_i+1}h_{P_j}(\v_{i,a})\right)+x_i\sum_{a=1}^{d_i+1}h_{P_i}(\v_{i,1})\\
=&\sum_{j\neq i}x_j\left(-\langle \y_j,\v_{i,a}\rangle\right)+x_i\left(-\langle\p_{i,2},\v_{i,1}\rangle-\sum_{a=2}^{d_i+1}\langle\p_{i,1},\v_{i,a}\rangle\right)\\
=&-\sum_{j\neq i}x_j\langle\y_j,\frac{1}{r}\one\rangle+x_i\left(-\langle\p_{i,2},\v_{i,1}\rangle-\left\langle\p_{i,1},\frac{1}{r}\one-\v_{i,1}\right\rangle\right)\\
=&-\frac{1}{r}\sum_{j\neq i}x_j+x_i\left(-\frac{1}{r}+\langle\p_{i,1},\v_{i,1}\rangle-\langle\p_{i,2},\v_{i,1}\rangle\right)\\
=&x_i-\frac{1}{r}\langle\x,\one\rangle.
\end{align*}
Taking $\z$ into account, we have
\[\sum_{a=1}^{d_i+1}h_{\x\P-\z}(\v_{i,a})=\sum_{a=1}^{d_i+1}h_{\x\P}(\v_{i,1})+\langle\z,\sum_{a=1}^{d_i+1}\v_{i,a}\rangle=x_i-\frac{1}{r}\langle\x,\one\rangle+\frac{1}{r}\langle\z,\one\rangle.\]
We then calculate $\tilde m_{\P}(\x,\z)$ by Definition~\ref{def:hyperplane-dual-volume}:
\begin{align*}
\tilde m_{\P}(\x,\z)=&\sum_{(a_1,\ldots,a_r)\in[d_1+1]\times\cdots\times[d_r+1]}\kappa(\P)\prod_{i:d_i>0}\prod_{a\neq a_i}\frac{1}{h_{\x\P-\z}(\v_{i,a})}\\
=&\kappa(\P)\prod_{i:d_i>0}\frac{\sum_{a=1}^{d_i+1}h_{\x\P-\z}(\v_{i,a})}{\prod_{a=1}^{d_i+1}h_{\x\P-\z}(\v_{i,a})}\\
=&\kappa(\P)\prod_{i:d_i>0}\frac{x_i-\langle\x,\one\rangle/r+\langle\z,\one\rangle/r}{\prod_{a=1}^{d_i+1}(h_{\x\P}(\v_{i,a})+\langle\v_{i,a},\z\rangle)}.
\end{align*}
When $\langle\x,\one\rangle=\langle\z,\one\rangle$, this expression simplifies as in the statement. 
\end{proof}

\begin{example}
Let $P_1=\conv\big(\p_{1,1}=(1,0,0),\ \p_{1,2}=(0,1,0)\big)$, $P_2=\conv\big(\p_{2,1}=(0,1,0),\ \p_2=(0,0,1)\big)$, $P_3=\p_{3,1}=(1,0,0)$ be two two line segments and a point in $H_1\subset\R^3$. Then $P=P_1+P_2+P_3$ is a parallelogram and $\Cone^*(P)$ is a cone with $4$ rays. To calculate $\v_{1,1}$, we need the following equations:
\begin{align*}
&\langle\v_{1,1},\p_{1,1}+\p_{2,1}+\p_{3,1}\rangle=\langle\v_{1,1},\p_{1,1}+\p_{2,2}+\p_{3,1}\rangle=1,\\
&\langle\v_{1,1},\p_{1,2}+\p_{2,1}+\p_{3,1}\rangle=\langle\v_{1,1},\p_{1,2}+\p_{2,2}+\p_{3,1}\rangle=0,\\
\end{align*}
and we calculate that $\v_{1,1}=(\frac{2}{3},-\frac{1}{3},-\frac{1}{3})$. Similarly, $\v_{1,2}=(-\frac{1}{3},\frac{2}{3},\frac{2}{3})$, $\v_{2,1}=(\frac{1}{3},\frac{1}{3},-\frac{2}{3})$, $\v_{2,2}=(0,0,1)$. We see $\v_{1,1}+\v_{1,2}=\v_{2,1}+\v_{2,2}=\frac{1}{3}\one$. Then $\kappa(\P)=|\det(\v_{1,1},\v_{2,1},\one)|=1$. The support functions are then calculated as
\begin{align*}
h_{\x\P}(\v_{1,1})=&\frac{1}{3}x_1+\frac{1}{3}x_2-\frac{2}{3}x_3,\\
h_{\x\P}(\v_{1,2})=&\frac{1}{3}x_1-\frac{2}{3}x_2+\frac{1}{3}x_3,\\
h_{\x\P}(\v_{2,1})=&-\frac{1}{3}x_1+\frac{2}{3}x_2-\frac{1}{3}x_3,\\
h_{\x\P}(\v_{2,2})=&0,
\end{align*}
and indeed $h_{\x\P}(\v_{i,1})+h_{\x\P}(\v_{i,2})=x_i-\frac{1}{3}\langle\x,\one\rangle$ for $i\in\{1,2\}$. Proposition~\ref{prop:{lemma}affine-fine-mixed-cell-dual-volume} gives
\[\tilde m_{\P}(\x,\z)=\frac{(x_1-\frac{1}{3}\langle\x,\one\rangle+\frac{1}{3}\langle\z,\one\rangle)(x_2-\frac{1}{3}\langle\x,\one\rangle+\frac{1}{3}\langle\z,\one\rangle)}{(\frac{1}{3}x_1{+}\frac{1}{3}x_2{-}\frac{2}{3}x_3{+}\langle\z,\v_{1,1}\rangle)(\frac{1}{3}x_1{-}\frac{2}{3}x_2{+}\frac{1}{3}x_3{+}\langle\z,\v_{1,2}\rangle)({-}\frac{1}{3}x_1{+}\frac{2}{3}x_2{-}\frac{1}{3}x_3{+}\langle\z,\v_{2,1}\rangle)(z_3)}.\]
\end{example}

\subsection{Fine mixed subdivisions in an affine hyperplane}
The following Proposition is immediate from Proposition~\ref{prop:hyperplane-mixed} and Theorem~\ref{thm: mixed subdivision formula}.
\begin{prop}\label{prop:affine-fine-mixed-subdivision}
Let $\P=(P_1,\ldots,P_r)$ be a sequence of polytopes in $H_1$ such that $P=P_1+\cdots+P_r$ is full-dimensional in $H_r$. Let $\mathcal{S}=\{\Q^{(1)},\ldots,\Q^{(N)}\}$ be a fine mixed subdivision of $\P$. Then \[\tilde m_{\P}(\x,\z)|_{\langle\x,\one\rangle=\langle\z,\one\rangle}=\sum_{\Q\in\mathcal{S}}\tilde m_{\Q}(\x,\z)|_{\langle\x,\one\rangle=\langle\z,\one\rangle}.\]
\end{prop}

\section{Zonotopes} \label{sec:zonotope}
\begin{defin}
A \emph{zonotope} $Z$ with $r$ \emph{zones} is a Minkowski sum of $r$ intervals.
\end{defin}
In this section, consider a sequence of intervals $\P=(P_1,\ldots,P_r)$ where $P_i=[-\p_i,\p_i]$ for $i=1,\ldots,r$. Assume $\P$ is full-dimensional, i.e. $\p_1,\ldots,\p_r$ span $\R^d$. Also assume for simplicity that $\p_i\neq0$ for all $i$ and $\p_i\neq\pm\p_j$ for $i\neq j$.

Fine mixed subdivisions of $P_1+\cdots+P_r$ are called \emph{zonotopal tilings}.

\subsection{Deletion-contraction formalism}
We will use deletion-contraction to compute the dual mixed volume function of zonotopes. We setup the inductive steps in this section. 

Let $P \subseteq \R^d$ be a full-dimensional polytope and $\p \in \R^d$ be a vector.  Let $H$ be the hyperplane normal to $\p$ and let $H_{>0}$ (resp. $H_{\geq 0}$) and $H_{<0}$ (resp. $H_{\leq 0}$) denote the corresponding closed (resp. open) halfspaces. Let 
$$
P(x) := P+x[-\p,\p]
$$
The aim is to give a recursive description of the rational function
$$
V(x,\z):=\VolC_\z(P(x)) = \VolC_\z(P+x[-\p,\p]).
$$
Let $V(\z) = V(0,\z) = \VolC_\z(P)$.
Let $P_{/\p}:={\rm proj}_\p P$ denote the orthogonal projection of $P$ into $H$, in the direction of $\p$. The goal is to describe $V(x,\z)$ in terms of $V(\z)$, the ``deletion" part, and $\VolC_{\z_0}(P_{/\p})$, the ``contraction" part.  Here, $\z_0$ denotes a point in the hyperplane $H \subseteq \R^d$, and the volume in $H$ is normalized as in Section~\ref{sub:normalization-affine}.

Let the denominator $B_{\z}$ of $V(\z)$ as defined in \eqref{eq:denominator of VolC} equal $B_+(\z)B_0(\z)B_-(\z)$ where the three factors correspond to the rays of the normal fan of $P$ that belong to $H_{>0}$, the hyperplane $H$, and $H_{<0}$ respectively.  

\begin{lemma}\label{lem:relprime}
    The three polynomials $B_+,B_0,B_-$ are pairwise relatively prime.
\end{lemma}
\begin{proof}
    Each of these polynomials is a product of affine-linear forms in $\z = (z_1,z_2,\ldots,z_d)$. 
 The linear parts of two such factors are proportional if and only if the corresponding rays of $\N(P)$ are parallel.  This happens only if we have two rays of $\N(P)$ that are opposite.  The fact that the affine-linear forms are not proportional then follows from the assumption that $P$ is full-dimensional.
\end{proof}

The cones of the normal fan $\N'$ of $P+[-\p,\p]$ are obtained from $\N(P)$ by intersecting with the cones $H_{\geq 0}$, $H$, and $H_{\leq 0}$. Let $D(\z)$ be the product of linear factors $h_{P-\z}(\v_i)$'s where $\v_i$'s are rays of $\N'$ that are not rays of $\N$.

\begin{theorem}\label{thm:project}
There exists a unique pair of rational functions $W_+=W_+(\z)$ and $W_- = W_-(\z)$ satisfying the following properties:

\begin{enumerate}
\item  $W_+(\z) + W_-(\z) = V(\z)$;
\item $W_{\pm}(\z) = A_{\pm}(\z)/\big(B_{\pm}(\z)B_0(\z)D(\z)\big)$ for some polynomial $A_{\pm}(\z)$;
\item writing $\z = \z_0 + t \p$ with $\z_0 \in H$, then $\lim_{t \to \infty} t W_{\pm}(\z_0,t) = \pm\Vol^\vee_{\z_0}(P_{/\p})/||\p||^2$.
\end{enumerate}

Consequently, we have $V(x,\z) = W_+(\z+x\p)+W_-(\z-x\p)$.
\end{theorem}
\begin{proof}
Let $\N_{+}$ (resp. $\N_{-}$) be the fan consisting of those cones of $\N'$ that belong to $H_{\geq 0}$ (resp. $H_{\leq 0}$). With notations in Definition~\ref{def:universal-fan-dual-volume}, we define 
$$
W_{\pm}:= f_{\N_{\pm}}(h_{P-\z}),
$$
It is clear that 
$$
W_+ + W_- = f_{\N_{+}}(h_{P-\z}) + f_{\N_{-}}(h_{P-\z}) = f_{\N'}(h_{P-\z}) = f_{\N}(h_{P-\z}) = V(\z),
$$
since $\N'$ is a refinement of $\N$.  It is also clear that $W_{\pm}$ satisfies property (2).  

Let us prove property (3).  By definition, the support function $h_{P_{/\p}}$ is the restriction of $h_P$ to the hyperplane $H$.  It follows that the normal fan of $P_{/\p}$, denoted as $\N_{/\p}$, consists of the cones $\N(P) \cap H$.  Equivalently, we have $\N_{/\p} = \N' \cap H$ as a complete fan in $H$.  The rational function $W_+$ is a sum over maximal cones of $\N_{+}$.  After triangulating, we may assume that the maximal cones are simplicial.  For a simplicial maximal cone $C$ spanned by $\v_1,\ldots,\v_d$ with a facet along $H$, we have  (c.f. Definition~\ref{def:universal-fan-dual-volume})
$$
f_C(h_{P-\z}) = \frac{|\det(\v_1,\ldots,\v_d)|}{\prod_{i=1}^d h_{P-\z}(\v_i)}
$$
where $h_{P-\z}(\v_i)=-\min_{\y\in P}\langle\y,\v_i\rangle+\langle\p,\v_i\rangle t+\langle\z_0,\v_i\rangle$. Since $C$ has a facet along $H$, we assume without loss of generality that $\v_2,\ldots,\v_d$ belong to $H$, meaning that $\langle\p,\v_i\rangle=0$ for $i=2,\ldots,d$. Thus, $h_{P-\z}(\v_1)$ depends on $t$, while $h_{P-\z}(\v_2),\ldots,h_{P-\z}(\v_d)$ do not depend on $t$. We now have 
\[\lim_{t\rightarrow\infty}tf_C(h_{P-\z})=\frac{|\det(\v_1,\ldots,\v_d)|}{\prod_{i=2}^d h_{P-\z}(\v_i)}\lim_{t\rightarrow\infty}\frac{t}{h_{P-\z}(\v_1)}=\frac{|\det(\v_1,\ldots,\v_d)|}{\prod_{i=2}^d h_{P-\z}(\v_i)}\frac{1}{\langle\p,\v_1\rangle}
\]
If we write $\v_1=k\p+\v_1'$ with $\v_1'\in H$, then $\langle\p,\v_1\rangle=k||\p||^2$ and $\det(\v_1,\ldots,\v_d)=k\det(\p,\v_2,\ldots,\v_d)$. Recall from Section~\ref{sub:normalization-affine} that $|\det(\p,\v_2,\ldots,\v_d)|$ is how we normalize the volume of the simplex spanned by $\v_2,\ldots,\v_d$ when we treat the affine subspace $H$ as the ambient vector space. Moreover, for $2\leq i\leq d$, as $\v_i$ belongs to $H$, $h_{P-\z}(\v_i)=h_{P/\p-\z_0}(\v_i)$. Continuing the above calculation, we have
\[\lim_{t\rightarrow\infty}tf_C(h_{P-\z})=\frac{|\det(\p,\v_2,\ldots,\v_d)|}{\prod_{i=2}^d h_{P/\p-\z_0}(\v_i)}\frac{1}{||\p||^2}=f_{C\cap H}(h_{P/\p-\z_0})\frac{1}{||\p||^2}.\]

On the other hand, if $C$ is a simplicial maximal cone where $H$ is not a facet, then at least two of the $h_{P-\z}(\v_1),\ldots,h_{P-\z}(\v_d)$ depend on $t$, resulting in $\lim_{t \to \infty} t f_C(h_{P-\z}) = 0$.  Summing over maximal cones, we obtain (3).

The stated formula $V(x,\z) = W_+(\z+x\p)+W_-(\z-x\p)$ follows from 
$$
h_{P+x[-\p,\p]}(\y) = \begin{cases} h_P(\y) + x \langle \p, \y \rangle & \mbox{if $\y \in H_{\geq 0}$} \\
h_P(\y) - x \langle \p, \y \rangle & \mbox{if $\y \in H_{\leq 0}$}.
\end{cases}
$$

It remains to show that the pair $(W_+,W_-)$ is unique.  Write 
$$
\frac{A_+}{B_+B_0D} + \frac{A_-}{B_-B_0D} = \frac{A}{B_+B_-B_0} = V(\z)
$$
and clear denominators to obtain the polynomial equality
\begin{equation}\label{eq:AA}
A_+B_- +A_-B_+= AD. 
\end{equation}
The choice of solution $(W_+,W_-)$
 is equivalent to the choice of the pair $(A_+,A_-)$ of polynomials.
Suppose that $(A'_+,A'_-)$ is another solution to \eqref{eq:AA}.  Then $(A_+-A'_+)B_- = (A'_--A_-)B_+$ and by \cref{lem:relprime} we conclude that $(A_+-A'_+)$ is divisible by $B_+$.  It follows that $W'_+ = \frac{A'_+}{B_+B_0D} = W_+ + \frac{E}{B_0D}$ for some polynomial $E = E(\z)$.  Writing all polynomials in the variables $(\z_0,t)$ we have that $B_0,D$ do not depend on $t$ but $E(\z)$ may or may not.  If $E$ is nonzero, the limit $\lim_{t \to \infty} t \frac{E}{B_0D}$ is therefore unbounded, so $W'_+$ cannot satisfy (3).  It follows that $(W_+,W_-)$ is unique.
 \end{proof}

\begin{example}
Consider the polytope $P=\conv((-2,-1),(0,1),(1,-1))$ in $\R^2$ with $\p=(1,0)$. Note that $||\p||^2=1$. The rays of its normal fan $\N$ are $\v_1=(-2,-1)$, $\v_2=(0,1)$, $\v_3=(1,-1)$ and there is an additional ray $\v_4=(0,-1)$ in the normal fan $\N'$ of $P=[-\p,\p]$, shown in Figure~\ref{fig:deletion-contraction-ex}.
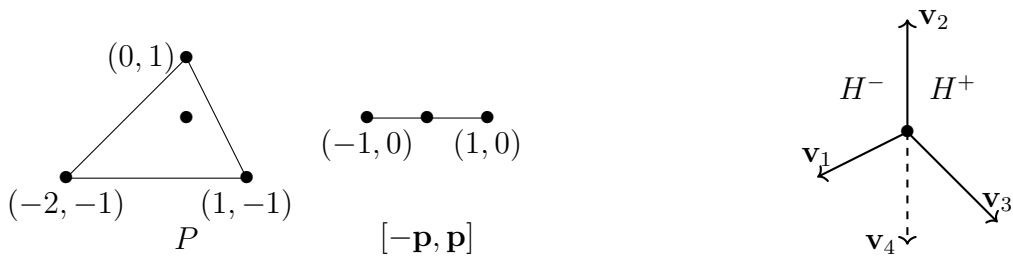
\begin{figure}[h!]
\centering
\begin{tikzpicture}[scale=0.8]
\node at (-2,-1) {$\bullet$};
\node at (0,1) {$\bullet$};
\node at (1,-1) {$\bullet$};
\node at (0,0) {$\bullet$};
\draw(-2,-1)--(0,1)--(1,-1)--(-2,-1);
\node[below] at (-2,-1) {$({-}2,{-}1)$};
\node[below] at (1,-1) {$(1,{-}1)$};
\node[left] at (0,1) {$(0,1)$};
\node at (0,-2) {$P$};

\node at (3,0) {$\bullet$};
\node at (4,0) {$\bullet$};
\node at (5,0) {$\bullet$};
\node[below] at (3,0) {$({-}1,0)$};
\node[below] at (5,0) {$(1,0)$};
\draw(3,0)--(5,0);
\node at (4,-2) {$[-\p,\p]$};
\end{tikzpicture}
\qquad\qquad\qquad\qquad
\begin{tikzpicture}[scale=0.6]
\node at (0,0) {$\bullet$};
\draw[->,thick](0,0)--(0,2.5);
\draw[->,thick](0,0)--(-2,-1);
\draw[->,thick](0,0)--(2,-2);
\draw[->,dashed,thick](0,0)--(0,-2.5);
\node[right] at (0,2.5) {$\v_2$};
\node[above] at (2,-2) {$\v_3$};
\node[left] at (0,-2.5) {$\v_4$};
\node[above] at (-2,-1) {$\v_1$};
\node at (1,1) {$H^+$};
\node at (-1,1) {$H^-$};
\end{tikzpicture}
\caption{Left: $P$ and $[-\p,\p]$. Right: $\N(P)$ (undashed) and $\N(P+[-\p,\p])$.}
\label{fig:deletion-contraction-ex}
\end{figure}

We calculate that $h_{P-\z}(\v_1)=1-2z_1-z_2=B^-(\z)$, $h_{P-\z}(\v_2)=1+z_2=B_0(\z)$, $h_{P-\z}(\v_3)=1+z_1-z_2=B^+(\z)$ and $h_{P-\z}(\v_4)=1-z_2=D(\z)$. We then have, by construction in the proof of Theorem~\ref{thm:project},
\begin{align*}
W_+=&f_{\N_+}(h_{P-\z})=\frac{1}{(1+z_1-z_2)(1+z_2)}+\frac{1}{(1+z_1-z_2)(1-z_2)},\\
W_-=&f_{\N_-}(h_{P-\z})=\frac{2}{(1-2z_1-z_2)(1+z_2)}+\frac{2}{(1-2z_1-z_2)(1-z_2)}.
\end{align*}
As a sanity check, the last two terms in $W_+$ and $W_-$ sum up to $f_C(h_{P-\z})$ where $C$ is the cone spanned by $\v_1$ and $\v_3$, so parts (1) and (2) and Theorem~\ref{thm:project} are satisfied. Write $\z=\z'+t\p$ so that $z_1=z_1'+t$ and $z_2=z_2'$. Then
\[\lim_{t\rightarrow\infty}tW_+(\z',t)=\lim_{t\rightarrow\infty}\frac{t}{(1+z'_1+t-z_2)(1+z'_2)}+\frac{t}{(1+z'_1+t-z_2)(1-z_2)}=\frac{1}{1+z_2'}+\frac{1}{1-z_2'}\]
On the other hand, the orthogonal projection $P/\p$ is a line segment from $-1$ to $1$, with coordinate $z_2'$. Its dual volume $\VolC_{\z'}(P/\p)$ is precisely the above expression (Example~\ref{ex:one-dimensional-dual-volume}). We can similarly calculate that $\lim_{t\rightarrow\infty}tW_{-}(\z',t)=-\VolC_{\z'}(P/\p)$.
\end{example}


Let $\P = ([-\p_1,\p_1],\ldots,[-\p_r,\p_r])$ be a zonotope.  We let $\P' = ([-\p_1,\p_1],\ldots,[-\p_{r-1},\p_{r-1}])$ be the deletion and $\P'' = \P/\p_r = ([-\p_1,\p_1],\ldots,[-\p_{r-1},\p_{r-1}]) \subseteq \R^d/\p_r$ be the contraction. 

\begin{cor}
    The dual mixed volume of the zonotope $\P$ depends only on the dual mixed volume functions of $\P'$ and $\P''$. 
\end{cor}

Recursively applying \cref{thm:project}, one can compute the dual mixed volume of any zonotope via contraction-deletion.

\subsection{Oriented matroids and the Bohne-Dress Theorem}
We need the language of oriented matroids. Readers are referred to \cite{oriented-matroid-book} for details on this subject matter.

Let $E=[r]$ be the \emph{ground set}. For a \emph{sign vector} $\epsilon\in\{0,+,-\}^E$, define its \emph{zero part} as $\epsilon^0:=\{i\in E\:|\:\epsilon_i=0\}$. Similarly, we define its \emph{positive part} $\epsilon^+$ and \emph{negative part} $\epsilon^-$. There is a partial order on signs given by $0<+$ and $0<-$, while the nonzero signs are incomparable. Let $\zeta, \xi$ be two sign vectors, their \emph{product} is defined as $(\zeta\circ\xi)_i=\xi_i$ if $\zeta_i<\xi_i$ and $(\zeta\circ\xi)_i=\zeta_i$ otherwise. The \emph{separation set} $S(\zeta,\xi):=\{i\in E\:|\: \zeta_i=-\xi_i\neq0\}$.
\begin{defin}
A set $\M\subset\{0,+,-\}^E$ is the set of \emph{covectors} of an \emph{oriented matroid} if
\begin{itemize}
\item $\mathbf{0}\in\M$;
\item $\epsilon\in\M$ implies $-\epsilon\in\M$;
\item $\zeta,\xi\in \M$ implies $\zeta\circ\xi\in \M$;
\item for $\zeta,\xi\in \M$ and $i\in S(\zeta,\xi)$, there exists $\epsilon\in\M$ such that $\epsilon_i=0$ and $\epsilon_j=(\zeta\circ\xi)_j=(\xi\circ\zeta)_j$ for all $j\notin S(\zeta,\xi)$.
\end{itemize}
\end{defin}
The \emph{deletion} $\M\backslash j$ and the \emph{contraction} $\M/ j$ are defined as
\begin{align*}
\M\backslash j:=&\{(\epsilon_1,\ldots,\epsilon_{j-1},\epsilon_{j+1},\ldots,\epsilon_r)\:|\: \epsilon\in\M\},\\
\M/ j:=&\{(\epsilon_1,\ldots,\epsilon_{j-1},\epsilon_{j+1},\ldots,\epsilon_r)\:|\: \epsilon\in\M,\epsilon_j=0\}.
\end{align*}
We say that $j\in E$ is a \emph{loop} of $\M$ if $\epsilon_j=0$ for all $\epsilon\in M$.
\begin{defin}\label{def:one-element-lifting}
Let $\M$ be an oriented matroid on $E$. A \emph{one-element lifting} of $\M$ is an oriented matroid $\tilde\M$ on $E\sqcup\{g\}$ such that $\widetilde\M/ g=\M$ and $g$ is not a loop of $\widetilde\M$.
\end{defin}

Given the data $\P=(P_1,\ldots,P_r)$ of a zonotope, we can naturally associate an oriented matroid as follows \[\M(\P):=\left\{\big(\text{sign}(\langle\p_1,\v\rangle),\ldots,\text{sign}(\langle\p_r,\v\rangle)\big)\:|\: \v\in\R^d\right\}\subset\{0,+,-\}^r.\]
Also, for a sign vector $\epsilon\in\{0,+,-\}^r$, write \[P_i(\epsilon):=\begin{cases}
P_i &\text{ if } \epsilon_i=0,\\
\p_i&\text{ if }\epsilon_i=+,\\
-\p_i&\text{ if }\epsilon_i=-.
\end{cases}\] Define $\P(\epsilon)=\sum_{i=1}^r P_i(\epsilon)$.
We say that $\P(\epsilon)$ has support $\epsilon^0$. Note that $\P(\epsilon)$ is a parallelotope exactly when $\{\p_i\:|\: i\in\epsilon^0\}$ is a basis of $\R^d$. For the following result, see \cite{bohne-dress-richter-gebert-ziegler}.

\begin{theorem}[Bohne-Dress Theorem] \label{thm:BD}
Zonotopal tilings of $\P=(P_1,\ldots,P_r)$ are in bijection with one-element liftings $\widetilde\M$ of the oriented matroid $\M(\P)$. Specifically, given such a lifting $\widetilde\M$ on $[r]\sqcup\{g\}$, \[\left\{\Q^{\epsilon}=(P_1(\epsilon),\ldots,P_r(\epsilon))\:\mid\: \epsilon\in\widetilde\M,\epsilon_g=+,|\epsilon^0|=d\right\}\]
is a zonotopal tiling, i.e. mixed subdivision, of $\P$.
\end{theorem}

When $\widetilde\M$ is a generic one-element lifting, the corresponding zonotopal tiling is fine.  Combining with \cref{thm: mixed subdivision formula} and \cref{prop:fine-mixed-dual-volume} we have the following result.
\begin{cor}\label{cor:BD}
For each one-element lifting $\widetilde\M$ of $\M(\P)$, we have 
$$
m_\P(\x) = \sum_{\epsilon} m_{\Q^\epsilon}(\x)
$$
where the summation is over sign vectors $\epsilon \in \widetilde\M$ satisfying $\epsilon_g = +$ and $|\epsilon^0| = d$.  When $\widetilde\M$ is generic, each $m_{\Q^\epsilon}(\x)$ is a monomial.
\end{cor}
It would be interesting to compare the formulae from \cref{cor:BD} with the one obtained from \cref{thm:project}.  We take one initial step in this direction.

Let $m_{\P'}(\x) = \sum_{\epsilon} m_{\Q^\epsilon}(\x)$ arise from \cref{cor:BD} from a fine zonotopal tiling, or equivalently, from a generic one-element lifting.  Then each $\Q^{\epsilon}$ is a parallelotope and 
$$
m_{\Q^\epsilon}(\x,\z) = \kappa(\epsilon)\prod_{i \in \epsilon^0}\frac{x_i}{(x_i^2/4-\langle\y,\v_i\rangle^2)}
$$ 
has been described in \cref{ex:par}, where we take $\y = \z-\sum_{j=1}^{r-1} \epsilon_j \p_j$, and $\{\v_1,\v_2,\ldots,\v_d\}$ is the dual basis to $\{\p_i \mid i \in \epsilon^0\}$.  We may apply \cref{thm:project} to each $\Q^\epsilon$ separately, and obtain a decomposition
$$
m_{\Q^\epsilon}(\x,\z) = m_{\Q^\epsilon,+}(\x,\z) + m_{\Q^\epsilon,-}(\x,\z).
$$
The following result states that these decompositions are compatible with the decomposition of $m_{\P'}(\x,\z)$.

\begin{lemma}
The decomposition of $m_{\P'}(\x,\z)$ from \cref{thm:project} is given by
$$
m_{\P'}(\x,\z) = \left(\sum_{\epsilon}m_{\Q^\epsilon,+}(\x,\z)\right) + \left(\sum_{\epsilon}m_{\Q^\epsilon,-}(\x,\z)\right).
$$
\end{lemma}
\begin{proof}
We showed in the proof of \cref{thm:project} that the decomposition $V(\z) = W_+(\z) + W_-(\z)$ is uniquely determined by the form of the denominators and the fact that the limit $\lim_{t \to \infty} W_{\pm}(\z_0, t)$ exists.  Both these properties can be directly verified for $\left(\sum_{\epsilon}m_{\Q^\epsilon,\pm}(\x,\z)\right)$.
\end{proof}

\section{Generalized permutohedra}\label{sec:genperm}
Readers are referred to \cite{postnikov-permutohedra} for details on generalized permutohedra, which are a larger family of polytopes than zonotopes that we have considered, and appear everywhere in algebra and combinatorics.
In this section, we switch our usual notation for dimension from $d$ to $n$ to be more consistent with the literature.

\subsection{Generalized permutohedra in an affine hyperplane}
In the literature, generalized permutohedra are defined to live in an affine hyperplane. 

Recall that $H_t = \{\y = (y_1,\ldots,y_n) \mid y_1+y_2 + \cdots + y_n = t\} \subseteq \R^n$.
\begin{defin}\label{def:generalized-permuohedra}
For $T\neq\emptyset\subset[n]$, define the following $(|T|-1)$-dimensional simplex \[\Delta_T := \conv(e_i \mid i \in T)\subseteq H_1.\]
Let $\P=(\Delta_T)_{T\in 2^{[n]}\setminus \{\emptyset \}}$ be the collection of all such simplices and $\x = (x_T)_{T\in 2^{[n]}\setminus \{\emptyset \}}$. Then the \emph{generalized permutohedron} is \[\x\P:=\sum_{T \in 2^{[n]} \setminus \{\emptyset\}} x_T \Delta_T \subseteq H_{\langle\x,\one\rangle}.\]
\end{defin}
In this section, we use different ways to calculate the hyperplane dual mixed volume of the generalized permutohedron to obtain identities. Our main focus is on $\tilde m_{\P}(\x,\z)$ in the case when $\langle\x,\one\rangle=0$ and $\z=\mathbf{0}$.

\begin{proposition}\label{prop:genperm}
When $\langle\x,\one\rangle=0$, or equivalently, substituting $x_{[n]}$ with $-\sum_{T\in 2^{[n]}\setminus \{[n]\}}x_T$, the hyperplane dual mixed volume of the generalized permutohedron becomes \[\tilde m_{\P}(\x)|_{\langle\x,\one\rangle=0}=(-1)^{n-1}\sum_{\sigma\in S_n}\prod_{a=1}^{n-1}\frac{1}{\displaystyle\sum_{\emptyset\neq T\subseteq\sigma[1:a]}x_T}.\]
\end{proposition}
\begin{proof}
Write $P=\sum_{T\in 2^{[n]}\setminus \{\emptyset\} }\Delta_T$ and we calculate $\Cone(P)^*$. The reduced normal fan of $P$ in $\R^n/\one$ is called the \emph{reduced braid fan} in the literature and it's structure is well-known. Its maximal cones are in bijection with permutations in $S_n$. Specifically, for a permutation $\sigma\in S_n$, its corresponding maximal cone $F_\sigma$ is generated by $\one_{\sigma[1:1]},\one_{\sigma[1:2]},\ldots,\one_{\sigma[1:n-1]}$ where $\sigma[i:j]=\{\sigma(i),\sigma(i+1),\ldots,\sigma(j)\}$, and $\one_{A}$ (for any subset $A\subseteq [n]$) is the indicator vector of $A$, i.e. $\langle\one_{A},\y\rangle=\sum_{a\in A}y_a$. 

In order to adapt to our language in Section~\ref{sec:affine}, we need to lift each $\one_{A}$, where $A\subset[n]$ is a proper subset, to an appropriate vector $\v_{A}=\one_{A}-b_A\one$ which is a generator of the boundary of $\Cone(P)^*$, for some $b_{A}\in\R$. We thus need
\begin{align*}
0=&\min_{\y\in P}\langle\v_{A},\y\rangle\\
=&\sum_{T\in 2^{A}\setminus \{\emptyset \} }(1-b_A)+\sum_{T\in 2^{[n]}\setminus \{ \emptyset \},\ T\nsubseteq A}(-b_A)\\
=&(2^{|A|}-1)(1-b_A)+(2^n-2^{|A|})(-b_A)
\end{align*}
which implies that $b_A=(2^{|A|}-1)/(2^n-1)$. Denote this number by $b_k=b_{A}$ where $|A|=k$, for $k=1,\ldots,n-1$. Note that the precise value of $b_k$ does not matter much. 

For a permutation $\sigma\in S_n$, we have \[|\det(\v_{\sigma[1:1]},\v_{\sigma[1:2]},\ldots,\v_{\sigma[1:n-1]},\one)|=|\det(\one_{\sigma[1:1]},\ldots,\one_{\sigma[1:n-1]},\one)|=1.\]

For a proper subset $A \subsetneq [n]$ of cardinality $k$, we can calculate that 
\begin{align*}
h_{\x\P}(\v_A)=&\sum_{T\in 2^{[n]}\setminus \{\emptyset \}}x_Th_{\Delta_T}(\v_A)\\
=&\sum_{T\in 2^{A}\setminus \{\emptyset \}}-(1-b_k)x_T+\sum_{T\in 2^{[n]}\setminus \{ \emptyset\},\ T\nsubseteq A}b_kx_T\\
=&b_k\langle\x,\one\rangle-\sum_{T\in 2^A \setminus \{ \emptyset\} }x_T,
\end{align*}
which equals $-\sum_{T\in 2^A \setminus \{ \emptyset\}}x_T$ when $\langle\x,\one\rangle=0$. By Definition~\ref{def:hyperplane-dual-mixed},
\begin{align*}
\tilde m_{\P}(\x)|_{\langle\x,\one\rangle=0}=&\sum_{\sigma\in S_n}\frac{|\det(\v_{\sigma[1:1]}, \ldots,\v_{\sigma[1:n-1]},\one)|}{h_{\x\P}(\v_{\sigma[1:1]})|_{\langle\x,\one\rangle=0} \cdots h_{\x\P}(\v_{\sigma[1:n-1]})|_{\langle\x,\one\rangle=0}}\\
=&\sum_{\sigma\in S_n}\frac{1}{\prod_{a=1}^{n-1} \big(-\sum_{\emptyset\neq T\subseteq\sigma[1:a]}x_T \big)}\\
=&(-1)^{n-1}\sum_{\sigma\in S_n}\prod_{a=1}^{n-1}\frac{1}{\sum_{\emptyset\neq T\subseteq\sigma[1:a]}x_T}.
\end{align*}
\end{proof}

\begin{example}
Let $n = 3$.  Then we have
\begin{align*}
\tilde m_{\P}(\x)|_{\langle\x,\one\rangle=0} &= \frac{1}{x_1(x_1+x_2+x_{12})} +  \frac{1}{x_1(x_1+x_3+x_{13})} +  \frac{1}{x_2(x_1+x_2+x_{12})}  \\
&+  \frac{1}{x_2(x_2+x_3+x_{23})} +  \frac{1}{x_3(x_1+x_3+x_{13})} +  \frac{1}{x_3(x_2+x_3+x_{23})}.
\end{align*}
\end{example}

\subsection{Subdivisions of generalized permutohedra}
One can apply the Cayley trick as in Section~\ref{sec:Cayley} to the generalized permutohedra and obtain fine mixed subdivisions of $\P$. We state results of Postnikov \cite{postnikov-permutohedra} and use them for our computation of the dual mixed volume of the generalized permutohedra.

Let $N=2^n-1$ and let $\{T_1,\ldots,T_N\}$ be all nonempty subsets of $[n]$. Let $K_{N,n}$ be the complete bipartite graph with vertices in one part are labeled $i=1,2,\ldots,N$, and vertices the other part are $\bar j=\bar 1,\bar2,\ldots,\bar n$. Let $G\subseteq K_{N,n}$ be the subgraph where the vertex $i$ is connected to $\bar j$ if $j\in T_i$.
\begin{defin}
For a collection of subsets $\J=(J_1,\ldots,J_N)$ where $J_i\subseteq T_i$, let $G_{\J}$ be the subgraph of $G$ where vertex $i$ is connected to $\bar j$ if $j\in J_i$.
\end{defin}
\begin{lemma}\cite[Lemma 14.7]{postnikov-permutohedra}
For any fine mixed subdivisions of $\P$, its fine mixed cells are in bijection with nonempty subsets $\J=(J_1,\ldots,J_N)$ such that $G_{\J}$ is a spanning tree of $G$. In this case, the corresponding fine mixed cell of $\P$ is $\Q_{\J}:=(\Delta_{J_i})_{i=1}^N$.
\end{lemma}
When $G_{\J}$ is a spanning tree, removing the vertex $i\in[N]$ and the edges incident to it in the graph $G_{\J}$ results in $|J_i|$ connected components, each containing exactly one element in $J_i$. For an edge $(i,\bar j)$ in $G_{\J}$, denote this connected component that contains $\bar j$ as $A_{i,\bar j}\sqcup \bar A_{i,\bar j}$ where $A_{i,\bar j}\subseteq[N]$ and $\bar A_{i,\bar j}\subseteq[\bar n]=\{\bar1,\bar2\ldots,\bar n\}$. Define \[h^{\J}_{i,\bar j}(\x):=\sum_{k\in A_{i,\bar j}}x_{T_k}.\]
In other words, $x_{T_k}$ is in the summation for $h^{\J}_{i,j}(\x)$ if and only if the path between $k$ and $i$ in $G_{\J}$ passes through $\bar j$.
Note that 
\[\sum_{j\in J_i}h^{\J}_{i,\bar j}(\x)=\sum_{T_i\neq T\in 2^{[n]}\setminus\emptyset}x_T.\]

\begin{proposition}\label{prop:permutohedra-cell-dual-volume}
Let $\J$ be a collection of subsets so that $G_{\J}$ is a spanning tree of $G$. When $\langle\x,\one\rangle=0$, the hyperplane dual mixed volume of the fine mixed cell $\Q_{\J}$ is \[\tilde m_{\Q_{\J}}(\x)|_{\langle\x,\one\rangle=0}=(-1)^{n-1}\prod_{i:|J_i|>1}\frac{-x_{T_i}}{\prod_{j\in J_i}h_{i,\bar j}^{\J}(\x)}.\]
\end{proposition}
\begin{proof}
Let $d_i=|J_i|-1$ for $i=1,\ldots,N$. In the current scenario, the majority of the $d_i$'s will be $0$, meaning that $J_i$ is a singleton. For $d_i>0$, let $J_i=\{c_1,\ldots,c_{d_i+1}\}\subseteq T_i$. We need to calculate the generator rays of $\Cone(Q_{\J})^*$, where $Q_{\J}=\sum_{i=1}^N \Delta_{J_i}$. Let $a=1,2,\ldots,d_i+1$, the linear equations to solve are given by
\[\left\langle\v_{i,a},\e_{c_a}+\sum_{i'\neq i}\e_{c_{i'}}\right\rangle=1,\qquad\left\langle\v_{i,a},\e_{c_{a'}}+\sum_{i'\neq i}\e_{c_{i'}}\right\rangle=0,\]
for any $c_{i'}\in J_{i'}$ and any $a'\neq a$ in $J_i$. Subtracting these two equations, we assume that $\v_{i,a}$ has value $-b$ at coordinate $c_{a'}$ for $a'\neq j$ and has value $1-b$ at coordinate $c_{a}$. Note that if some $J_{i'}$ has cardinality greater than $1$, then $\v_{i,a}$ must take on the same value at all coordinates in $J_{i'}$. For $j\in J_i$, recall that $A_{i,\bar j}\sqcup \bar A_{i,\bar j}$ is the connected component of $G_{\J}$ taken away the vertex $i$ that contains $\bar j$. The connectivity of $\bar A_{i,\bar j}$ means that $\v_{i,a}$ takes the same value on coordinates in $\bar A_{i,\bar j}$. As a summary, $\v_{i,a}$ equals $1-b$ at coordinates in $\bar A_{i,\bar c_{a}}$ and equals $-b$ at all other coordinates. Write $\v_{i,\bar j}$ for $\v_{i,a}$ where $j=c_a$.

The calculation that $\kappa(\Q_{\J})=1$ is a straightforward graph-theoretic arguments, which we establish in the following claim, with self-contained notations.
\begin{claim}
Let $G$ be a bipartite tree with parts $U$ and $V$. For any vertex $a$, we denote the set of its neighbors by $D(a)$. For an edge $(u,v)$ of $G$, define $h_{u,v}\in\R^{V}$ as a vector with $\{0,1\}$ entries. The entry at coordinate $v'\in V$ equals $1$ if the unique path between $v$ and $v'$ in $G$ passes through $u$, and $0$ otherwise. For each $u\in U$, we select $(|D(u)|-1)$ of its neighbors, meaning we choose all its neighbors $v$ except for one. The collection of all such vectors $h_{u,v}$, along with $\one$, form a matrix with determinant $\pm1$.
\end{claim}
\begin{proof}[Proof of Claim]
We start with a sanity check on the dimension of this matrix. $G$ is a tree, so it has $|U|+|V|-1$ edges. Except $\one$, we have collected a total of \[\sum_{u\in U}(|D(u)|-1)=\left(\sum_{u\in U}|D(u)|\right)-|U|=|E(G)|-|U|=|V|-1\]
vectors in $\R^{V}$, so the matrix is in $\R^{V\times V}$ and we can thus take the determinant. Also notice that, for a fixed $u\in U$, $\sum_{v\in D(u)}h_{u,v}=\one$ so the choices of its neighbors will not matter. We now do induction on $|U|+|V|$.

The base case where $|U|=|V|=1$ is evident. Since $G$ is a tree, it must contain a leaf. If $u\in U$ is a leaf, this vector simply does not contribute to the calculation of the determinant, and we can remove it and reduce it to the case of $G\setminus u$. If $v\in V$ is a leaf, let $u\in U$ be connected to it. Then $h_{u,v}$ equals $1$ at coordinate $v$ and equals $0$ elsewhere. In the determinant calculation, we must then choose this entry, and reduce it to the case of the tree $G\setminus v$. In either case, we are done by the inductive hypothesis.
\end{proof}

It remains to compute some support functions. Now $h_{\Delta_{J_k}}(\v_{i,\bar j})=-\min_{t\in J_k}\langle\v_{i,\bar j},\e_t\rangle$, which equals $-1+b$ if $J_k$ is a subset of $\bar A_{i,\bar j}$, and $b$ otherwise. The condition that $J_k$ is contained inside $\bar A_{i,\bar j}$ is equivalent to $k\in A_{i,\bar j}$. As a result, \[h_{\x\Q_{\J}}(\v_{i,\bar j})=\sum_{k=1}^{N}x_{T_k}h_{\Delta_{J_k}}(\v_{i,\bar j})=b\langle\x,\one\rangle-\sum_{k\in A_{i,\bar j}}x_{T_k}=b\langle\x,\one\rangle-h_{i,\bar j}^{\J}(\x).\]

Note that $\sum_{i}|J_i|=2^n+n-2$ so $\sum_{i:|J_i|>1}(|J_i|-1)=\sum_{i}(|J_i|-1)=n-1$ and thus $\sum_{i:|J_i|>1}|J_i|=n-1+\sum_{i:|J_i|>1}1$. Now by Proposition~\ref{prop:{lemma}affine-fine-mixed-cell-dual-volume}, we have
\begin{align*}
\tilde m_{\Q_{\J}}(\x)|_{\langle\x,\one\rangle=0}=&\kappa(\Q_{\J})\frac{\prod_{i:|J_i|>1}x_{T_i}}{\prod_{i:|J_i|>1}\prod_{a=1}^{|J_i|}h_{\x\Q_{\J}}(\v_{i,a})}\\
=&\frac{\prod_{i:|J_i|>1}x_{T_i}}{\prod_{i:|J_i|>1}\prod_{j\in J_i}(-h_{i,\bar j}^{\J}(\x))}=(-1)^{n-1}\prod_{i:|J_i|>1}\frac{-x_{T_i}}{\prod_{j\in J_i}h_{i,\bar j}^{\J}(\x)}.
\end{align*}
\end{proof}
One very nice property of the formula given in Proposition~\ref{prop:permutohedra-cell-dual-volume} is that every factor $(-x_{T_i})/\prod_{j\in J_i}h_{i,\bar j}^{\J}(\x)$ equals a subtraction-free rational function once we have substituted $x_{[n]}=-\sum_{T\subset[n]}x_T$. This is because if $T_i=[n]$, then the numerator $-x_{T_i}=-x_{[n]}$ becomes positive and none of the $h_{i,\bar j}^{\J}(\x)$ contains $x_{[n]}$; and if $T_i\neq[n]$, then the numerator is negative but there is exactly one $h_{i,\bar j}^{\J}(\x)$ in its denominator with $x_{[n]}$ in it and thus the negative signs cancel. 

\begin{example}\label{ex:permutohedra-cell-tree}
Let $n=3$ and order the nonempty subsets of $[n]$ as $(1,2,3,12,13,23,123)$. We can take $\J=(1,2,3,1,1,23,12)$ with the spanning tree $G_{\J}$ shown in Figure~\ref{fig:permutohedron-fine-cell-spanning-tree}.
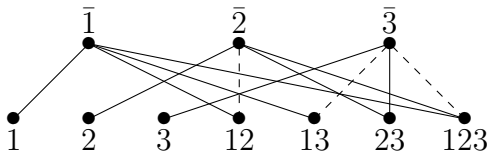
\begin{figure}[h!]
\centering
\begin{tikzpicture}[scale=1.0]
\node at (0,0) {$\bullet$};
\node at (1,0) {$\bullet$};
\node at (2,0) {$\bullet$};
\node at (3,0) {$\bullet$};
\node at (4,0) {$\bullet$};
\node at (5,0) {$\bullet$};
\node at (6,0) {$\bullet$};
\node[below] at (0,0) {$1$};
\node[below] at (1,0) {$2$};
\node[below] at (2,0) {$3$};
\node[below] at (3,0) {$12$};
\node[below] at (4,0) {$13$};
\node[below] at (5,0) {$23$};
\node[below] at (6,0) {$123$};
\node at (1,1) {$\bullet$};
\node at (3,1) {$\bullet$};
\node at (5,1) {$\bullet$};
\node[above] at (1,1) {$\bar1$};
\node[above] at (3,1) {$\bar2$};
\node[above] at (5,1) {$\bar3$};
\draw(0,0)--(1,1);
\draw(1,0)--(3,1);
\draw(2,0)--(5,1);
\draw(3,0)--(1,1);
\draw(4,0)--(1,1);
\draw(5,0)--(3,1);
\draw(5,0)--(5,1);
\draw(6,0)--(1,1);
\draw(6,0)--(3,1);
\draw[thin,dashed](3,0)--(3,1);
\draw[thin,dashed](4,0)--(5,1);
\draw[thin,dashed](6,0)--(5,1);
\end{tikzpicture}
\caption{A spanning tree $G_{\J}$}
\label{fig:permutohedron-fine-cell-spanning-tree}
\end{figure}

In this case, $J_{23}=\{2,3\}$ and $J_{123}=\{1,2\}$ are the only subsets with cardinality greater than $1$. If we remove $123$, then one connected component is $\{\bar1,1,12,13\}$ and the other connected component is $\{\bar2,\bar3,2,3,23\}$ so \[h_{123,\bar1}^{\J}(\x)=x_1+x_{12}+x_{13},\qquad h_{123,\bar2}^{\J}(\x)=x_2+x_3+x_{23}.\]
Similarly, removing $23$ also results in two connected components, and we have \[h_{23,\bar2}^{\J}(\x)=x_1+x_2+x_{12}+x_{13}+x_{123},\qquad h_{23,\bar3}^{\J}(\x)=x_3.\]
By Proposition~\ref{prop:permutohedra-cell-dual-volume}, we see that 
\begin{align*}
\tilde m_{\Q_\J}(\x)|_{\langle\x,\one\rangle=0}=&\frac{-x_{123}}{(x_1+x_{12}+x_{13})(x_2+x_3+x_{23})}\cdot\frac{-x_{23}}{(x_1+x_2+x_{12}+x_{13}+x_{123})(x_3)}\\
=&\frac{x_1+x_2+x_3+x_{12}+x_{13}+x_{23}}{(x_1+x_{12}+x_{13})(x_2+x_3+x_{23})}\cdot\frac{x_{23}}{(x_3+x_{23})(x_3)},
\end{align*}
with the substitution $x_{123}=-x_1-x_2-x_3-x_{12}-x_{13}-x_{23}$.
\end{example}

\begin{example}
When $n=3$, consider a fine mixed subdivision of $\P$: $\J^{(1)}=(1,2,3,1,1,23,12)$, $\J^{(2)}=(1,2,3,12,1,23,2)$, $\J^{(3)}=(1,2,3,2,13,23,2)$, $\J^{(4)}=(1,2,3,2,13,3,23)$, $\J^{(5)}=(1,2,3,12,13,3,3)$, $\J^{(6)}=(1,2,3,1,1,3,123)$, $\J^{(7)}=(1,2,3,12,1,3,23)$ shown in Figure~\ref{fig:fine-mixed-subdivision-generalized-permutohedron}. The cell $\J^{(1)}$ has been considered in Example~\ref{ex:permutohedra-cell-tree}.
\begin{figure}[h!]
\centering
\begin{tikzpicture}[scale=0.8]
\node at (0,1) {$\bullet$};
\node at (0,-1) {$\bullet$};
\node at (-1.732,0) {$\bullet$};
\node at (2,1) {$\bullet$};
\node at (2,-1) {$\bullet$};
\node at (-1,-2.732) {$\bullet$};
\node at (-2.732,-1.732) {$\bullet$};
\node at (-1,2.732) {$\bullet$};
\node at (-2.732,1.732) {$\bullet$};
\node[left] at (-1.732,0) {$e_3$};
\node[right] at (0,1) {$e_1$};
\node[right] at (0,-1) {$e_2$};
\node[right] at (2,-1) {$e_2$};
\node[right] at (2,1) {$e_1$};
\node[above] at (-2.732,-1.732) {$e_3$};
\node[right] at (-1,-2.732) {$e_2$};
\node[below] at (-2.732,1.732) {$e_3$};
\node[right] at (-1,2.732) {$e_1$};
\draw(0,1)--(0,-1)--(-1.732,0)--(0,1);
\draw(2,1)--(2,-1);
\draw(-1,-2.732)--(-2.732,-1.732);
\draw(-1,2.732)--(-2.732,1.732);
\node at (1,0) {$+$};
\node at (-1.366,-1.366) {$+$};
\node at (-1.366,1.366) {$+$};
\end{tikzpicture}
\qquad\qquad
\begin{tikzpicture}[scale=0.8]
\def\a{0.8};
\def\b{1.3};
\def\c{0.7};
\node at (0,0) {$\bullet$};
\node at (0,2) {$\bullet$};
\node at (-1.732,1) {$\bullet$};
\node at (\c*1.732,-\c) {$\bullet$};
\node at (\c*1.732,2-\c) {$\bullet$};
\node at (0,-2*\a) {$\bullet$};
\node at (\c*1.732,-\c-2*\a) {$\bullet$};
\node at (-\b*1.732,-2*\a-\b) {$\bullet$};
\node at (-\b*1.732+\c*1.732,-\c-2*\a-\b) {$\bullet$};
\node at (-1.732,-2*\a+1) {$\bullet$};
\node at (-\b*1.732-1.732,-2*\a-\b+1) {$\bullet$};
\node at (-\b*1.732-1.732,-\b+1) {$\bullet$};
\draw(0,2)--(\c*1.732,2-\c);
\draw(-1.732,1)--(\c*1.732,-\c);
\draw(-1.732,-2*\a+1)--(\c*1.732,-\c-2*\a);
\draw(-\b*1.732-1.732,-2*\a-\b+1)--(-\b*1.732+\c*1.732,-\c-2*\a-\b);
\draw(\c*1.732,2-\c)--(\c*1.732,-\c-2*\a);
\draw(0,2)--(0,-2*\a);
\draw(-1.732,1)--(-1.732,-2*\a+1);
\draw(-\b*1.732-1.732,-\b+1)--(-\b*1.732-1.732,-2*\a-\b+1);
\draw(-\b*1.732-1.732,-\b+1)--(0,2);
\draw(-\b*1.732-1.732,-2*\a-\b+1)--(-1.732,-2*\a+1);
\draw(-\b*1.732,-2*\a-\b)--(0,-2*\a);
\draw(-\b*1.732+\c*1.732,-\c-2*\a-\b)--(\c*1.732,-\c-2*\a);
\node at (-1/1.732,1) {$\J^{(6)}$};
\node at (\c*0.866,1-0.5*\c) {$\J^{(1)}$};
\node at (\c*0.866,-0.5*\c-\a) {$\J^{(2)}$};
\node at (\c*0.866-\b*0.866,-2*\a-0.5*\b-0.5*\c) {$\J^{(3)}$};
\node at (-\b*0.866-0.866,-2*\a-0.5*\b+0.5) {$\J^{(4)}$};
\node at (-\b*0.866-1.732,-\a-0.5*\b+1) {$\J^{(5)}$};
\node at (-0.866,-\a+0.5) {$\J^{(7)}$};
\end{tikzpicture}
\caption{A fine mixed subdivision of the generalized permutohedron}
\label{fig:fine-mixed-subdivision-generalized-permutohedron}
\end{figure}
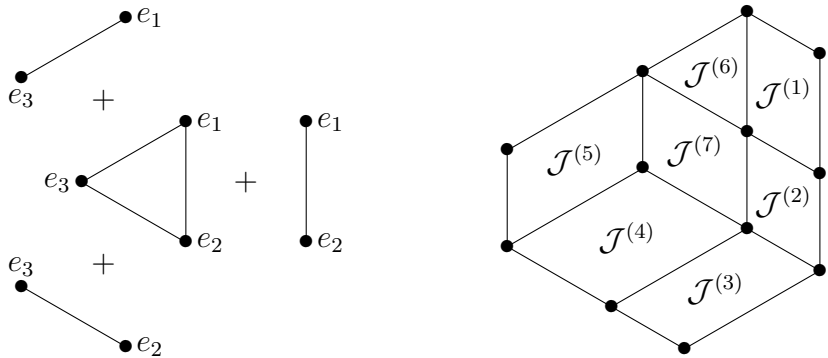

We then obtain the equality on $\tilde m_{\P}(\x)|_{\langle\x,\one\rangle=0}$, as a rational function in $\{x_T\}_{T\subset[n]\text{ is proper}}$:
\begin{align*}
& \frac{1}{x_1(x_1+x_2+x_{12})} +  \frac{1}{x_1(x_1+x_3+x_{13})} +  \frac{1}{x_2(x_1+x_2+x_{12})}  \\
&+  \frac{1}{x_2(x_2+x_3+x_{23})} +  \frac{1}{x_3(x_1+x_3+x_{13})} +  \frac{1}{x_3(x_2+x_3+x_{23})}\\
=&\frac{x_{23}}{(x_3{+}x_{23})(x_3)}\frac{x_1{+}x_2{+}x_3{+}x_{12}{+}x_{13}{+}x_{23}}{(x_1{+}x_{12}{+}x_{13})(x_{2}{+}x_3{+}x_{23})}+\frac{x_{12}}{(x_1{+}x_{13})(x_1{+}x_{12}{+}x_{13})}\frac{x_{23}}{(x_3{+}x_{23})(x_{3})}\\
&+\frac{x_{13}}{(x_1)(x_1{+}x_{13})}\frac{x_{23}}{(x_1{+}x_3{+}x_{13}{+}x_{23})(x_1{+}x_3{+}x_{13})}+\frac{x_{13}}{(x_1)(x_1{+}x_{13})}\frac{x_{1}{+}x_2{+}x_3{+}x_{12}{+}x_{13}{+}x_{23}}{(x_1{+}x_3{+}x_{13}{+}x_{23})(x_2{+}x_{12})}\\
&+\frac{x_{12}}{(x_2{+}x_{12})(x_2)}\frac{x_{13}}{(x_1{+}x_2{+}x_{12})(x_1{+}x_2{+}x_{12}{+}x_{13})}+\frac{x_{1}{+}x_2{+}x_3{+}x_{12}{+}x_{13}{+}x_{23}}{(x_1{+}x_{12}{+}x_{13})(x_2)(x_3{+}x_{23})}\\
&+\frac{x_{12}}{(x_1{+}x_{13})(x_1{+}x_{12}{+}x_{13})}\frac{x_{1}{+}x_2{+}x_3{+}x_{12}{+}x_{13}{+}x_{23}}{(x_1{+}x_2{+}x_{12}{+}x_{13})(x_3{+}x_{23})}.
\end{align*}
\end{example}

We assume that the nonempty subsets of $[n]$ have been ordered so that $T_N = [n]$.  Write $A^{(k)}_{i,\bar j}$ for the set $A_{i,\bar j}$ defined using $\J^{(k)}$.

\begin{cor}\label{cor:Jm1}
Let $(\J^{(1)}, \ldots, \J^{(m)})$ be a fine mixed subdivision of $\P$.  Then 
$$
\sum_{\sigma\in S_n}\prod_{a=1}^{n-1}\frac{1}{\displaystyle\sum_{\emptyset\neq T\subseteq\sigma[1:a]}x_T} = \left.\sum_{k=1}^m \prod_{i:|J^{(k)}_i|>1}\frac{-x_{T_i}}{\prod_{j\in J^{(k)}_i}h_{i,\bar j}^{\J^{(k)}}(\x)}\right|_{x_{[n]}= -\sum_{T\subset[n]} x_T}.
$$
\end{cor}

Let 
$$
\alpha_{i,\bar j}^{(k)}:= \begin{cases} |A_{i,\bar j}^{(k)}| & \mbox{if $N \notin A_{i,\bar j}^{(k)}$} \\
2^n - 1 - |A_{i,\bar j}^{(k)}| & \mbox{otherwise}.
\end{cases}
$$
Specializing $x_T = 1$ for $T \neq [n]$ gives
\begin{cor}\label{cor:Jm2}
Let $(\J^{(1)}, \ldots, \J^{(m)})$ be a fine mixed subdivision of $\P$.  Then
$$
\frac{n!}{(2^1 - 1)(2^2-1) \cdots(2^{n-1}-1)} = (2^n-2)\sum_{k \mid |\J^{(k)}_{N}|>1} \prod_{i \mid |J_i|>1}\prod_{j \in J_i}\frac{1}{|\alpha^{(k)}_{i,\bar j}|} + \sum_{k \mid |\J^{(k)}_{N}|<1} \prod_{i \mid |J_i|>1}\prod_{j \in J_i}\frac{1}{|\alpha^{(k)}_{i,\bar j}|}.
$$
\end{cor}

\section{Associahedra}\label{sec:associahedra}
In this section, we work with the setup in Section~\ref{sec:genperm}.
\begin{defin}\label{def:associahedra}
For $n\in\mathbb{Z}_{>0}$, let $\P=(\Delta_{[i,j]})_{1\leq i\leq j\leq n}$. The \emph{associahedron} is \[\x\P:=\sum_{1\leq i\leq j\leq n}x_{ij}\Delta_{[i,j]},\]
which is a special case of the generalized permutohedron (Definition~\ref{def:generalized-permuohedra}).
\end{defin}
Note that here we use $x_{ij}$, where $1\leq i\leq j\leq n$ to index a variable. Definition~\ref{def:associahedra} is commonly referred to as the \emph{Loday realization} \cite{loday} of the associahedra. Our presentation of the material largely follows \cite{postnikov-permutohedra}. 

First, we need to introduce some Catalan objects. 
\begin{defin}
A \emph{plane binary tree} is a tree such that each node has at most $1$ left child and at most $1$ right child. Note that if a node has only one child, we specify it to be either the left child or the right child. Let $\PB(n)$ denote the set of plane binary trees with $n$ nodes.
\end{defin}
\begin{defin}
A \emph{planar cubic tree} with $n$ leaves is a planar tree drawn in a disk with $n$ boundary vertices labeled $1$ through $n$ in counterclockwise order and $n-3$ interior vertices inside the disk such that each boundary vertex has degree $1$ and each interior vertex has degree $3$. Let $\PC(n)$ denote the set of planar cubic trees with $n$ leaves. 
\end{defin}
It is well-known that $|\PB(n)|=|\PC(n+2)|=\frac{1}{n+1}{2n+1\choose n}=C_n$, the $n$-th Catalan number. For any $B\in \PB(n)$, there is a unique way to label its vertices from $1$ to $n$ such that if a node is labeled $i$, then all of its descendants from its left branch are labeled with numbers less than $i$ and all of its descendants from its right branch are labeled with numbers greater than $i$. For $B\in\PB(n)$ and $i\in[n]$, all the nodes that are descendants of $i$, including itself, form an interval $[l_B(i),r_B(i)]$.
\begin{lemma}[Corollary 8.2 of \cite{postnikov-permutohedra}]
Let $x_{ij}>0$ for $1\leq i\leq j\leq n$. The set of vertices of $\x\P$ is in bijection with $\PB(n)$. Specifically, given $B\in \PB(n)$, denote its corresponding vertex as $\p_B$. Then its $k$-th coordinate equals \[\langle\p_B,\e_k\rangle=\sum_{l_B(k)\leq i\leq k\leq j\leq r_B(k)}x_{ij}.\]
\end{lemma}

\begin{ex}
Let $B\in\PB(8)$ as in Figure~\ref{fig:plane-binary-tree-example}. Then $[l_B(1),r_B(1)]=[1,3]$, $[l_B(2),r_B(2)]=\{2\}$, $[l_B(3),r_B(3)]=[2,3]$, $[l_B(4),r_B(4)]=[1,8]$, $[l_B(5),r_B(5)]=\{5\}$, $[l_B(6),r_B(6)]=[5,8]$, $[l_B(7),r_B(7)]=\{7\}$, $[l_B(8),r_B(8)]=[7,8]$. If $x_{ij}=1$ for all variables, then $B$ is corresponding to the vertex $\p_B=(3,1,2,20,1,6,1,2)$ of $\x\P$.
\begin{figure}[h!]
\centering
\begin{tikzpicture}[scale=0.5]
\node at (0,0) {$\bullet$};
\node at (-4,-1.5) {$\bullet$};
\node at (-2,-3) {$\bullet$};
\node at (-3,-4.5) {$\bullet$};
\node at (4,-1.5) {$\bullet$};
\node at (2,-3) {$\bullet$};
\node at (6,-3) {$\bullet$};
\node at (5,-4.5) {$\bullet$};
\node[below] at (0,0) {$4$};
\node[left] at (-4,-1.5) {$1$};
\node[right] at (-2,-3) {$3$};
\node[left] at (-3,-4.5) {$2$};
\node[right] at (4,-1.5) {$6$};
\node[below] at (2,-3) {$5$};
\node[right] at (5,-4.5) {$7$};
\node[right] at (6,-3) {$8$};
\draw(0,0)--(-4,-1.5)--(-2,-3)--(-3,-4.5);
\draw(0,0)--(4,-1.5)--(2,-3);
\draw(4,-1.5)--(6,-3)--(5,-4.5);
\node at (0,-5.5) {};
\end{tikzpicture}
\qquad
\begin{tikzpicture}[scale=0.7]
\draw (0,0) circle (2);
\coordinate (b10) at (90:2);
\coordinate (b1) at (120:2);
\coordinate (b2) at (150:2);
\coordinate (b3) at (180:2);
\coordinate (b4) at (210:2);
\coordinate (b5) at (240:2);
\coordinate (b6) at (270:2);
\coordinate (b7) at (300:2);
\coordinate (b8) at (330:2);
\coordinate (b9) at (0:2);
\coordinate (i4) at (0,1.3);
\coordinate (i1) at (-0.7,1.1);
\coordinate (i2) at (-1.5,0.2);
\coordinate (i3) at (-0.4,0.6);
\coordinate (i5) at (0,0);
\coordinate (i6) at (1.0,1.0);
\coordinate (i7) at (0.5,-0.5);
\coordinate (i8) at (1.4,0);
\node at (b1) {$\bullet$};
\node at (b2) {$\bullet$};
\node at (b3) {$\bullet$};
\node at (b4) {$\bullet$};
\node at (b5) {$\bullet$};
\node at (b6) {$\bullet$};
\node at (b7) {$\bullet$};
\node at (b8) {$\bullet$};
\node at (b9) {$\bullet$};
\node at (b10) {$\bullet$};
\node at (b1) {$\bullet$};
\node at (i1) {$\bullet$};
\node at (i2) {$\bullet$};
\node at (i3) {$\bullet$};
\node at (i4) {$\bullet$};
\node at (i5) {$\bullet$};
\node at (i6) {$\bullet$};
\node at (i7) {$\bullet$};
\node at (i8) {$\bullet$};
\node[above] at (b10) {$10$};
\node[above] at (b1) {$1$};
\node[left] at (b2) {$2$};
\node[left] at (b3) {$3$};
\node[left] at (b4) {$4$};
\node[below] at (b5) {$5$};
\node[below] at (b6) {$6$};
\node[below] at (b7) {$7$};
\node[right] at (b8) {$8$};
\node[right] at (b9) {$9$};
\draw(b10)--(i4)--(i1)--(i3)--(i2);
\draw(i4)--(i6)--(i8)--(i7)--(b7);
\draw(i6)--(i5)--(b5);
\draw(b1)--(i1);
\draw(b2)--(i2);
\draw(b3)--(i2);
\draw(b4)--(i3);
\draw(b6)--(i5);
\draw(b8)--(i7);
\draw(b9)--(i8);
\end{tikzpicture}
\caption{A plane binary tree $B\in\PB(8)$ and a planar cubic tree $\varphi(B)\in\PC(10)$}
\label{fig:plane-binary-tree-example}
\end{figure}
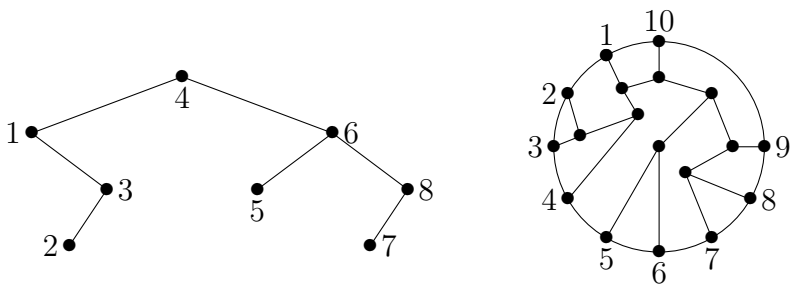
\end{ex}

We need some additional notations in order to describe the normal fan $\N(P)$ of the associahedron $\x\P$. Let $B$ be a plane binary tree and let $e\in B$ be an edge of it. Deleting $e$ from this tree results in two connected components $[n]=L_{B,e}\sqcup U_{B,e}$ where $U_{B,e}$ contains the root of $B$. Let $\v_{B,e}$ be the indicator vector of $L_{B,e}$. In other words, the $k$-th coordinate of $\v_{B,e}$ equals $1$ if $k\in L_{B,e}$, and $0$ if $k\in U_{B,e}$. Let $\tilde\v_{B,e}$ be the projection of $\v_{B,e}$ onto $H=\{\y\in\R^n\:|\:\langle\y,\one\rangle=0\}$.

The following result can be deduced from \cite[Proposition 2.6]{loday}.

\begin{prop}
    \label{prop:normal-fan-associahedron}
The maximal cones in the normal fan of the associahedron $\x\P$ are in bijection with $\PB(n)$. Specifically, given $B\in \PB(n)$, its corresponding maximal cone can be written as $C_{B}:=\mathrm{span}_{\R_{\geq 0}}\{\tilde\v_{B,e}\:|\: e\in B\}$.
\end{prop}

\begin{prop}\label{prop:assoc}
When $\langle\x,\one\rangle=0$, we have \[\tilde m_{\P}(\x)|_{\langle\x,\one\rangle=0}=(-1)^{n-1}\sum_{B\in\PB(n)}\prod_{e\in B}\frac{1}{\sum_{i,j\in L_{B,e}}x_{ij}}.\]
\end{prop}
\begin{proof}
In Definition~\ref{def:hyperplane-dual-mixed}, we notice that when $\z=\mathbf{0}$ and $\langle\x,\one\rangle=0$, $\tilde m_{\P}(x)$ stays unchanged when we increase $\v$ by a multiple of $\one$. For simplicity, we work with $\v_{B,e}$'s.

We first show that for any $B\in\PB(n)$, the vectors $\v_{B,e}$'s for all edges $e\in B$ together with $\one$ have determinant $\pm1$. Use induction on $n$, where the base case $n=1$ is clear. Let $k$ be the root of $B$ and label the edges as $e^{(1)},\ldots,e^{(k-1)},e^{(k+1)},\ldots,e^{(n)}$ where the edge $e^{(i)}$ is the edge whose lower endpoint is the vertex $i$. In particular, $\v_{e^{(k')}}=(1,1,\ldots,1,0,\ldots,0)$ with $k-1$ ones in the vector, where $k'$ is the left child of the root $k$. Replace $\one$ by $\one-\v_{e^{(k')}}$ the $n\times n$ matrix becomes block diagonal, which equals $\pm1$ by induction hypothesis. If $k$ does not have a left child, let $k'$ be its left child so that $\v_{e^{(k')}}=(0,1,\ldots,1)$. Replace $\one$ by $\one-\v_{e^{(k')}}$ and we obtain a block diagonal matrix where one block has size $1$. The determinant equals $\pm1$ as desired.

Next, we compute the support functions. The vectors of interest are all of the form $\v_{A}$, where $A\subseteq [n]$ and $\v_{A}$ is its indicator function. It is straightforward that $h_{\Delta_S}(\v_{A})$ equals $0$ if $S\cap ([n]\setminus A)\neq\emptyset$ and equals $-1$ if $S\subseteq A$. So $h_{\x\P}(\v_A)=-\sum_{[i,j]\subseteq A}x_{ij}$. Note that each $L_{B,e}$ is an interval, so $[i,j]\subseteq L_{B,e}$ is equivalent to $i,j\in L_{B,e}$. Putting together, by Definition~\ref{def:hyperplane-dual-mixed}, we conclude the desired proposition statement.
\end{proof}

We now translate various statistics to planar cubic trees, by defining a bijection $\varphi:\PB(n)\rightarrow\PC(n+2)$ as follows. Each plane binary tree has a depth-first search ordering on its nodes, by traversing from its root and always going to the left child first (see \cite{postnikov-permutohedra} for details). Under this order, whenever we encounter a node with degree $d$, we turn it into an interior vertex and connect it to $3-d$ boundary vertices to form a planar cubic tree. If the root of $B$ has two children, then in $\varphi(B)$ it is connected to the boundary vertex labeled $n+2$; if the root of $B$ is $1$ which has only a right child, then in $\varphi(B)$ it is connected to the boundary vertices $n+1$ and $n+2$; if the root of $B$ is $n$ which has only a left child, then it $\varphi(B)$ it is connected to the boundary vertices $n+2$ and $1$; the rest of labels in $\varphi(B)$ is extended cyclically. An example is seen in Figure~\ref{fig:plane-binary-tree-example}.

For a planar cubic tree $T\in\PC(n)$ and an interior edge $e\in T$, removing $e$ from $T$ result in a partition of the boundary vertices by connectivity, denoted as $L_{T,e}\sqcup U_{T,e}$ where $n\in U_{T,e}$. By planarity, $L_{T,e}$ is an interval and $U_{T,e}$ is a cyclic interval. Under the map $\varphi$, edges of $B\in\PB(n)$ naturally correspond to interior edges of $\varphi(B)\in\PC(n+2)$. The following lemma is a simple check that we leave for the readers.
\begin{lemma}
The map $\varphi$ is a bijection. Moreover, if $L_{B,e}=[i,j]$, then $L_{\varphi(B),e}=[i,j+1]$.
\end{lemma}
\begin{example}\label{ex:pentagon}
Let $n=3$ and consider the contribution of each maximal cones in $\N(P)$ towards the hyperplane dual mixed volume function $\tilde m_{\P}(\x)|_{\langle\x,\one\rangle=0}$:

\begin{center}
\begin{tabular}{c|c|c}
$\PB(3)$ & $\PC(5)$ & summand  \\ \hline
\centeredtab{\begin{tikzpicture}[scale=0.5]
\node at (0,0) {$\bullet$};
\node at (-1,-1) {$\bullet$};
\node at (1,-1) {$\bullet$};
\node[above] at (0,0) {$2$};
\node[below] at (-1,-1) {$1$};
\node[below] at (1,-1) {$3$};
\draw(-1,-1)--(0,0)--(1,-1);
\end{tikzpicture}} & 
\centeredtab{\begin{tikzpicture}[scale=0.6]
\draw (0,0) circle (1);
\coordinate (b1) at (162:1);
\coordinate (b2) at (234:1);
\coordinate (b3) at (306:1);
\coordinate (b4) at (18:1);
\coordinate (b5) at (90:1);
\coordinate (i1) at (0,0);
\coordinate (i2) at (198:0.6);
\coordinate (i3) at (-18:0.6);
\node at (b1) {$\bullet$};
\node at (b2) {$\bullet$};
\node at (b3) {$\bullet$};
\node at (b4) {$\bullet$};
\node at (b5) {$\bullet$};
\node at (i1) {$\bullet$};
\node at (i2) {$\bullet$};
\node at (i3) {$\bullet$};
\node[above] at (90:1) {$5$};
\node[left] at (162:1) {$1$};
\node[left] at (234:1) {$2$};
\node[right] at (306:1) {$3$};
\node[right] at (18:1) {$4$};
\draw(b5)--(i1);
\draw(b1)--(i2);
\draw(b2)--(i2);
\draw(b3)--(i3);
\draw(b4)--(i3);
\draw(i2)--(i1)--(i3);
\end{tikzpicture}} & \centeredtab{$\displaystyle{\frac{1}{x_{11}x_{33}}}$} \\\hline
\centeredtab{\begin{tikzpicture}[scale=0.5]
\node at (0,0) {$\bullet$};
\node at (-1,-1) {$\bullet$};
\node at (0,-2) {$\bullet$};
\draw(0,0)--(-1,-1)--(0,-2);
\node[right] at (0,0) {$3$};
\node[left] at (-1,-1) {$1$};
\node[right] at (0,-2) {$2$};
\end{tikzpicture}} & \centeredtab{\begin{tikzpicture}[scale=0.6]
\draw (0,0) circle (1);
\coordinate (b1) at (234:1);
\coordinate (b2) at (306:1);
\coordinate (b3) at (18:1);
\coordinate (b4) at (90:1);
\coordinate (b5) at (162:1);
\coordinate (i1) at (0,0);
\coordinate (i2) at (270:0.6);
\coordinate (i3) at (54:0.6);
\node at (b1) {$\bullet$};
\node at (b2) {$\bullet$};
\node at (b3) {$\bullet$};
\node at (b4) {$\bullet$};
\node at (b5) {$\bullet$};
\node at (i1) {$\bullet$};
\node at (i2) {$\bullet$};
\node at (i3) {$\bullet$};
\node[above] at (90:1) {$5$};
\node[left] at (162:1) {$1$};
\node[left] at (234:1) {$2$};
\node[right] at (306:1) {$3$};
\node[right] at (18:1) {$4$};
\draw(b5)--(i1);
\draw(b1)--(i2);
\draw(b2)--(i2);
\draw(b3)--(i3);
\draw(b4)--(i3);
\draw(i2)--(i1)--(i3);
\end{tikzpicture}} & \centeredtab{$\displaystyle{\frac{1}{x_{22}(x_{11}+x_{12}+x_{22})}}$}  \\\hline
\centeredtab{\begin{tikzpicture}[scale=0.5]
\node at (0,0) {$\bullet$};
\node at (1,-1) {$\bullet$};
\node at (2,-2) {$\bullet$};
\draw(0,0)--(2,-2);
\node[right] at (0,0) {$1$};
\node[right] at (1,-1) {$2$};
\node[right] at (2,-2) {$3$};
\end{tikzpicture}} & \centeredtab{\begin{tikzpicture}[scale=0.6]
\draw (0,0) circle (1);
\coordinate (b1) at (306:1);
\coordinate (b2) at (18:1);
\coordinate (b3) at (90:1);
\coordinate (b4) at (162:1);
\coordinate (b5) at (234:1);
\coordinate (i1) at (0,0);
\coordinate (i2) at (342:0.6);
\coordinate (i3) at (126:0.6);
\node at (b1) {$\bullet$};
\node at (b2) {$\bullet$};
\node at (b3) {$\bullet$};
\node at (b4) {$\bullet$};
\node at (b5) {$\bullet$};
\node at (i1) {$\bullet$};
\node at (i2) {$\bullet$};
\node at (i3) {$\bullet$};
\node[above] at (90:1) {$5$};
\node[left] at (162:1) {$1$};
\node[left] at (234:1) {$2$};
\node[right] at (306:1) {$3$};
\node[right] at (18:1) {$4$};
\draw(b5)--(i1);
\draw(b1)--(i2);
\draw(b2)--(i2);
\draw(b3)--(i3);
\draw(b4)--(i3);
\draw(i2)--(i1)--(i3);
\end{tikzpicture}} & \centeredtab{$\displaystyle{\frac{1}{x_{33}(x_{22}+x_{23}+x_{33})}}$}  \\\hline
\centeredtab{\begin{tikzpicture}[scale=0.5]
\node at (0,0) {$\bullet$};
\node at (-1,-1) {$\bullet$};
\node at (-2,-2) {$\bullet$};
\draw(0,0)--(-1,-1)--(-2,-2);
\node[right] at (0,0) {$3$};
\node[right] at (-1,-1) {$2$};
\node[right] at (-2,-2) {$1$};
\end{tikzpicture}} & \centeredtab{\begin{tikzpicture}[scale=0.6]
\draw (0,0) circle (1);
\coordinate (b1) at (18:1);
\coordinate (b2) at (90:1);
\coordinate (b3) at (162:1);
\coordinate (b4) at (234:1);
\coordinate (b5) at (306:1);
\coordinate (i1) at (0,0);
\coordinate (i2) at (54:0.6);
\coordinate (i3) at (198:0.6);
\node at (b1) {$\bullet$};
\node at (b2) {$\bullet$};
\node at (b3) {$\bullet$};
\node at (b4) {$\bullet$};
\node at (b5) {$\bullet$};
\node at (i1) {$\bullet$};
\node at (i2) {$\bullet$};
\node at (i3) {$\bullet$};
\node[above] at (90:1) {$5$};
\node[left] at (162:1) {$1$};
\node[left] at (234:1) {$2$};
\node[right] at (306:1) {$3$};
\node[right] at (18:1) {$4$};
\draw(b5)--(i1);
\draw(b1)--(i2);
\draw(b2)--(i2);
\draw(b3)--(i3);
\draw(b4)--(i3);
\draw(i2)--(i1)--(i3);
\end{tikzpicture}} & \centeredtab{$\displaystyle{\frac{1}{x_{11}(x_{11}+x_{12}+x_{22})}}$}  \\\hline
\centeredtab{\begin{tikzpicture}[scale=0.5]
\node at (0,0) {$\bullet$};
\node at (1,-1) {$\bullet$};
\node at (0,-2) {$\bullet$};
\draw(0,0)--(1,-1)--(0,-2);
\node[left] at (0,0) {$1$};
\node[right] at (1,-1) {$3$};
\node[left] at (0,-2) {$2$};
\end{tikzpicture}} & \centeredtab{\begin{tikzpicture}[scale=0.6]
\draw (0,0) circle (1);
\coordinate (b1) at (90:1);
\coordinate (b2) at (162:1);
\coordinate (b3) at (234:1);
\coordinate (b4) at (306:1);
\coordinate (b5) at (18:1);
\coordinate (i1) at (0,0);
\coordinate (i2) at (136:0.6);
\coordinate (i3) at (270:0.6);
\node at (b1) {$\bullet$};
\node at (b2) {$\bullet$};
\node at (b3) {$\bullet$};
\node at (b4) {$\bullet$};
\node at (b5) {$\bullet$};
\node at (i1) {$\bullet$};
\node at (i2) {$\bullet$};
\node at (i3) {$\bullet$};
\node[above] at (90:1) {$5$};
\node[left] at (162:1) {$1$};
\node[left] at (234:1) {$2$};
\node[right] at (306:1) {$3$};
\node[right] at (18:1) {$4$};
\draw(b5)--(i1);
\draw(b1)--(i2);
\draw(b2)--(i2);
\draw(b3)--(i3);
\draw(b4)--(i3);
\draw(i2)--(i1)--(i3);
\end{tikzpicture}} & \centeredtab{$\displaystyle{\frac{1}{x_{22}(x_{22}+x_{23}+x_{33})}}$}
\end{tabular}
\end{center}
\end{example}

\subsection{Relation to planar $\phi^3$-amplitude at tree level}
We relate \cref{prop:assoc} to a rational function $\mathcal{A}^{\phi_3}_n(\mathbf{s}_{ij})$ appearing in physics, called the $\phi^3$-amplitude; see \cite{AHLstringy,ABHY}.  This rational function was a main motivation for our study of dual mixed volumes.

Scattering amplitudes are functions that compute the outcome of scattering experiments in particle physics.  Traditionally, they are computed as the sum over Feynman diagrams which depend on the choice of particles and their interactions.  In ``planar $\phi^3$-theory", the Feynman diagrams are planar cubic trees.  For $n$-particle scattering, the amplitude\footnote{The full amplitude has a pertubative expansion; we only consider the first term which is a sum over trees.  The later terms involve graphs with cycles.} $\mathcal{A}^{\phi_3}_n(\mathbf{s}_{ij})$ is the following sum over $\PC(n)$:
\begin{equation}\label{eq:phi3}
\mathcal{A}^{\phi_3}_n := \sum_{T \in \PC(n)} \prod_{e \in E(T)} \frac{1}{X_e}
\end{equation}
where the product is over the interior edge set $E(T)$ of $T$, and 
$$
X_e = \sum_{a \leq i < j \leq b-1} s_{ij}
$$
where one of the connected components of $T \setminus e$ contains the leaves labeled $a,a+1,\ldots,b-1$, considered cyclically. 
 Here, the \emph{Mandelstam variables} $s_{ij}$ satisfy the relations
\begin{equation}\label{eq:Man}
s_{ii} = 0, \qquad s_{ij}= s_{ji}, \qquad \sum_{j} s_{ij} = 0 \mbox{ for each }i. 
\end{equation}
The quantity $X_e$ has the physical interpretation as the square of the momentum travelling along the edge $e$. Comparing \eqref{eq:phi3} with \cref{prop:assoc}, we obtain the following result.

\begin{prop}
Up to a sign, the dual mixed volume $\tilde m_{\P}(\x)|_{\langle\x,\one\rangle=0}$ of the associahedron is equal to the planar $\phi^3$-amplitude $\mathcal{A}^{\phi_3}_n$ under the substitution
$x_{ij} \mapsto s_{i,j+1}$.
\end{prop}

It would be interesting to investigate the implications of \cref{thm: mixed subdivision formula} for $A_n^{\phi^3}$.

\begin{example}
    Continuing \cref{ex:pentagon}, the substitution $x_{ij} \mapsto s_{i,j+1}$ gives the five terms
$$
\frac{1}{s_{12}s_{34}}, \frac{1}{s_{23}(s_{12}+s_{13}+s_{23})}, \frac{1}{s_{34}(s_{23}+s_{24}+s_{34})}, \frac{1}{s_{12}(s_{12}+s_{13}+s_{23})}, \frac{1}{s_{23}(s_{23}+s_{24}+s_{34})}
$$
Using the identity
$$
\sum_{a \leq i < j \leq b-1} s_{ij} = \sum_{b \leq i < j \leq a-1} s_{ij}
$$
that can be deduced from \eqref{eq:Man}, we can rewrite the last four terms to get
$$
A_5^{\phi^3} = \frac{1}{s_{12}s_{34}}+ \frac{1}{s_{23}s_{45}}+ \frac{1}{s_{34}s_{15}}+ \frac{1}{s_{12}s_{45}}+ \frac{1}{s_{23}s_{15}}.
$$
\end{example}

\begin{remark}
    The generalized associahedra appearing in cluster algebras have Minkowski decompositions analogous to the associahedron; see for example \cite{AHLcluster}.  Thus generalized associahedra also possess natural dual mixed volume functions.
\end{remark}

\section{Further directions}\label{sec:further}
Let $A_\P(\x)$ denote the numerator of the dual mixed volume function, as in \eqref{eq:DMVratio}.  It would be interesting to investigate what polynomials can appear as $A_\P(\x)$ and compare to the mixed volume polynomial $\Vol_\P(\x)$ defined in \eqref{eq:MV}.

The degree of the homogeneous polynomial $A_\P(\x)$ is equal to $f-d$ where $f$ is the number of facets of $P=P_1+\cdots+P_r$.

\begin{question}
What can we say about the support of the polynomial $A_\P(\x)$?
\end{question}

We expect the coefficients of $A_\P(\x)$ to satisfy interesting positivity properties.  For example, the numerators in \cref{sec:genperm} and \cref{sec:associahedra} are (up to an overall sign) positive polynomials.  Aluffi \cite{Aluffi} showed that certain adjoint polynomials (see \cref{sec:adjoint}) are \emph{sectionally log-concave}, relating them to covolume polynomials and Lorentzian polynomials.  We expect that these results can be generalized to dual mixed volumes.  We view such results as dual analogues of the Alexandrov-Fenchel inequality \eqref{eq:AF} for mixed volumes.

\begin{problem}
    Find conditions on $\P=(P_1,\ldots,P_r)$ to guarantee that $A_\P(\x)$ is a polynomial with positive coefficients or satisfies sectional log-concavity.
\end{problem}

The following is a dual (but easier) version of the Brunn-Minkowski inequality first established by Firey in \cite{Firey}, also see \cite[p.~207]{Barvinokconvexity}), which in its easiest form states that
$$
\Vol((A+B)/2) \geq \Vol(A)^{1/2} \Vol(B)^{1/2}
$$  
for two non-empty compact sets $A,B$.   

\begin{prop}\label{prop:dualBM}\cite{Firey}
Let $P$ and $Q$ be polytopes containing $\0$ in the interior. We have
$$
\VolC((P+Q)/2)^2 \leq \VolC(P)\VolC(Q).
$$
\end{prop}
\begin{proof}
By \cref{cor:intform}, we have
$$
\VolC(P)=\int_{\R^d}e^{-h_P(\v)}d\v.
$$  
H\"older's inequality states that
$$
\int_{\R^d} |f(\v) g(\v)| d\v \leq \left(\int_{\R^d} |f(\v)|^2 d\v \right)^{1/2} \left(\int_{\R^d} |g(\v)|^2 d\v\right)^{1/2}.
$$
Thus
    $$
    \VolC((P+Q)/2) = \int_{\R^n}|e^{-h_P(\v)/2-h_Q(\v)/2}|d\v \leq  \left(\int_{\R^n}|e^{-h_P(\v)}|d\v\right)^{1/2} \left(\int_{\R^n}|e^{-h_Q(\v)}|d\v\right)^{1/2}
    $$
giving the stated inequality.
\end{proof}
This result generalizes to closed convex sets.

It would be interesting to understand the relation between \cref{prop:dualBM} and the positivity (and other) properties of the coefficients of $A_\P(\x)$.


\bibliographystyle{alpha}
\bibliography{note}

\newcommand{\etalchar}[1]{$^{#1}$}
\begin{thebibliography}{AHBHY18}

\bibitem[AHBHY18]{ABHY}
Nima Arkani-Hamed, Yuntao Bai, Song He, and Gongwang Yan.
\newblock {Scattering Forms and the Positive Geometry of Kinematics, Color and
  the Worldsheet}.
\newblock {\em JHEP}, 05:096, 2018.

\bibitem[AHBL17]{ABL}
Nima Arkani-Hamed, Yuntao Bai, and Thomas Lam.
\newblock Positive geometries and canonical forms.
\newblock {\em J. High Energy Phys.}, (11):039, front matter+121, 2017.

\bibitem[AHHL21a]{AHLcluster}
Nima Arkani-Hamed, Song He, and Thomas Lam.
\newblock Cluster configuration spaces of finite type.
\newblock {\em SIGMA Symmetry Integrability Geom. Methods Appl.}, 17:Paper No.
  092, 41, 2021.

\bibitem[AHHL21b]{AHLstringy}
Nima Arkani-Hamed, Song He, and Thomas Lam.
\newblock Stringy canonical forms.
\newblock {\em J. High Energy Phys.}, (2):Paper No. 069, 59, 2021.

\bibitem[Alu24]{Aluffi}
Paolo Aluffi.
\newblock Lorentzian polynomials, {S}egre classes, and adjoint polynomials of
  convex polyhedral cones.
\newblock {\em Adv. Math.}, 437:Paper No. 109440, 37, 2024.

\bibitem[Bar]{Barvinokprivatecommunication}
Alexander Barvinok.
\newblock The dual volume valuation on polyhedra.

\bibitem[Bar02]{Barvinokconvexity}
Alexander Barvinok.
\newblock {\em A course in convexity}, volume~54 of {\em Graduate Studies in
  Mathematics}.
\newblock American Mathematical Society, Providence, RI, 2002.

\bibitem[BLVS{\etalchar{+}}99]{oriented-matroid-book}
Anders Bj\"orner, Michel Las~Vergnas, Bernd Sturmfels, Neil White, and
  G\"unter~M. Ziegler.
\newblock {\em Oriented matroids}, volume~46 of {\em Encyclopedia of
  Mathematics and its Applications}.
\newblock Cambridge University Press, Cambridge, second edition, 1999.

\bibitem[Fil92]{filliman}
P.~Filliman.
\newblock The volume of duals and sections of polytopes.
\newblock {\em Mathematika}, 39(1):67--80, 1992.

\bibitem[Fir61]{Firey}
William~J. Firey.
\newblock Polar means of convex bodies and a dual to the {B}runn-{M}inkowski
  theorem.
\newblock {\em Canadian J. Math.}, 13:444--453, 1961.

\bibitem[Gae]{Gaetz}
Christian Gaetz.
\newblock Positive geometries learning seminar. lecture 3: canonical forms of
  polytopes from adjoints.

\bibitem[HRS00]{HRS}
Birkett Huber, J\"{o}rg Rambau, and Francisco Santos.
\newblock The {C}ayley trick, lifting subdivisions and the {B}ohne-{D}ress
  theorem on zonotopal tilings.
\newblock {\em J. Eur. Math. Soc. (JEMS)}, 2(2):179--198, 2000.

\bibitem[KR20]{Kohn-Ranestad}
Kathl\'en Kohn and Kristian Ranestad.
\newblock Projective geometry of {W}achspress coordinates.
\newblock {\em Found. Comput. Math.}, 20(5):1135--1173, 2020.

\bibitem[Kup03]{Kuperburg}
Greg Kuperberg.
\newblock A generalization of {F}illiman duality.
\newblock {\em Proc. Amer. Math. Soc.}, 131(12):3893--3899, 2003.

\bibitem[Lam24]{lam2022invitation}
Thomas Lam.
\newblock An invitation to positive geometries.
\newblock In {\em Open Problems in Algebraic Combinatorics}, volume 110 of {\em
  Proceedings of Symposia in Pure Mathematics}, pages 159--180. Amer. Math.
  Soc., Providence, RI, 2024.

\bibitem[Lod04]{loday}
Jean-Louis Loday.
\newblock Realization of the {S}tasheff polytope.
\newblock {\em Arch. Math. (Basel)}, 83(3):267--278, 2004.

\bibitem[McM89]{mcmullen-polytope-algebra}
Peter McMullen.
\newblock The polytope algebra.
\newblock {\em Adv. Math.}, 78(1):76--130, 1989.

\bibitem[Pos09]{postnikov-permutohedra}
Alexander Postnikov.
\newblock Permutohedra, associahedra, and beyond.
\newblock {\em Int. Math. Res. Not. IMRN}, (6):1026--1106, 2009.

\bibitem[RGZ94]{bohne-dress-richter-gebert-ziegler}
J\"urgen Richter-Gebert and G\"unter~M. Ziegler.
\newblock Zonotopal tilings and the {B}ohne-{D}ress theorem.
\newblock In {\em Jerusalem combinatorics '93}, volume 178 of {\em Contemp.
  Math.}, pages 211--232. Amer. Math. Soc., Providence, RI, 1994.

\bibitem[San05]{Santos}
Francisco Santos.
\newblock The {C}ayley trick and triangulations of products of simplices.
\newblock In {\em Integer points in polyhedra---geometry, number theory,
  algebra, optimization}, volume 374 of {\em Contemp. Math.}, pages 151--177.
  Amer. Math. Soc., Providence, RI, 2005.

\bibitem[Sch13]{Schneider}
Rolf Schneider.
\newblock {\em Convex Bodies: The Brunn–Minkowski Theory}.
\newblock Encyclopedia of Mathematics and its Applications. Cambridge
  University Press, 2 edition, 2013.

\bibitem[Sta81]{stanley-combinatorial-Aleksandrov-Fenchel}
Richard~P. Stanley.
\newblock Two combinatorial applications of the {A}leksandrov-{F}enchel
  inequalities.
\newblock {\em J. Combin. Theory Ser. A}, 31(1):56--65, 1981.

\bibitem[War96]{Warren}
Joe Warren.
\newblock Barycentric coordinates for convex polytopes.
\newblock {\em Adv. Comput. Math.}, 6(2):97--108, 1996.

\end{thebibliography}
\end{document}